\documentclass[12pt]{article}
\usepackage{amsbsy,amsfonts,amsmath,amssymb,amsthm,bbm,euscript,eepic,rotating}
\usepackage{graphicx}

\def\UseSection{
        \numberwithin{equation}{section}
    \theoremstyle{plain}
        \newtheorem{theorem}    {Theorem}[section]
        \DefineTheorems 
}

\def\DefineTheorems{
    
    \newtheorem{lemma}      [theorem] {Lemma}
    
    \newtheorem{prop}       [theorem] {Proposition}
    
    \newtheorem{cor}        [theorem] {Corollary}

    \theoremstyle{definition}
    \newtheorem{defn}       [theorem] {Definition}

    \theoremstyle{definition}

}

\newcommand{\bt}   {\begin{theorem}}
\newcommand{\et}   {\end  {theorem}}
\newcommand{\bl}   {\begin{lemma}}
\newcommand{\el}   {\end  {lemma}}
\newcommand{\bp}   {\begin{prop}}
\newcommand{\ep}   {\end  {prop}}
\newcommand{\bc}   {\begin{cor}}
\newcommand{\ec}   {\end  {cor}}
\newcommand{\bd}   {\begin{defn}}
\newcommand{\ed}   {\end  {defn}}

\newcommand{\ba}   {\begin{array}}
\newcommand{\ea}   {\end  {array}}
\newcommand{\be}   {\begin{enumerate}}
\newcommand{\ee}   {\end  {enumerate}}
\newcommand{\bi}   {\begin{itemize}}
\newcommand{\ei}   {\end  {itemize}}

\def\eq#1\en{\begin{equation}#1\end{equation}}
\def\eqsplit#1\ensplit{
    \begin{equation}\begin{split}#1\end{split}\end{equation}
    }
\def\eqalign#1\enalign{
    \begin{align}#1\end{align}
    }
\def\eqmul#1\enmul{
    \begin{multline}#1\end{multline}
    }
\newcommand{\eqarrstar} {\begin{eqnarray*}}
\newcommand{\enarrstar} {\end{eqnarray*}}
\newcommand{\eqarray}   {\begin{eqnarray}}
\newcommand{\enarray}   {\end{eqnarray}}

\newcommand{\lbeq}[1]  {\label{e:#1}}
\newcommand{\refeq}[1] {\eqref{e:#1}}    

\makeatletter
\newcommand{\labelcounter}[2]{{%
    \stepcounter{#1}
    \protected@write\@auxout{}%
    {\string\newlabel{#2}{{\csname the#1\endcsname}{\thepage}}}%
    {\ref{#2}}
    }}
\makeatother
\newcommand{\sss}   { \scriptscriptstyle }

\newcommand{\Cbold} {{\mathbb C}}

\newcommand{\Nbold} {{\mathbb N}}

\newcommand{\Rbold} {{\mathbb R}}

\newcommand{\Zbold} {{\mathbb Z}}

\newcommand{\avec}  {\boldsymbol{a}}

\newcommand{\ovec}  {\boldsymbol{o}}

\newcommand{\svec}  {\boldsymbol{s}}

\newcommand{\uvec}  {\boldsymbol{u}}
\newcommand{\vvec}  {\boldsymbol{v}}
\newcommand{\wvec}  {\boldsymbol{w}}
\newcommand{\xvec}  {\boldsymbol{x}}
\newcommand{\yvec}  {\boldsymbol{y}}
\newcommand{\zvec}  {\boldsymbol{z}}

\newcommand{\Acal}   {\mathcal{A}}

\newcommand{\Ccal}   {\mathcal{C}}

\newcommand{\Rd}    {{ {\Rbold}^d}}
\newcommand{\Zd}    {{ {\Zbold}^d }}

\newcommand{\spose}[1] {{\hbox to 0pt{#1\hss}} }
\newcommand{\ltapprox} {\mathrel{\spose{\lower 3pt\hbox{$\mathchar"218$}}
 \raise 2.0pt\hbox{$\mathchar"13C$}}}
\newcommand{\gtapprox} {\mathrel{\spose{\lower 3pt\hbox{$\mathchar"218$}}
 \raise 2.0pt\hbox{$\mathchar"13E$}}}

\UseSection  
\theoremstyle{plain}
\newtheorem{thm}{Theorem}[section]
\newtheorem{lem}[thm]{Lemma}

\newcommand{\bb}{\underline{b}}
\newcommand{\bA}{{\bf A}}
\newcommand{\bC}{{\bf C}}
\newcommand{\cA}{{\cal A}}

\newcommand{\cC}{{\cal C}}
\newcommand{\cE}{{\cal E}}

\newcommand{\Cg}{C\!_{\sss K}}

\newcommand{\cL}{{\cal L}}

\newcommand{\cR}{{\cal R}}
\newcommand{\daw}{\downarrow}
\newcommand{\db}{\Longrightarrow}

\newcommand{\BDcup}[1]{~\underset{#1}{\Dot{\bigcup}}~}
\newcommand{\Dcup}[2]{~\Dot{\cup}_{#1}^{#2}~}
\newcommand{\DDcup}{~\Dot{\cup}}

\newcommand{\ddsum}{\sideset{_{}^{}}{_{}^{\bullet}}\sum}
\newcommand{\dsum}{\sum^\bullet}

\newcommand{\dpst}{\displaystyle}

\newcommand{\ind}[1]{\mathbbm{1}{\scriptstyle #1}}
\newcommand{\indic}[1]{\ind{\{ #1 \}}}
\newcommand{\lamb}{\lambda}
\newcommand{\lambc}{\lamb_{\rm c}}
\newcommand{\lambce}{\lambc^{\sss(\vep)}}

\newcommand{\mN}{{\mathbb N}}
\newcommand{\mP}{{\mathbb P}}

\newcommand{\mR}{{\mathbb R}}
\newcommand{\mZ}{{\mathbb Z}}

\newcommand{\nn}{\nonumber}
\newcommand{\nnmb}{\nonumber}

\newcommand{\tb}{\overline{b}}

\newcommand{\uaw}{\uparrow}

\newcommand{\vep}{\varepsilon}

\newcommand{\wD}{\hat{D}}

\newcommand{\wpr}{\hat{p}_\vep}
\newcommand{\wtau}{\hat{\tau}}
\newcommand{\wpi}{\hat{\pi}}
\newcommand{\Zp}{\mZ_+}

 \newcommand{\smallsup}[1] {{\scriptscriptstyle{({#1}})}}
 
 \newcommand{\R}{\Rbold}
\newcommand{\Z}{\Zbold}
\newcommand{\N}{\Nbold}
\newcommand{\conn}{\longrightarrow}

\newcommand{\dbc}{\Longrightarrow}
\newcommand{\ct}[1]     { \stackrel{#1}{\conn} }

\newcommand{\ctx}[1]     { \xrightarrow{#1} }

\newcommand{\shift}   {\!\!\!\!}

\newcommand{\lupa}{\begin{rotate}{-10}$\scriptstyle\Longleftarrow$
\end{rotate}}
\newcommand{\rupa}{\begin{rotate}{10}$\scriptstyle\Longrightarrow$
\end{rotate}}
\newcommand{\btr}{$\scriptstyle\blacktriangle$}
\newcommand{\wtr}{$\scriptstyle\vartriangle$}

\newcommand{\sT}{{\sss T}}
\newcommand{\sN}{{\sss N}}
\newcommand{\sR}{{\sss R}}
\newcommand{\sstar}{{\scriptstyle\,\star\,}}

\title{Gaussian scaling for the critical spread-out contact process\\
   above the upper critical dimension}

\author{
    Remco van der Hofstad\thanks{Department of Mathematics and
        Computer Science, Eindhoven University of Technology,
        P.O.\ Box 513, 5600 MB Eindhoven, The Netherlands.
        {\tt r.w.v.d.hofstad@TUE.nl}}\\
    Akira Sakai\thanks{EURANDOM, P.O.\ Box 513, 5600 MB Eindhoven,
        The Netherlands.  {\tt sakai@eurandom.tue.nl}}
    }

\date{August 7, 2003\footnote{Revised September 24, 2004}}

\begin{document}
\maketitle

\begin{abstract}
We consider the critical spread-out contact process in $\Zd$ with $d\geq1$,
whose infection range is denoted by $L\geq1$.  The two-point function
$\tau_t(x)$ is the probability that $x\in\Zd$ is infected at time $t$
by the infected individual located at the origin $o\in\Zd$ at time 0.
We prove Gaussian behaviour for the two-point function with $L\geq L_0$ for
some finite $L_0=L_0(d)$ for $d>4$.  When $d\leq4$, we also perform a local
mean-field limit to obtain Gaussian behaviour for $\tau_{\sss tT}(x)$ with
$t>0$ fixed and $T\to\infty$ when the infection range depends on $T$ in
such a way that $L_{\sT}=LT^b$ for any $b>(4-d)/2d$.

The proof is based on the lace expansion and an adaptation of the inductive
approach applied to the discretized contact process.  We prove the existence
of several critical exponents and show that they take on their respective
mean-field values.  The results in this paper provide crucial ingredients
to prove convergence of the finite-dimensional distributions for the contact
process towards those for the canonical measure of super-Brownian motion,
which we defer to a sequel of this paper.

The results in this paper also apply to oriented percolation, for which we
reprove some of the results in \cite{hs01} and extend the results to the
local mean-field setting described above when $d\leq4$.
\end{abstract}

\section{Introduction and results}\label{s:intro}
\subsection{Introduction}\label{ss:intro}
The contact process is a model for the spread of an infection among
individuals in the $d$-dimensional integer lattice $\mZ^d$.  We suppose
that the origin $o\in\Zd$ is the only infected individual at time 0, and
that every infected individual may infect a healthy individual at a distance
less than $L\geq1$.  We refer to this model as the {\it spread-out} contact
process.  The rate of infection is denoted by $\lamb$, and it is well known
that there is a phase transition in $\lamb$ (see e.g., \cite{Ligg99}).

Sakai \cite{s00,s01} has proved that when $d>4$, the sufficiently spread-out
contact process has several critical exponents which are equal to those of
branching random walk.  The proof by Sakai uses the lace expansion for the
time-discretized contact process, and the main ingredient is the proof of
the so-called infrared bound uniformly in the time discretization.  Thus,
we can think of his results as proving Gaussian {\it upper bounds} for the
two-point function of the critical contact process.  Since these Gaussian
upper bounds imply the so-called triangle condition in \cite{bw98}, it
follows that certain critical exponents take on their mean-field values,
i.e., the values for branching random walk.  These values also agree with
the critical exponents appearing on the tree.  See
\cite[Chapter~I.4]{Ligg99} for an extensive account of the contact process
on a tree.

Recently, van~der~Hofstad and Slade \cite{hs01} proved that for all
$r\geq2$, the $r$-point functions for sufficiently spread-out critical
oriented percolation with spatial dimension $d>4$ converge to those of the
canonical measure of super-Brownian motion when we scale space by $n^{1/2}$,
where $n$ is the largest temporal component among the $r$ points, and then
take $n\uaw\infty$.  That is, the finite-dimensional distributions of the
critical oriented percolation cluster when it survives up to time $n$
converge to those of the canonical measure of super-Brownian motion.
The proof in \cite{hs01} is based on the lace expansion and the inductive
method of \cite{hs02}.  Important ingredients in \cite{hs01} are detailed
asymptotics and estimates of the oriented percolation two-point function.
The proof for the higher-point functions then follows by deriving a lace
expansion for the $r$-point functions together with an induction argument
in $r$.

In this paper, we prove the two-point function results for the contact
process via a time discretization.  The discretized contact process is
oriented percolation in $\Zd\times\vep\Zp$ with $\vep\in(0,1]$, and the
proof uses the same strategy as applied to oriented percolation with
$\vep=1$, i.e., an application of the lace expansion and the inductive
method.  However, to obtain the results for $\vep\ll1$, we use a
different lace expansion from the two expansions used in
\cite[Sections~3.1--3.2]{hs01}, and modify the induction hypotheses of
\cite{hs02} to incorporate the $\vep$-dependence.  In order to extend
the results from infrared bounds (as in \cite{s01}) to precise asymptotics (as in \cite{hs01}), it is
imperative to prove that the properly scaled lace expansion coefficients
converge to a certain continuum limit.  We can think of this continuum
limit as giving rise to a lace expansion in continuous time, even though
our proof is not based on the arising partial differential equation.
In the proof that the continuum limit exists, we make heavy use of
convergence results in \cite{bg91} which show that the discretized contact
process converges to the original continuous-time contact process.


In a sequel to this paper \cite{hsa04}, we use the results proved here as a
key ingredient in the proof that the finite-dimensional distributions of
the critical contact process above four dimensions converge to those of the
canonical measure of super-Brownian motion, as was proved in \cite{hs01}
for oriented percolation.

\subsection{The spread-out contact process and main results}\label{ss:results}
We define the spread-out contact process as follows.  Let $\bC_t\subset\Zd$
be the set of infected individuals at time $t\in\mR_+$, and let
$\bC_0=\{o\}$.  An infected site $x$ recovers in a small time interval
$[t,t+\vep]$ with probability $\vep+o(\vep)$ independently of $t$, where
$o(\vep)$ is a function that satisfies $\lim_{\vep\to0}o(\vep)/\vep=0$.
In other words, $x\in\bC_t$ recovers at rate 1.  A healthy site $x$ gets
infected, depending on the status of its neighbours, at rate
$\lamb\sum_{y\in\bC_t}D(x-y)$, where $\lamb\geq0$ is the infection rate and
$D(x-y)$ represents the strength of the interaction between $x$ and $y$.
We denote by $\mP^\lamb$ the associated probability measure.

The function $D$ is a probability distribution over $\Zd$ that is symmetric
with respect to the lattice symmetries, and satisfies certain assumptions
that involve a parameter $L\geq 1$ which serves to spread out the
infections and will be taken to be large.  In particular, we require that
there are $L$-independent constants $C,C_1,C_2\in(0,\infty)$ such that
$D(o)=0$, $\sup_{x\in\Zd}D(x)\leq CL^{-d}$ and $C_1L\leq\sigma\leq C_2L$,
where $\sigma^2$ is the variance of $D$:
\begin{align}\lbeq{sgm-def}
\sigma^2=\sum_{x\in\Zd}|x|^2D(x),
\end{align}
where $|\cdot|$ denotes the Euclidean norm on $\R^d$.  Moreover, we require
that there is a $\Delta>0$ such that
\begin{align}\lbeq{Deltadef}
\sum_{x\in\Zd}|x|^{2+2\Delta}D(x)\leq CL^{2+2\Delta}.
\end{align}
See Section~\ref{sss:genass} for the precise assumptions on $D$.  A simple
example of $D$ is the uniform distribution over the cube of side length
$2L$, excluding its center:
\begin{align}\lbeq{Dexample}
D(x)=\frac{\indic{0<\|x\|_\infty\leq L}}{(2L+1)^d-1},
\end{align}
where $\|x\|_\infty=\sup_i|x_i|$ for $x=(x_1,\dots,x_d)$.

The {\it two-point function} is defined as
\begin{align}\lbeq{2pt-def}
\tau_t^\lamb(x)=\mP^\lamb(x\in\bC_t)\qquad(x\in\Zd,~t\in\mR_+).
\end{align}
In words, $\tau_t^\lamb(x)$ is the probability that at time $t$, the
individual located at $x\in\Zd$ is infected due to the infection located
at $o\in\Zd$ at time 0.

By an extension of the results in \cite{bg91, gh01} to the spread-out
contact process, there exists a unique critical value $\lambc\in(0,\infty)$
such that
\begin{align}\lbeq{lambc-def}
\chi(\lambda)=\int_0^\infty\!\!dt~\hat\tau_t^\lamb(0)\begin{cases}
 <\infty,&\text{if }\lamb<\lambc,\\ =\infty,&\text{if }\lamb\geq\lambc,
 \end{cases}&&
\theta(\lambda)\equiv\lim_{t\uaw\infty}\mP^\lamb(\bC_t\ne\varnothing)
 \begin{cases}
 =0,&\text{if }\lamb\leq\lambc,\\ >0,&\text{if }\lamb>\lambc,
 \end{cases}
\end{align}
where we denote the Fourier transform of a summable function
$f:\Zd\mapsto\mR$ by
\begin{align}
\hat f(k)=\sum_{x\in\Zd}f(x)\;e^{ik\cdot x}\qquad(k\in[-\pi,\pi]^d).
\end{align}
We next describe our results for the sufficiently spread-out contact
process at $\lamb=\lambc$ for $d>4$.

\subsubsection{Results above four dimensions}\label{sss:resultsd>4}
We now state the results for the two-point function. In the statements,
$\sigma$ and $\Delta$ are defined in \refeq{sgm-def}--\refeq{Deltadef},
and we write $\|f\|_\infty=\sup_{x\in\Zd}|f(x)|$ for a function $f$ on
$\Zd$.

\begin{thm}\label{thm:2pt}
Let $d>4$ and $\delta\in(0,1\wedge\Delta\wedge\frac{d-4}2)$.  There is an
$L_0=L_0(d)$ such that, for $L\geq L_0$, there are positive and finite
constants $v=v(d,L)$, $A=A(d,L)$, $C_1=C_1(d)$ and $C_2=C_2(d)$ such that
\begin{gather}
\hat\tau_t^{\lambc}(\tfrac{k}{\sqrt{v\sigma^2t}})=A\,e^{-\frac{|k|^2}{2d}}\,
 \big[1+O(|k|^2(1+t)^{-\delta})+O((1+t)^{-(d-4)/2})\big],\lbeq{tauasy}\\
\frac1{\wtau_t^{\lambc}(0)}\sum_{x\in\Zd}|x|^2\tau_t^{\lambc}(x)=v\sigma^2t
 \,\big[1+O((1+t)^{-\delta})\big],\lbeq{taugyr}\\
C_1L^{-d}(1+t)^{-d/2}\leq\|\tau_t^{\lambc}\|_\infty\leq e^{-t}+C_2L^{-d}
 (1+t)^{-d/2},\lbeq{tausup}
\end{gather}
with the error estimate in \refeq{tauasy} uniform in $k\in\mR^d$ with
$|k|^2/\log(2+t)$ sufficiently small.
\end{thm}

The above results correspond to \cite[Theorem~1.1]{hs01}, where the
two-point function for sufficiently spread-out critical oriented percolation
with $d>4$ was proved to obey similar behaviour.  The proof in \cite{hs01}
is based on the inductive method of \cite{hs02}.  We apply a modified
version of this induction method to prove Theorem~\ref{thm:2pt}.  The proof
also reveals that
\begin{align}\lbeq{estimates}
\lambc=1+O(L^{-d}),&& A=1+O(L^{-d}),&& v=1+O(L^{-d}).
\end{align}
In a sequel to this paper \cite{hsa04b}, we will investigate the critical
point in more detail and prove that
\begin{align}\lbeq{betaexp}
\lambc-1=\sum_{n=2}^\infty D^{*n}(o)+O(L^{-2d}),
\end{align}
holds for $d>4$, where $D^{*n}$ is the $n$-fold convolution of $D$ in
$\Zd$.  In particular, when $D$ is defined by \refeq{Dexample}, we obtain
(see \cite[Theorem~1.2]{hsa04b})
\begin{align}\lbeq{DP99}
\lambc-1=L^{-d}\sum_{n=2}^\infty U^{\star n}(o)+O(L^{-d-1}),
\end{align}
where $U$ is the uniform probability density over $[-1,1]^d\subset\mR^d$,
and $U^{\star n}$ is the $n$-fold convolution of $U$ in $\mR^d$.
The above expression was already obtained in \cite{dp99},
but with a weaker error estimate.

Let $\gamma$ and $\beta$ be the critical exponents for the quantities in
\refeq{lambc-def}, defined as
\begin{align}
\chi(\lamb)\sim(\lambc-\lamb)^{-\gamma}\quad(\lamb<\lambc),&&
\theta(\lamb)\sim(\lamb-\lambc)^\beta\quad(\lamb>\lambc),
\end{align}
where we use ``$\sim$'' in an appropriate sense.  For example, the
strongest form of $\chi(\lamb)\sim(\lambc-\lamb)^{-\gamma}$ is that
there is a $C\in(0,\infty)$ such that
\begin{align}\lbeq{gammasharp}
\chi(\lamb)=[C+o(1)]\,(\lambc-\lamb)^{-\gamma},
\end{align}
where $o(1)$ tends to 0 as $\lamb\uaw\lambc$.  Other examples are
the weaker form
\begin{align}\lbeq{critexpobds}
\exists\,C_1,C_2\in(0,\infty):\quad C_1(\lambc-\lamb)^{-\gamma}\leq
 \chi(\lamb)\leq C_2(\lambc-\lamb)^{-\gamma},
\end{align}
and the even weaker form
\begin{align}
\chi(\lamb)=(\lamb-\lambc)^{-\gamma+o(1)}.
\end{align}
See also \cite[p.70]{Ligg99} for various ways to define the critical
exponents.

As discussed for oriented percolation in \cite[Section~1.2.1]{hs01},
\refeq{tauasy} and \refeq{tausup} imply finiteness at $\lamb=\lambc$
of the triangle function
\begin{align}\lbeq{triangle-def}
\triangledown(\lamb)=\int_0^\infty dt\int_0^tds\sum_{x,y\in\Zd}\tau_t^\lamb
 (y)\,\tau_{t-s}^\lamb(y-x)\,\tau_s^\lamb(x).
\end{align}
Extending the argument in \cite{ny93} for oriented percolation to the
continuous-time setting,  we conclude that $\triangledown(\lambc)<\infty$
implies the triangle condition of \cite{an84,ba91,bw98}, under which
$\gamma$ and $\beta$ are both equal to 1 in the form given in
\refeq{critexpobds}, independently of the value of $d$ \cite{bw98}.
Since these $d$-independent values also arise on the tree
\cite{Scho98,wu95}, we call them the mean-field values.  The results
\refeq{tauasy}--\refeq{taugyr} also show that the critical exponents
$\nu$ and $\eta$, defined as
\begin{align}
\frac1{\hat\tau_t^{\lambc}(0)}\sum_{x\in\Zd}|x|^2\tau_t^{\lambc}(x)\sim
 t^{2\nu},&& \hat\tau_t^{\lambc}(0)\sim t^\eta,
\end{align}
take on the mean-filed values $\nu=1/2$ and $\eta=0$, in the stronger form
given in \refeq{gammasharp}.  The result $\eta=0$ proves that the statement
in \cite[Proposition 4.39]{Ligg99} on the tree also holds for sufficiently
spread-out contact process on $\Zd$ for $d>4$.  See the remark below
\cite[Proposition 4.39]{Ligg99}.  Furthermore, following from bounds
established in the course of the proof of Theorem~\ref{thm:2pt},
we can extend the aforementioned result of \cite{bw98}, i.e.,
$\gamma=1$ in the form given in \refeq{critexpobds}, to the precise
asymptotics as in \refeq{gammasharp}.  We will prove this in Section~\ref{sec-outlinedisc}.

So far, $d>4$ is a sufficient condition for the mean-field behaviour for
the spread-out contact process.  It has been shown, using the hyperscaling
inequalities in \cite{s02}, that $d\geq4$
is also a necessary condition for the mean-field behaviour.  Therefore,
the upper critical dimension for the spread-out contact process is 4,
and one can expect log corrections in $d=4$.

In \cite{hsa04}, we will investigate the higher-point functions of the
critical spread-out contact process for $d>4$.  These higher-point functions
are defined for $\vec t\in[0,\infty)^{r-1}$ and $\vec x\in\Z^{d(r-1)}$ by
\begin{align}\lbeq{rpt-def}
\tau_{\vec t}^\lamb(\vec x)=\mP^\lamb(x_i\in\bC_{t_i}~\forall i=1,\dots,r-1).
\end{align}
The proof will be based on a lace expansion that expresses the $r$-point
function in terms of $s$-point functions with $s<r$.  On the arising
equation, we will then perform induction in $r$, with the results for
$r=2$ given by Theorem \ref{thm:2pt}.  We discuss the extension to the
higher point functions in somewhat more detail in Section~\ref{ss:expans},
where we discuss the lace expansion.  In order to bound the lace expansion
coefficients for the higher point functions, the upper bounds in
\refeq{tauasy} for $k=0$ and in \refeq{tausup} are crucial.

\subsubsection{Results below and at four dimensions}\label{sss:resultsd<=4}
We also consider the low-dimensional case, i.e., $d\leq4$.  In this case,
the contact process is believed {\it not} to exhibit the mean-field behaviour
as long as $L$ remains finite, and Gaussian asymptotics are not expected to
hold in this case.  However, we can prove local Gaussian behaviour when the
range grows in time as
\begin{align}\lbeq{Lt-def}
L_{\sT}=L_1T^b\qquad(T\geq1),
\end{align}
where $L_1\geq1$ is the initial infection range.  We denote by
$\sigma_{\sT}^2$ the variance of $D$ in this situation.  We assume that
\begin{align}\lbeq{alphadef}
\alpha=bd+\frac{d-4}2>0.
\end{align}
Our main result is the following.

\begin{thm}\label{thm:2pt-lowdim}
Let $d\leq4$ and $\delta\in(0,1\wedge\Delta\wedge\alpha)$.  Then, there is
a $\lamb_{\sT}=1+O(T^{-\mu})$ for some $\mu\in(0,\alpha-\delta)$ such that,
for sufficiently large $L_1$, there are positive and finite constants
$C_1=C_1(d)$ and $C_2=C_2(d)$ such that, for every $0<t\leq\log T$,
\begin{gather}
\hat\tau_{Tt}^{\lamb_T}(\tfrac{k}{\sqrt{\sigma_T^2Tt}})=e^{-\frac
 {|k|^2}{2d}}\,\big[1+O(T^{-\mu})+O(|k|^2(1+Tt)^{-\delta})\big],
 \lbeq{tauasy-lowdim}\\
\frac1{\hat\tau_{Tt}^{\lamb_T}(0)}\sum_{x\in\Zd}|x|^2\tau_{Tt}^{\lamb_T}
 (x)=\sigma_{\sT}^2Tt\,\big[1+O(T^{-\mu})+O((1+Tt)^{-\delta})\big],
 \lbeq{taugyr-lowdim}\\
C_1L_{\sT}^{-d}(1+Tt)^{-d/2}\leq\|\tau_{Tt}^{\lamb_T}\|_\infty
 \leq e^{-Tt}+C_2L_{\sT}^{-d}(1+Tt)^{-d/2},\lbeq{tausup-lowdim}
\end{gather}
with the error estimate in \refeq{tauasy-lowdim} uniform in $k\in\mR^d$
with $|k|^2/\log(2+Tt)$ sufficiently small.
\end{thm}

The upper bound on $t$ in the statement can be replaced by any slowly
varying function.  However, we use $\log T$ to make the statement more
concrete.  The proof of Theorem~\ref{thm:2pt-lowdim} follows the same
steps as the proof of Theorem~\ref{thm:2pt}.

First, we give a heuristic explanation of how \refeq{alphadef} arises.
Recall that, for $d>4$, $\triangledown(\lambc)<\infty$ is a sufficient
condition for the mean-field behaviour.  For $d\leq4$, since $\triangledown(\lamb_{\sT})$ cannot be defined in full space-time
as in \refeq{triangle-def}, we modify the triangle function as
\begin{align}
\triangledown_{\rm ld}(\lamb_{\sT})=\int_0^{T\log T}\!\!dt\int_0^tds
 \sum_{x,y\in\Zd}\tau_t^{\lamb_{\sT}}(y)\,\tau_{t-s}^{\lamb_{\sT}}(y-x)\,
 \tau_s^{\lamb_{\sT}}(x).
\end{align}
Using the upper bounds in \refeq{tauasy-lowdim} for $k=0$ and in
\refeq{tausup-lowdim}, we obtain
\begin{align}
\triangledown_{\rm ld}(\lamb_{\sT})\leq C^2\int_0^{T\log T}\!\!dt\int_0^t
 ds~(e^{-tT}+C_2L_{\sT}^{-d}T^{-d/2})\leq O(T^{-2})+O(T^{2-bd-d/2}\log^2T),
\end{align}
which is finite for all $T$ whenever $bd>\frac{4-d}2$.  We can find a
similar argument in \cite[Section~14]{Slad04}.

Next, we compare the ranges needed in our results and in the results of
Durrett and Perkins \cite{dp99}, in which the convergence of the rescaled
contact process to super-Brownian motion was proved.  As in \refeq{alphadef}
we need $bd>\frac{4-d}{2}$, while in \cite{dp99} $bd=1$ for all $d\geq3$.
For $d=2$, which is a critical case in the setting of \cite{dp99}, the model
with range $L_{\sT}^2=T\log T$ was also investigated.  In comparison, we
are allowed to use ranges that grow to infinity slower than the ranges in
\cite{dp99} when $d\geq 3$, but the range for $d=2$ in our results needs to
be larger than that in \cite{dp99}.  It would be of interest to investigate
whether Theorem~\ref{thm:2pt-lowdim} holds when $L_{\sT}^2=T\log T$ (or
even smaller) by adapting our proofs.

Finally, we give a conjecture on the asymptotics of $\lamb_{\sT}$ as
$T\uaw\infty$.  The role of $\lamb_{\sT}$ is a sort of critical value for
the contact process in the finite-time interval $[0,T\log T]$, and hence
$\lamb_{\sT}$ approximates the real critical value $\lamb_{{\rm c},\sT}$
that also converges to 1 in the mean-field limit $T\uaw\infty$.  We believe
that the leading term of $\lamb_{{\rm c},\sT}-1$, say $c_{\sT}$, is equal
to that of $\lamb_{\sT}-1$.  As we will discuss below
in Section~\ref{ss:lowdim}, $\lamb_{\sT}$ satisfies a type of recursion
relation \refeq{lambT-def}.  We expect that, for $d\leq4$, we may employ
the methods in \cite{hsa04b} to obtain
\begin{align}\lbeq{asymlambT}
\lamb_{\sT}=1+[1+O(T^{-\mu})]\int_0^{T\log T}\!\!dt~\int_{[-\pi,\pi]^d}
 \frac{d^dk}{(2\pi)^d}~\hat D^2_{\sT}(k)\,e^{-[1-\hat D_{\sT}(k)]t},
\end{align}
where $D_{\sT}$ equals $D$ with range $L_{\sT}$.  (In fact, the exponent
$\mu$ could be replaced by any positive number strictly smaller than
$\alpha$.)  The integral with respect to $t\in\mR_+$ converges when $d>2$,
and hence we may obtain for sufficiently large $T$ that
\begin{align}\lbeq{asymlambT3,4}
\lamb_{\sT}&=1+[1+O(T^{-\mu})]\bigg[\int_{[-\pi,\pi]^d}\frac{d^dk}
 {(2\pi)^d}~\frac{\hat D_{\sT}^2(k)}{1-\hat D_{\sT}(k)}+O(T^{-bd-
 \frac{d-2}2})\bigg]\nn\\
&=1+\sum_{n=2}^\infty D_{\sT}^{*n}(o)+O(L_{\sT}^{-d-\frac\mu{b}
 \wedge\frac{d-2}{2b}}),
\end{align}
where we use \refeq{Lt-def} and the fact that the sum in
\refeq{asymlambT3,4} is $O(L_{\sT}^{-d})$.  Based on our belief mentioned
above, this would be a stronger result than the result in \cite{dp99} when
$d=3,4$, where $c_{\sT}=\sum_{n=2}^\infty D_{\sT}^{*n}(o)$.  However, to
prove this conjecture, we may require serious further work using block
constructions used in \cite{dp99}.

\section{Outline of the proof}\label{ss:outline}
In this section, we provide an outline of the proof of our main results.
This section is organized as follows.  In Section~\ref{ss:discretize},
we explain what the discretized contact process is, and state the results
for the discretized contact process.  These results apply in particular
to oriented percolation, which is a special example of the discretized
contact process.  In Section~\ref{ss:expans}, we briefly explain the
lace expansion for the discretized contact process, and state the
bounds on the lace expansion coefficients in Section~\ref{sec-ble}.
In Section~\ref{sec-indimp}, we explain how to use induction to prove
the asymptotics for the discretized contact process.  In
Section~\ref{sec-outlinedisc}, we state the results concerning
the continuum limit, and show that the results for the discretized
contact process together with the continuum limit imply the main
results in Theorems~\ref{thm:2pt}--\ref{thm:2pt-lowdim}.

\subsection{Discretization}\label{ss:discretize}
By the graphical representation, the contact process can be
constructed as follows.  We consider $\Zd\times\mR_+$ as
space-time.  Along each time line $\{x\}\times\mR_+$, we place
points according to a Poisson process with intensity 1,
independently of the other time lines.  For each ordered pair of
distinct time lines from $\{x\}\times\mR_+$ to $\{y\}\times\mR_+$,
we place directed bonds $((x,\,t),(y,\,t))$, $t\geq0$, according to
a Poisson process with intensity $\lamb\,D(y-x)$,
independently of the other Poisson processes.  A site $(x,s)$ is
said to be {\it connected to} $(y,t)$ if either $(x,s)= (y,t)$ or
there is a non-zero path in $\Zd\times\mR_+$ from $(x,s)$ to
$(y,t)$ using the Poisson bonds and time line segments traversed
in the increasing time direction without traversing the Poisson
points.  The law of $\bC_t$ defined in Section~\ref{ss:results} is
equivalent to that of $\{x\in \Zd:(o,\,0)$ is connected to
$(x,\,t)\}$. See also \cite[Section~I.1]{Ligg99}.

\begin{figure}
\begin{center}
\end{center}
\begin{center}
\includegraphics[scale = 0.45]{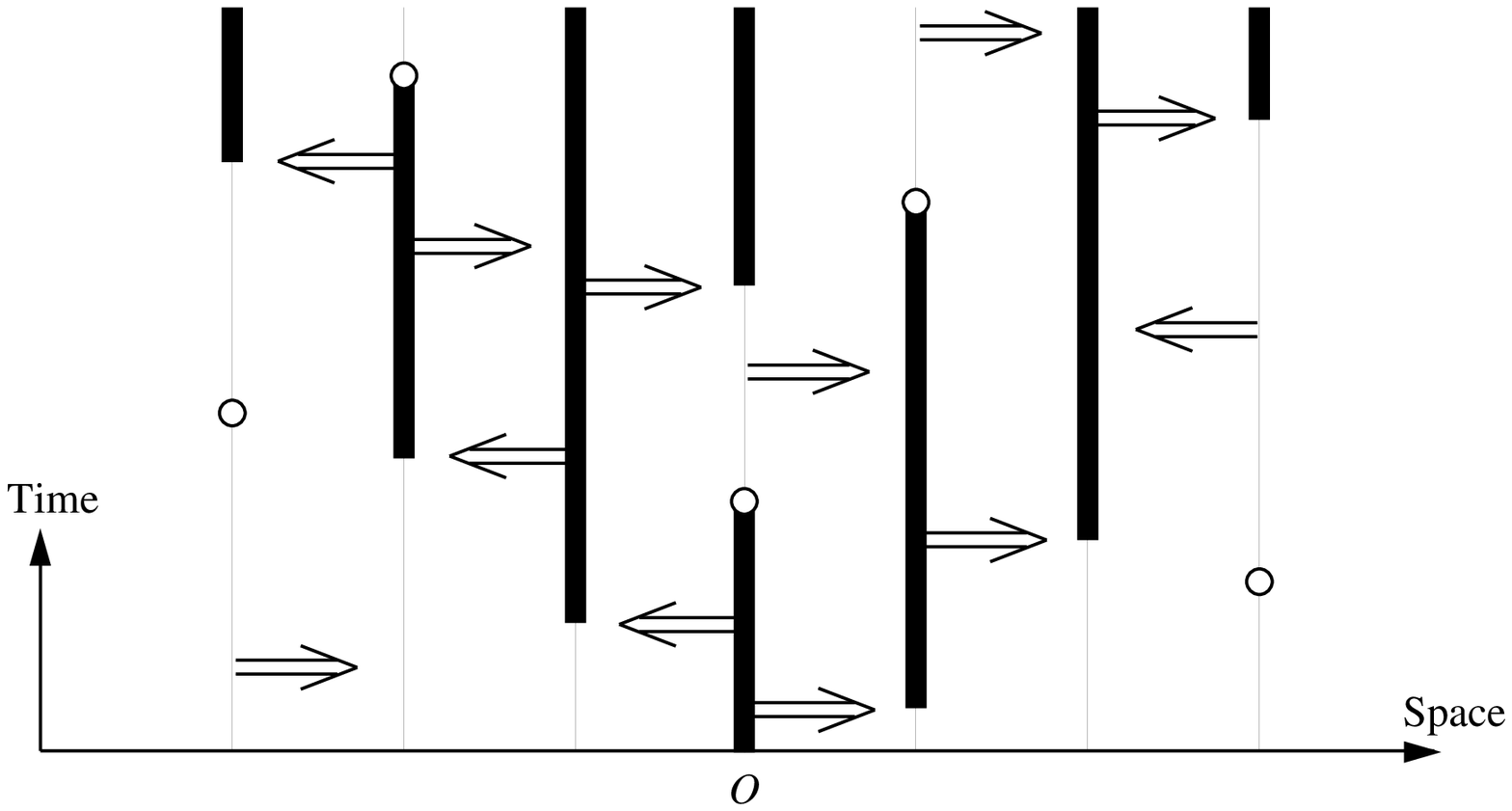}\hskip 1cm
\includegraphics[scale = 0.45]{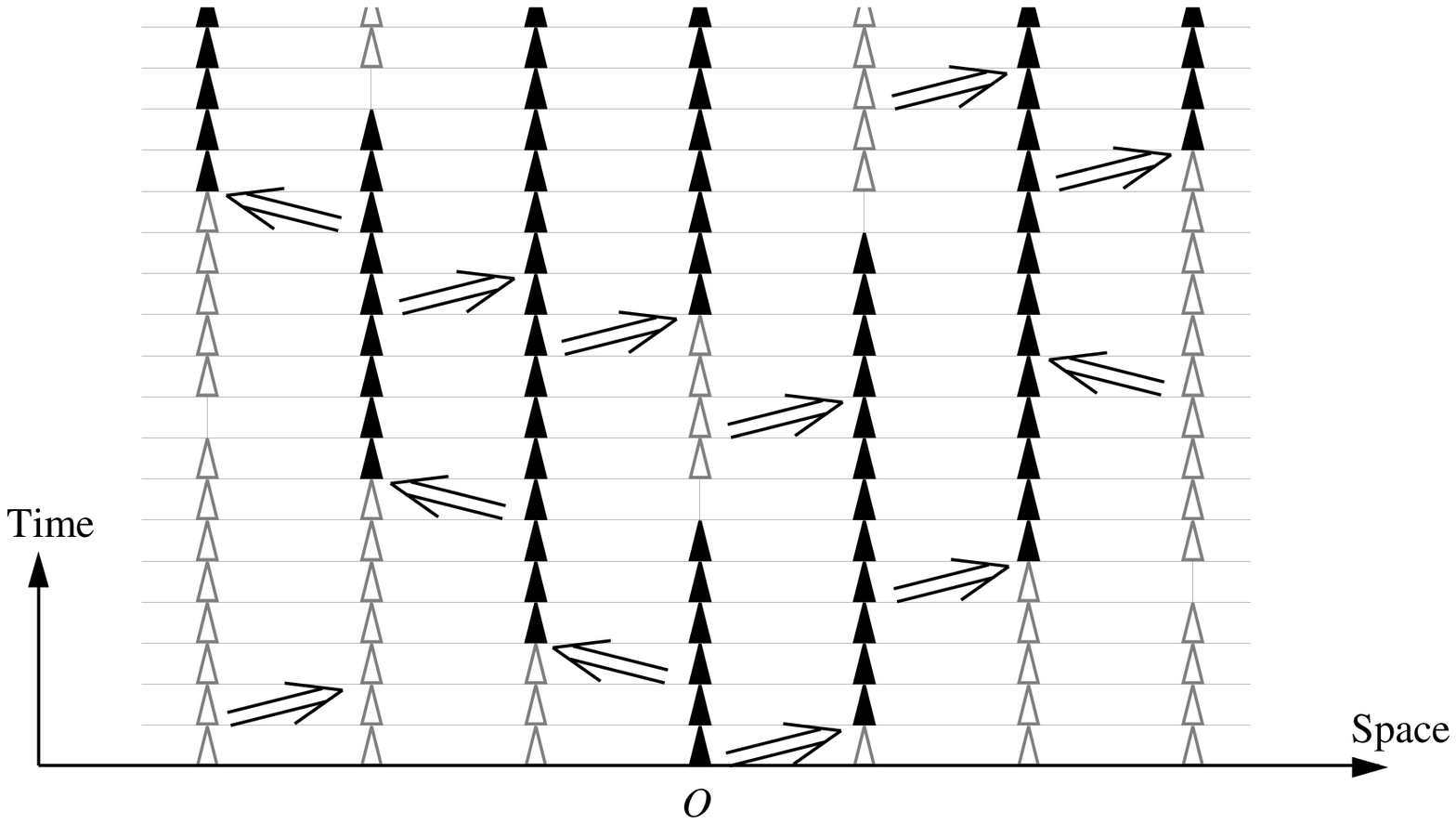}
\label{fig-1}
\caption{Graphical representation of the contact process and the discretized
contact process.}
\end{center}
\end{figure}

Inspired by this percolation structure in space-time
and following \cite{s01}, we consider an oriented percolation
approximation in $\Zd\times\vep\Zp$ to the contact process, where
$\vep \in(0,1]$ is a discretization parameter.  We call
this approximation the {\it discretized contact process},
and it is defined as follows. A directed pair
$b=((x,t),(y,t+\vep))$ of sites in $\Zd\times\vep\Zp$ is called a
{\it bond}.  In particular, $b$ is a {\it temporal bond} if
$x=y$, otherwise $b$ is a {\it spatial bond}.  Each bond is either
{\it occupied} or {\it vacant} independently of the other bonds,
and a bond $b=((x,t),(y,t+\vep))$ is occupied with probability
\begin{align}\lbeq{bprob}
p_\vep(y-x)=\begin{cases}
 1-\vep,&\text{if }x=y,\\ \lamb\vep\,D(y-x),&\text{if }x\ne y,
 \end{cases}
\end{align}
provided that $\|p_{\vep}\|_\infty\leq1$.  We denote the associated
probability measure by $\mP_\vep^\lamb$.  It is proved in \cite{bg91}
that $\mP_\vep^\lamb$ weakly converges to $\mP^\lamb$ as $\vep\daw0$.
See Figure~\ref{fig-1} for a graphical representation of the contact
process and the discretized contact process.  As explained in more
detail in Section~\ref{ss:expans}, we prove our main results by proving
the results first for the discretized contact process, and then taking
the continuum limit when $\vep\daw0$.

We also emphasize that the discretized contact process with
$\vep=1$ is equivalent to oriented percolation, for which
$\lamb\in[0,\|D\|_\infty^{-1}]$ is the expected number of
occupation bonds per site.

We denote by $(x,s)\conn(y,t)$ the event that $(x,s)$ is {\it
connected to} $(y,t)$, i.e., either $(x,s)=(y,t)$ or there is a
non-zero path in $\Zd\times\vep\Zp$ from $(x,s)$ to $(y,t)$
consisting of occupied bonds.  The {\it two-point function} is
defined as
    \eq\lbeq{disc2pt-def}
        \tau_{t;\vep}^\lamb(x)=\mP_\vep^\lamb((o,0)\conn(x,t)).
    \en
Similarly to \refeq{lambc-def}, the discretized contact process
has a critical value $\lambda_{\rm c}^{\smallsup{\vep}}$ satisfying
\begin{align}\lbeq{lambcvep-def}
\vep\sum_{t\in\vep\Zp}\hat\tau_{t;\vep}^\lamb(0)\begin{cases}
 <\infty,&\mbox{if }\lamb<\lambc^{\sss(\vep)},\\
 =\infty,&\mbox{if }\lamb\geq\lambc^{\sss(\vep)},
 \end{cases}&&
\lim_{t\uaw\infty}\mP^\lamb_\vep(\bC_t\ne\varnothing)\begin{cases}
 =0,&\mbox{if }\lamb\leq\lambc^{\sss(\vep)},\\
 >0,&\mbox{if }\lamb>\lambc^{\sss(\vep)}.
 \end{cases}
\end{align}

The main result for the discretized contact process with $\vep\in(0,1]$
is the following theorem:

\begin{prop}[Discretized results for $d>4$]\label{thm-disc}
Let $d>4$ and $\delta\in(0,1\wedge\Delta\wedge\frac{d-4}2)$.  Then,
there is an $L_0=L_0(d)$ such that, for $L\geq L_0$, there are positive
and finite constants $v^{\sss(\vep)}=v^{\sss(\vep)}(d,L)$,
$A^{\sss(\vep)}=A^{\sss(\vep)}(d,L)$, $C_1(d)$ and $C_2(d)$ such that
\begin{gather}
\hat\tau_{t;\vep}^{\lambc^{(\vep)}}(\tfrac{k}{\sqrt{v^{\smallsup{\vep}}
 \sigma^2t}})=A^{\sss(\vep)}e^{-\frac{|k|^2}{2d}}\,\big[1+O(|k|^2(1+t)^{-
 \delta})+O((1+t)^{-(d-4)/2})\big],\lbeq{tauasyvep}\\
\frac1{\hat\tau_{t;\vep}^{\lambc^{(\vep)}}(0)}\sum_{x\in\Zd}|x|^2\tau_{t;
 \vep}^{\lambc^{(\vep)}}(x)=v^{\sss(\vep)}\sigma^2t\,\big[1+O((1+t)^{-
 \delta})\big],\lbeq{taugyrvep}\\
C_1L^{-d}(1+t)^{-d/2}\leq\|\tau_{t;\vep}^{\lambc^{(\vep)}}\|_\infty\leq
 (1-\vep)^{t/\vep}+C_2L^{-d}(1+t)^{-d/2},\lbeq{tausupvep}
    \end{gather}
where all error terms are uniform in $\vep\in(0,1]$.  The error estimate
in \refeq{tauasyvep} is uniform in $k\in\mR^d$ with $|k|^2/\log(2+t)$
sufficiently small.
\end{prop}

Proposition \ref{thm-disc} is the discrete analog of
Theorem~\ref{thm:2pt}.  The uniformity in $\vep$ of the error terms is
crucial, as this will allow us to take the limit $\vep\daw0$ and to
conclude the results in Theorem~\ref{thm:2pt} from the corresponding
statements in Proposition~\ref{thm-disc}.  In particular,
Proposition~\ref{thm-disc} applied to oriented percolation
(i.e., $\vep=1$) reproves \cite[Theorem 1.1]{hs01}.

The discretized version of Theorem~\ref{thm:2pt-lowdim} is given
in the following proposition:

\begin{prop}[Discretized results for $d\leq4$]\label{thm-disc2}
Let $d\leq4$ and $\delta\in(0,1\wedge\Delta\wedge\alpha)$.  Then, there is
a $\lamb_{\sT}=1+O(T^{-\mu})$ for some $\mu\in(0,\alpha-\delta)$ such that,
for sufficiently large $L_1$, there are positive and finite constants
$C_1=C_1(d)$ and $C_2=C_2(d)$ such that, for every $0<t\leq\log T$,
\begin{gather}
\hat\tau_{Tt;\vep}^{\lamb_T}(\tfrac{k}{\sqrt{\sigma_T^2Tt}})=e^{-\frac
 {|k|^2}{2d}}\,\big[1+O(T^{-\mu})+O(|k|^2(1+Tt)^{-\delta})\big],
 \lbeq{tauasy-lowdimdis}\\
\frac1{\hat\tau_{Tt;\vep}^{\lamb_T}(0)}\sum_{x\in\Zd}|x|^2\tau_{Tt;\vep}
 ^{\lamb_T}(x)=\sigma_{\sT}^2Tt\,\big[1+O(T^{-\mu})+O((1+Tt)^{-\delta})
 \big],\lbeq{taugyr-lowdimdis}\\
C_1L_{\sT}^{-d}(1+Tt)^{-d/2}\leq\|\tau_{Tt;\vep}^{\lamb_T}\|_\infty\leq
 (1-\vep)^{Tt/\vep}+C_2 L_{\sT}^{-d}(1+Tt)^{-d/2},\lbeq{tausup-lowdimdis}
\end{gather}
where all error terms are uniform in $\vep\in (0,1]$, and the error
estimate in \refeq{tauasy-lowdimdis} is uniform in $k\in\mR^d$ with
$|k|^2/\log(2+Tt)$ sufficiently small.
\end{prop}

Note that Proposition \ref{thm-disc2} applies also to oriented percolation,
for which $\vep=1$.

\subsection{Expansion}
\label{ss:expans}

The proof of Proposition~\ref{thm-disc} makes use of the lace expansion,
which is an expansion for the two-point function.  We postpone the
derivation of the expansion to Section~\ref{s:lace}, and here we provide
only a brief motivation. We also motivate why we discretize time for the
contact process.

We make use of the convolution of functions, which is defined for
absolutely summable functions $f,g$ on $\Zd$ by
\begin{align}
(f*g)(x)=\sum_{y\in\Zd}f(y)\,g(x-y).
\end{align}

We first motivate the basic idea underlying the expansion, similarly as in
\cite[Section~2.1.1]{hs01}, by considering the much simpler corresponding
expansion for continuous-time random walk. For continuous-time random walk
making jumps from $x$ to $y$ at rate $\lambda D(y-x)$ with killing rate
$1-\lambda$, we have the partial differential equation
\begin{align}\lbeq{SRWx}
\partial_tq_t^\lamb(x)=\lamb\,(D*q_t^\lamb)(x)-q_t^\lamb(x),
\end{align}
where $q_t^\lamb(x)$ is the probability that continuous-time random walk
started at $o\in\Zd$ is at $x\in\Zd$ at time $t$.  By taking the Fourier
transform, we obtain
\begin{align}\lbeq{SRW}
\partial_t\hat q_t^\lamb(k)=-[1-\lamb\hat D(k)]\,\hat q_t^\lamb(k).
\end{align}
In this simple case, the above equation is readily solved to yield that
\begin{align}\lbeq{solSRW}
\hat q_t^\lamb(k)=e^{-[1-\lamb\hat D(k)]t}.
\end{align}
We see that $\lambda=1$ is the critical value, and the central limit
theorem at $\lambda=\lambda_{\rm c}=1$ follows by a Taylor expansion
of $1-\hat{D}(k)$ for small $k$, yielding
\begin{align}
\hat q_t^1\big(\tfrac{k}{\sqrt{\sigma^2 t}}\big)=e^{-\frac{|k|^2}{2d}}
 \,[1+o(1)],
\end{align}
where $|k|^2=\sum_{j=1}^d k_i^2$ (recall also \refeq{sgm-def}).

The above solution is quite specific to continuous-time random walk.
When we would have a more difficult function on the right-hand side of
\refeq{SRW}, such as $-[1-\lamb\hat D(k)]\,\hat q_{t-1}^\lamb(k)$,
it would be much more involved to solve the above equation, even though
one would expect that the central limit theorem at the critical value
still holds.

A more robust proof of central limit behaviour uses induction
in time $t$. Since time is continuous, we first discretize time.
The two-point function for discretized continuous-time random
walk is defined by setting $q_{0;\vep}^\lamb(x)=\delta_{0,x}$ and
(recall \refeq{bprob})
    \eq \lbeq{2-pt.smp.1}
    q_{t;\vep}^\lambda(x)=p_{\varepsilon}^{*t/\vep}(x)
    \quad \quad (t \in \varepsilon\N).
    \en
To obtain a recursion relation for $q_{t;\vep}^\lambda(x)$, we
simply observe that by independence of the underlying random walk
    \eq
    \lbeq{2-pt.smp.2}
    q_{t;\vep}^\lambda(x) =  (p_{\varepsilon} * q_{t-\varepsilon;\vep}^\lambda)(x)
    \quad \quad (t \in \varepsilon\N).
    \en
We can think of this as a simple version of the lace expansion, applied
to random walk, which has no interaction.

For the discretized continuous-time random walk,
we can use induction in $n$ for all $t=n\vep$. If we can further
show that the arising error terms are uniform in $\vep$, then we can take the
continuum limit $\vep \downarrow 0$ afterwards, and obtain the result
for the continuous-time model. The above proof is more robust, and
can for instance be used to deal with the situation where the right-hand side of \refeq{SRW}
equals $-[1-\lambda\hat{D}(k)] \hat q_{t-1}^{\lambda} (k).$ This robustness
of the proof is quite valuable when we wish to apply it to the contact
process.


The identity \refeq{2-pt.smp.2} can be solved using the Fourier transform
to give
\begin{align}
\hat q_{t;\vep}^\lamb(k)=\hat p_\vep(k)^{t/\vep}=[1-\vep+\lamb\vep\hat
 D(k)]^{t/\vep}=e^{-[1-\lamb\hat D(k)]t+O(t\vep[1-\hat D(k)]^2)}.
\end{align}
We note that the limit of
$[\hat q_{t;\vep}(k)-\hat q_{t-\vep;\vep}(k)]/\vep$ exists and equals
\refeq{SRW}.  In order to obtain the central limit theorem, we divide
$k$ by $\sqrt{\sigma^2t}$.  Then, uniformly in $\varepsilon>0$, we have
\begin{align}\lbeq{qvepconv}
\hat q_{t;\vep}^1\big(\tfrac{k}{\sqrt{\sigma^2t}}\big)=e^{-\frac{|k|^2}
 {2d}+O(|k|^{2+2\Delta}t^{-\Delta})+O(\vep|k|^4t^{-1})}.
\end{align}
Therefore, the central limit theorem holds uniformly in $\vep>0$.

We follow Mark Kac's adagium: ``Be wise, discretize!'' for two reasons.
Firstly, discretizing time allows us to obtain an expansion as in
\refeq{SRWx}, and secondly, it allows us to analyse the arising equation.
The lace expansion, which is explained in more detail below, can be used
for the contact process to produce an equation of the form
\begin{align}\lbeq{tautxcont}
\partial_t\hat\tau_t^\lamb(k)=-[1-\lamb\hat D(k)]\,\hat\tau_t^\lamb(k)
 +\int_0^tds~\hat\pi_s^\lamb(k)\,\hat\tau_{t-s}^\lamb(k),
\end{align}
where $\hat\pi_s^\lamb$ are certain expansion coefficients.  In order to
derive the equation \refeq{tautxcont}, we use that the discretized contact
process is oriented percolation, for which lace expansions have been
derived in the literature \cite{hs01,ny93,ny95,s00,s01}.  Clearly, the
equation \refeq{tautxcont} is much more complicated than the corresponding
equation for simple random walk in \refeq{SRWx}.  Therefore, a simple
solution to the equation as in \refeq{solSRW} is impossible.  We see no way
to analyse the partial differential equation in \refeq{tautxcont} other
than to discretize time combined with induction.  It would be of interest
to investigate whether \refeq{tautxcont} can be used directly.

\begin{figure}[t]
\begin{center}
\vskip6cm
\setlength{\unitlength}{0.00067in}
{\hskip-4cm
\begin{picture}(400,0)
{\put(-1400,0){\makebox(0,0)[lb]{\wtr}}
\put(-1400,100){\makebox(0,0)[lb]{\wtr}}
\put(-1400,200){\makebox(0,0)[lb]{\wtr}}
\put(-1400,300){\makebox(0,0)[lb]{\wtr}}
\put(-1400,400){\makebox(0,0)[lb]{\wtr}}
\put(-1400,500){\makebox(0,0)[lb]{\wtr}}
\put(-1400,600){\makebox(0,0)[lb]{\wtr}}
\put(-1400,800){\makebox(0,0)[lb]{\wtr}}
\put(-1400,900){\makebox(0,0)[lb]{\wtr}}
\put(-1400,1000){\makebox(0,0)[lb]{\wtr}}
\put(-1400,1100){\makebox(0,0)[lb]{\wtr}}
\put(-1400,1200){\makebox(0,0)[lb]{\wtr}}
\put(-1400,1300){\makebox(0,0)[lb]{\wtr}}
\put(-1400,1400){\makebox(0,0)[lb]{\wtr}}
\put(-1400,1500){\makebox(0,0)[lb]{\wtr}}
\put(-1400,1600){\makebox(0,0)[lb]{\wtr}}
\put(-1400,1700){\makebox(0,0)[lb]{\wtr}}
\put(-1400,1800){\makebox(0,0)[lb]{\wtr}}
\put(-1400,1900){\makebox(0,0)[lb]{\wtr}}
\put(-1400,2000){\makebox(0,0)[lb]{\wtr}}
\put(-1400,2100){\makebox(0,0)[lb]{\wtr}}
\put(-1400,2200){\makebox(0,0)[lb]{\wtr}}
\put(-1400,2300){\makebox(0,0)[lb]{\wtr}}
\put(-1400,2400){\makebox(0,0)[lb]{\wtr}}
\put(-1400,2500){\makebox(0,0)[lb]{\wtr}}
\put(-1400,2600){\makebox(0,0)[lb]{\wtr}}
\put(-1400,2800){\makebox(0,0)[lb]{\wtr}}
\put(-1400,2900){\makebox(0,0)[lb]{\wtr}}
\put(-1400,3000){\makebox(0,0)[lb]{\wtr}}
\put(-1400,3100){\makebox(0,0)[lb]{\wtr}}

\put(-1000,0){\makebox(0,0)[lb]{\wtr}}
\put(-1000,100){\makebox(0,0)[lb]{\wtr}}
\put(-1000,200){\makebox(0,0)[lb]{\wtr}}
\put(-1000,300){\makebox(0,0)[lb]{\wtr}}
\put(-1000,500){\makebox(0,0)[lb]{\wtr}}
\put(-1000,600){\makebox(0,0)[lb]{\wtr}}
\put(-1000,700){\makebox(0,0)[lb]{\wtr}}
\put(-1000,800){\makebox(0,0)[lb]{\wtr}}
\put(-1000,900){\makebox(0,0)[lb]{\wtr}}
\put(-1000,1000){\makebox(0,0)[lb]{\wtr}}
\put(-1000,1200){\makebox(0,0)[lb]{\wtr}}
\put(-1000,1300){\makebox(0,0)[lb]{\wtr}}
\put(-1000,1400){\makebox(0,0)[lb]{\wtr}}
\put(-1000,1500){\makebox(0,0)[lb]{\wtr}}
\put(-1000,1600){\makebox(0,0)[lb]{\wtr}}
\put(-1000,1700){\makebox(0,0)[lb]{\wtr}}
\put(-1000,1800){\makebox(0,0)[lb]{\btr}}
\put(-1000,1900){\makebox(0,0)[lb]{\btr}}
\put(-1000,2000){\makebox(0,0)[lb]{\btr}}
\put(-1000,2200){\makebox(0,0)[lb]{\wtr}}
\put(-1000,2300){\makebox(0,0)[lb]{\wtr}}
\put(-1000,2400){\makebox(0,0)[lb]{\wtr}}
\put(-1000,2500){\makebox(0,0)[lb]{\wtr}}
\put(-1000,2600){\makebox(0,0)[lb]{\wtr}}
\put(-1000,2700){\makebox(0,0)[lb]{\wtr}}
\put(-1000,2800){\makebox(0,0)[lb]{\wtr}}
\put(-1000,2900){\makebox(0,0)[lb]{\wtr}}
\put(-1000,3000){\makebox(0,0)[lb]{\wtr}}
\put(-1000,3100){\makebox(0,0)[lb]{\wtr}}

\put(-600,0){\makebox(0,0)[lb]{\wtr}}
\put(-600,100){\makebox(0,0)[lb]{\wtr}}
\put(-600,200){\makebox(0,0)[lb]{\wtr}}
\put(-600,300){\makebox(0,0)[lb]{\wtr}}
\put(-600,400){\makebox(0,0)[lb]{\wtr}}
\put(-600,500){\makebox(0,0)[lb]{\wtr}}
\put(-600,600){\makebox(0,0)[lb]{\wtr}}
\put(-600,700){\makebox(0,0)[lb]{\wtr}}
\put(-600,800){\makebox(0,0)[lb]{\wtr}}
\put(-600,900){\makebox(0,0)[lb]{\wtr}}
\put(-600,1000){\makebox(0,0)[lb]{\wtr}}
\put(-600,1100){\makebox(0,0)[lb]{\wtr}}
\put(-600,1200){\makebox(0,0)[lb]{\wtr}}
\put(-600,1400){\makebox(0,0)[lb]{\wtr}}
\put(-600,1500){\makebox(0,0)[lb]{\wtr}}
\put(-600,1600){\makebox(0,0)[lb]{\btr}}
\put(-600,1700){\makebox(0,0)[lb]{\btr}}
\put(-600,1800){\makebox(0,0)[lb]{\btr}}
\put(-600,1900){\makebox(0,0)[lb]{\btr}}
\put(-600,2000){\makebox(0,0)[lb]{\btr}}
\put(-600,2100){\makebox(0,0)[lb]{\btr}}
\put(-600,2200){\makebox(0,0)[lb]{\btr}}
\put(-600,2300){\makebox(0,0)[lb]{\btr}}
\put(-600,2400){\makebox(0,0)[lb]{\btr}}
\put(-600,2500){\makebox(0,0)[lb]{\btr}}
\put(-600,2600){\makebox(0,0)[lb]{\btr}}
\put(-600,2700){\makebox(0,0)[lb]{\btr}}
\put(-600,2800){\makebox(0,0)[lb]{\btr}}
\put(-600,2900){\makebox(0,0)[lb]{\btr}}
\put(-600,3000){\makebox(0,0)[lb]{\btr}}
\put(-600,3100){\makebox(0,0)[lb]{\btr}}

\put(-200,0){\makebox(0,0)[lb]{\wtr}}
\put(-200,100){\makebox(0,0)[lb]{\wtr}}
\put(-200,200){\makebox(0,0)[lb]{\btr}}
\put(-200,300){\makebox(0,0)[lb]{\btr}}
\put(-200,400){\makebox(0,0)[lb]{\btr}}
\put(-200,500){\makebox(0,0)[lb]{\btr}}
\put(-200,700){\makebox(0,0)[lb]{\wtr}}
\put(-200,800){\makebox(0,0)[lb]{\wtr}}
\put(-200,900){\makebox(0,0)[lb]{\wtr}}
\put(-200,1000){\makebox(0,0)[lb]{\btr}}
\put(-200,1100){\makebox(0,0)[lb]{\btr}}
\put(-200,1200){\makebox(0,0)[lb]{\btr}}
\put(-200,1300){\makebox(0,0)[lb]{\btr}}
\put(-200,1400){\makebox(0,0)[lb]{\btr}}
\put(-200,1500){\makebox(0,0)[lb]{\btr}}
\put(-200,1600){\makebox(0,0)[lb]{\btr}}
\put(-200,1700){\makebox(0,0)[lb]{\btr}}
\put(-200,1800){\makebox(0,0)[lb]{\btr}}
\put(-200,1900){\makebox(0,0)[lb]{\btr}}
\put(-200,2000){\makebox(0,0)[lb]{\btr}}
\put(-200,2100){\makebox(0,0)[lb]{\btr}}
\put(-200,2200){\makebox(0,0)[lb]{\btr}}
\put(-200,2300){\makebox(0,0)[lb]{\btr}}
\put(-200,2400){\makebox(0,0)[lb]{\btr}}
\put(-200,2500){\makebox(0,0)[lb]{\btr}}
\put(-200,2600){\makebox(0,0)[lb]{\btr}}
\put(-200,2700){\makebox(0,0)[lb]{\btr}}
\put(-200,2800){\makebox(0,0)[lb]{\btr}}
\put(-200,2900){\makebox(0,0)[lb]{\btr}}
\put(-200,3000){\makebox(0,0)[lb]{\btr}}
\put(-200,3100){\makebox(0,0)[lb]{\btr}}

\put(200,0){\makebox(0,0)[lb]{\btr}}
\put(200,100){\makebox(0,0)[lb]{\btr}}
\put(200,200){\makebox(0,0)[lb]{\btr}}
\put(200,300){\makebox(0,0)[lb]{\btr}}
\put(200,400){\makebox(0,0)[lb]{\btr}}
\put(200,500){\makebox(0,0)[lb]{\btr}}
\put(200,600){\makebox(0,0)[lb]{\btr}}
\put(200,700){\makebox(0,0)[lb]{\btr}}
\put(200,800){\makebox(0,0)[lb]{\btr}}
\put(200,900){\makebox(0,0)[lb]{\btr}}
\put(200,1000){\makebox(0,0)[lb]{\btr}}
\put(200,1100){\makebox(0,0)[lb]{\btr}}
\put(200,1200){\makebox(0,0)[lb]{\btr}}
\put(200,1300){\makebox(0,0)[lb]{\btr}}
\put(200,1400){\makebox(0,0)[lb]{\btr}}
\put(200,1500){\makebox(0,0)[lb]{\btr}}
\put(200,1600){\makebox(0,0)[lb]{\btr}}
\put(200,1800){\makebox(0,0)[lb]{\wtr}}
\put(200,1900){\makebox(0,0)[lb]{\wtr}}
\put(200,2000){\makebox(0,0)[lb]{\wtr}}
\put(200,2100){\makebox(0,0)[lb]{\wtr}}
\put(200,2200){\makebox(0,0)[lb]{\wtr}}
\put(200,2300){\makebox(0,0)[lb]{\wtr}}
\put(200,2400){\makebox(0,0)[lb]{\wtr}}
\put(200,2500){\makebox(0,0)[lb]{\wtr}}
\put(200,2600){\makebox(0,0)[lb]{\wtr}}
\put(200,2700){\makebox(0,0)[lb]{\btr}}
\put(200,2800){\makebox(0,0)[lb]{\btr}}
\put(200,2900){\makebox(0,0)[lb]{\btr}}
\put(200,3000){\makebox(0,0)[lb]{\btr}}
\put(200,3100){\makebox(0,0)[lb]{\btr}}

\put(600,0){\makebox(0,0)[lb]{\wtr}}
\put(600,100){\makebox(0,0)[lb]{\wtr}}
\put(600,200){\makebox(0,0)[lb]{\wtr}}
\put(600,300){\makebox(0,0)[lb]{\btr}}
\put(600,400){\makebox(0,0)[lb]{\btr}}
\put(600,500){\makebox(0,0)[lb]{\btr}}
\put(600,600){\makebox(0,0)[lb]{\btr}}
\put(600,700){\makebox(0,0)[lb]{\btr}}
\put(600,900){\makebox(0,0)[lb]{\wtr}}
\put(600,1000){\makebox(0,0)[lb]{\wtr}}
\put(600,1100){\makebox(0,0)[lb]{\wtr}}
\put(600,1200){\makebox(0,0)[lb]{\wtr}}
\put(600,1300){\makebox(0,0)[lb]{\wtr}}
\put(600,1400){\makebox(0,0)[lb]{\wtr}}
\put(600,1500){\makebox(0,0)[lb]{\wtr}}
\put(600,1600){\makebox(0,0)[lb]{\btr}}
\put(600,1700){\makebox(0,0)[lb]{\btr}}
\put(600,1800){\makebox(0,0)[lb]{\btr}}
\put(600,2000){\makebox(0,0)[lb]{\wtr}}
\put(600,2100){\makebox(0,0)[lb]{\wtr}}
\put(600,2200){\makebox(0,0)[lb]{\wtr}}
\put(600,2300){\makebox(0,0)[lb]{\wtr}}
\put(600,2400){\makebox(0,0)[lb]{\wtr}}
\put(600,2500){\makebox(0,0)[lb]{\wtr}}
\put(600,2600){\makebox(0,0)[lb]{\wtr}}
\put(600,2700){\makebox(0,0)[lb]{\wtr}}
\put(600,2900){\makebox(0,0)[lb]{\wtr}}
\put(600,3000){\makebox(0,0)[lb]{\wtr}}
\put(600,3100){\makebox(0,0)[lb]{\wtr}}

\put(1000,0){\makebox(0,0)[lb]{\wtr}}
\put(1000,100){\makebox(0,0)[lb]{\wtr}}
\put(1000,200){\makebox(0,0)[lb]{\wtr}}
\put(1000,300){\makebox(0,0)[lb]{\wtr}}
\put(1000,400){\makebox(0,0)[lb]{\btr}}
\put(1000,500){\makebox(0,0)[lb]{\btr}}
\put(1000,600){\makebox(0,0)[lb]{\btr}}
\put(1000,700){\makebox(0,0)[lb]{\btr}}
\put(1000,900){\makebox(0,0)[lb]{\wtr}}
\put(1000,1000){\makebox(0,0)[lb]{\wtr}}
\put(1000,1100){\makebox(0,0)[lb]{\wtr}}
\put(1000,1200){\makebox(0,0)[lb]{\wtr}}
\put(1000,1300){\makebox(0,0)[lb]{\wtr}}
\put(1000,1500){\makebox(0,0)[lb]{\wtr}}
\put(1000,1600){\makebox(0,0)[lb]{\wtr}}
\put(1000,1700){\makebox(0,0)[lb]{\wtr}}
\put(1000,1800){\makebox(0,0)[lb]{\wtr}}
\put(1000,1900){\makebox(0,0)[lb]{\wtr}}
\put(1000,2000){\makebox(0,0)[lb]{\wtr}}
\put(1000,2100){\makebox(0,0)[lb]{\wtr}}
\put(1000,2200){\makebox(0,0)[lb]{\wtr}}
\put(1000,2300){\makebox(0,0)[lb]{\wtr}}
\put(1000,2400){\makebox(0,0)[lb]{\wtr}}
\put(1000,2500){\makebox(0,0)[lb]{\wtr}}
\put(1000,2600){\makebox(0,0)[lb]{\wtr}}
\put(1000,2700){\makebox(0,0)[lb]{\wtr}}
\put(1000,2800){\makebox(0,0)[lb]{\wtr}}
\put(1000,3000){\makebox(0,0)[lb]{\wtr}}
\put(1000,3100){\makebox(0,0)[lb]{\wtr}}

\put(1400,0){\makebox(0,0)[lb]{\wtr}}
\put(1400,100){\makebox(0,0)[lb]{\wtr}}
\put(1400,200){\makebox(0,0)[lb]{\wtr}}
\put(1400,300){\makebox(0,0)[lb]{\wtr}}
\put(1400,400){\makebox(0,0)[lb]{\wtr}}
\put(1400,500){\makebox(0,0)[lb]{\wtr}}
\put(1400,600){\makebox(0,0)[lb]{\wtr}}
\put(1400,700){\makebox(0,0)[lb]{\btr}}
\put(1400,800){\makebox(0,0)[lb]{\btr}}
\put(1400,900){\makebox(0,0)[lb]{\btr}}
\put(1400,1000){\makebox(0,0)[lb]{\btr}}
\put(1400,1200){\makebox(0,0)[lb]{\wtr}}
\put(1400,1300){\makebox(0,0)[lb]{\wtr}}
\put(1400,1400){\makebox(0,0)[lb]{\wtr}}
\put(1400,1600){\makebox(0,0)[lb]{\wtr}}
\put(1400,1700){\makebox(0,0)[lb]{\wtr}}
\put(1400,1800){\makebox(0,0)[lb]{\wtr}}
\put(1400,1900){\makebox(0,0)[lb]{\wtr}}
\put(1400,2000){\makebox(0,0)[lb]{\wtr}}
\put(1400,2100){\makebox(0,0)[lb]{\wtr}}
\put(1400,2200){\makebox(0,0)[lb]{\wtr}}
\put(1400,2300){\makebox(0,0)[lb]{\wtr}}
\put(1400,2400){\makebox(0,0)[lb]{\wtr}}
\put(1400,2500){\makebox(0,0)[lb]{\wtr}}
\put(1400,2700){\makebox(0,0)[lb]{\wtr}}
\put(1400,2800){\makebox(0,0)[lb]{\wtr}}
\put(1400,2900){\makebox(0,0)[lb]{\wtr}}
\put(1400,3000){\makebox(0,0)[lb]{\wtr}}
\put(1400,3100){\makebox(0,0)[lb]{\wtr}}

\put(1800,0){\makebox(0,0)[lb]{\wtr}}
\put(1800,100){\makebox(0,0)[lb]{\wtr}}
\put(1800,200){\makebox(0,0)[lb]{\wtr}}
\put(1800,300){\makebox(0,0)[lb]{\wtr}}
\put(1800,500){\makebox(0,0)[lb]{\wtr}}
\put(1800,600){\makebox(0,0)[lb]{\wtr}}
\put(1800,700){\makebox(0,0)[lb]{\wtr}}
\put(1800,800){\makebox(0,0)[lb]{\wtr}}
\put(1800,900){\makebox(0,0)[lb]{\wtr}}
\put(1800,1000){\makebox(0,0)[lb]{\wtr}}
\put(1800,1100){\makebox(0,0)[lb]{\wtr}}
\put(1800,1200){\makebox(0,0)[lb]{\wtr}}
\put(1800,1300){\makebox(0,0)[lb]{\wtr}}
\put(1800,1400){\makebox(0,0)[lb]{\wtr}}
\put(1800,1500){\makebox(0,0)[lb]{\wtr}}
\put(1800,1600){\makebox(0,0)[lb]{\wtr}}
\put(1800,1700){\makebox(0,0)[lb]{\wtr}}
\put(1800,1800){\makebox(0,0)[lb]{\wtr}}
\put(1800,1900){\makebox(0,0)[lb]{\wtr}}
\put(1800,2000){\makebox(0,0)[lb]{\wtr}}
\put(1800,2200){\makebox(0,0)[lb]{\wtr}}
\put(1800,2300){\makebox(0,0)[lb]{\wtr}}
\put(1800,2400){\makebox(0,0)[lb]{\wtr}}
\put(1800,2600){\makebox(0,0)[lb]{\wtr}}
\put(1800,2700){\makebox(0,0)[lb]{\wtr}}
\put(1800,2800){\makebox(0,0)[lb]{\wtr}}
\put(1800,2900){\makebox(0,0)[lb]{\wtr}}
\put(1800,3000){\makebox(0,0)[lb]{\wtr}}
\put(1800,3100){\makebox(0,0)[lb]{\wtr}}

\put(-100,930){\makebox(0,0)[lb]{\lupa}}
\put(-100,130){\makebox(0,0)[lb]{\lupa}}
\put(-100,1430){\makebox(0,0)[lb]{\lupa}}
\put(340,280){\makebox(0,0)[lb]{\rupa}}
\put(740,380){\makebox(0,0)[lb]{\rupa}}
\put(-60,2580){\makebox(0,0)[lb]{\rupa}}
\put(1140,580){\makebox(0,0)[lb]{\rupa}}
\put(-500,1540){\makebox(0,0)[lb]{\lupa}}
\put(-900,1740){\makebox(0,0)[lb]{\lupa}}
\put(-500,2140){\makebox(0,0)[lb]{\lupa}}
\put(300,530){\makebox(0,0)[lb]{\lupa}}
\put(340,1480){\makebox(0,0)[lb]{\rupa}}

\put(100,-250){\makebox(0,0)[lb]{$\scriptstyle (o,0)$}}
\put(-700,3250){\makebox(0,0)[lb]{$\scriptstyle (x,t)$}}

\put(4000, 200){\ellipse{200}{400}}
\put(4000, 750){\ellipse{200}{400}}
\put(4000, 1300){\ellipse{200}{400}}
\put(4000, 1850){\ellipse{200}{400}}
\put(4000, 2400){\ellipse{200}{400}}
\put(4000, 2950){\ellipse{200}{400}}
\path(4000,400)(4000,550)
\path(4000,950)(4000,1100)
\path(4000,1500)(4000,1650)
\path(4000,2050)(4000,2200)
\path(4000,2600)(4000,2750)
\put(3850,-250){\makebox(0,0)[lb]{$\scriptstyle (o,0)$}}
\put(3850,3250){\makebox(0,0)[lb]{$\scriptstyle (x,t)$}}

}
\end{picture}
}
\end{center}
\caption{\label{fig-2pt1}(a) A configuration for the discretized contact
process.  Open triangles $\triangle$ denote occupied temporal bonds that
are {\it not} connected from $(o,0)$, while closed triangles
$\blacktriangle$ denote occupied temporal bonds that are connected from
$(o,0)$.  The arrows denote occupied spatial bonds, which represent the
spread of the infection to neighbouring sites.  (b) Schematic depiction
of the configuration connecting $(o,0)$ and $(x,t)$ as a ``string of
sausages.''}
\end{figure}
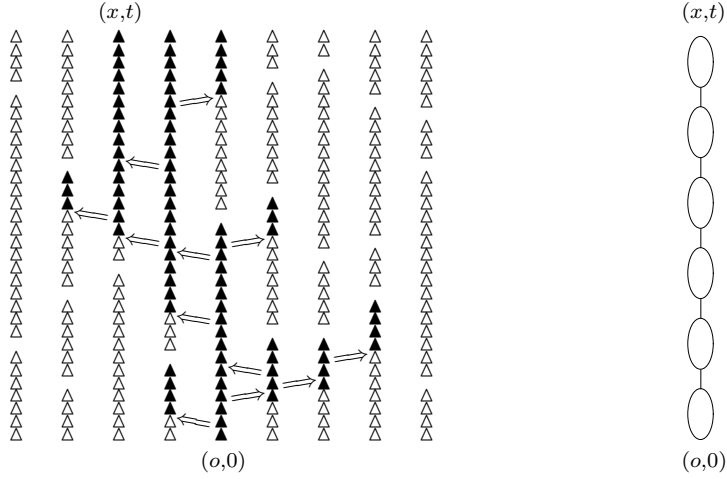

We next explain the expansion for the discretized contact
process in more detail, following the explanation in
\cite[Section~2.1.1]{hs01}.  For the discretized contact process,
we will regard the part of the oriented percolation cluster
connecting $(o,0)$ to $(x,t)$ as a ``string of sausages.''
An example of such a cluster is shown in Figure~\ref{fig-2pt1}.
The difference between oriented percolation and random walk
resides in the fact that for oriented percolation, there can be
multiple paths of occupied bonds connecting $(o,0)$ to $(x,t)$.
However, for $d>4$, each of those paths passes through the same
\emph{pivotal bonds}, which are the essential bonds for the
connection from $(o,0)$ to $(x,t)$. More precisely, a bond is
pivotal for the connection from $(o,0)$ to $(x,t)$ when
$(o,0)\conn (x,t)$ in the possibly modified configuration in
which the bond is made occupied, and $(o,0)$ is not connected to $(x,t)$
in the possibly modified configuration in which the bond is made
vacant (see also Definition \ref{def-1} below). In the strings-and-sausages
picture, the strings are the pivotal bonds, and the sausages are the
parts of the cluster from $(o,0)$ in between the subsequent
pivotal bonds. We expect that there are of the order $t/\vep$
pivotal bonds. For instance, the first black triangle indicates
that $(o,0)$ is connected to $(o,\vep)$, and this bond is pivotal
for the connection from $(o,0)$ to $(x,t)$.

Using this picture, we can think of the oriented percolation two-point
function as a kind of random walk two-point function with a distribution
describing the statistics of the sausages, taking steps in both space and
time.  Due to the nature of the pivotal bonds, each sausage avoids the
backbone from the endpoint of that sausage to $(x,t)$, so that any
connected path between the sausages is via the pivotal bonds between these
sausages.  Therefore, there is a kind of repulsive interaction between the
sausages.  The main part of our proof shows that this interaction is weak
for $d>4$.

Fix $\lambda\geq 0$.  As we will prove in Section~\ref{s:lace}
below, the generalisation of \refeq{2-pt.smp.2} to the
discretized contact process takes the form
    \eq \lbeq{tautx}
    \tau_{t;\vep}^\lamb(x) =
    \ddsum_{s=0}^{t-\varepsilon}
    (\pi_{s;\vep}^\lamb * p_{\varepsilon} * \tau_{t-s-\varepsilon;\vep}^\lamb)(x)
    + \pi_{t;\vep}^\lamb(x)
    \quad\quad (t \in \varepsilon\N),
    \en
where we use the notation $\dsum$ to denote sums over $\vep\Zp$ and
the coefficients $\pi_{t;\vep}^\lamb(x)$ will be defined in
Section~\ref{s:lace}.  In particular, $\pi_{t;\vep}^\lamb(x)$
depends on $\lamb$, is invariant under the lattice symmetries, and
$\pi_{0;\vep}^\lamb(x)=\delta_{o,x}$ and $\pi_{\vep;\vep}^\lamb(x)=0$.
Note that for $t=0,\vep$, we have $\tau_{0;\vep}^\lamb(x)=\delta_{0,x}$
and $\tau_{\vep;\vep}^\lamb(x)=p_\vep(x)$, which is consistent with
\refeq{tautx}.

Together with the initial values $\pi_{0;\vep}^\lamb(x)=\delta_{o,x}$ and
$\pi_{\varepsilon;\vep}^\lamb(x) = 0$, the identity \refeq{tautx} gives
an inductive definition of the sequence $\pi_{t;\vep}^\lamb(x)$
for $t \geq 2\varepsilon$ with $t \in \varepsilon\Z_+$. However, to
analyse the recursion relation \refeq{tautx}, it will be crucial to have a useful
representation for $\pi_{t;\vep}^\lamb(x)$, and this is provided in
Section~\ref{s:lace}. Note that \refeq{2-pt.smp.2} is of the form
\refeq{tautx} with $\pi_{t;\vep}^\lamb(x) = \delta_{o,x}\delta_{0,t}$,
so that we can think of the coefficients $\pi_{t;\vep}^\lamb(x)$
for $t \geq 2\varepsilon$ as quantifying the repulsive interaction
between the sausages in the ``string of sausages'' picture.

Our proof will be based on showing that
$\frac1{\vep^2}\pi_{t;\vep}^\lamb(x)$ for $t\geq2\vep$ is small at
$\lamb=\lambc^{\sss(\vep)}$ if $d>4$ and both $t$ and $L$ are large,
uniformly in $\vep>0$.  Based on this fact, we can rewrite the Fourier
transform of \refeq{tautx} as
\begin{align}
\frac{\hat\tau_{t;\vep}^\lamb(k)-\hat\tau_{t-\vep;\vep}^\lamb(k)}\vep
 =\frac{\hat p_\vep(k)-1}\vep\,\hat\tau_{t;\vep}^\lamb(k)+\vep\ddsum_{s
 =\vep}^{t-\vep}\frac{\hat\pi_{s;\vep}^\lamb(k)}{\vep^2}\,\hat p_\vep(k)
 \,\hat\tau_{t-s-\vep;\vep}^\lamb(k)+\frac{\hat\pi_{t;\vep}^\lamb(k)}\vep.
\end{align}
Assuming convergence of $\frac1{\vep^2}\hat\pi_{s;\vep}^\lamb(k)$ to
$\hat\pi_s^\lamb(k)$, which will be shown in Section~\ref{sec-outlinedisc},
we obtain \refeq{tautxcont}.  Therefore, \refeq{tautxcont} is regarded as
a small perturbation of \refeq{SRW} when $d>4$ and $L\gg1$, and this will
imply the central limit theorem for the critical two-point function.

Now we briefly explain the expansion coefficients $\pi_{t;\vep}^\lamb(x)$.
In Section~\ref{s:lace}, we will obtain the expression
\begin{align}\lbeq{pitN}
    \pi_{t;\vep}^\lamb(x) =\sum_{N=0}^{\infty} (-1)^N
    \pi_{t;\vep}^{\smallsup{N}}(x),
\end{align}
where we suppress the dependence of $\pi_{t;\vep}^{\sss(N)}(x)$ on
$\lamb$.  The idea behind the proof of \refeq{pitN} is the following.
Let
\begin{align}\lbeq{pi0def}
\pi_{t;\vep}^{\sss(0)}(x)=\mP_\vep^\lamb((o,0)\db(x,t))
\end{align}
denote the contribution to $\tau^\lamb_{t;\vep}(x)$ from configurations
in which there are no pivotal bonds, so that
\begin{align}\lbeq{tausplit1}
\tau^\lamb_{t;\vep}(x)=\pi_{t;\vep}^{\sss(0)}(x)+\sum_b\mP_\vep^\lamb
 \big(b\text{ first occupied and pivotal bond for }(o,0)\conn(x,t)\big),
\end{align}
where the sum over $b$ is over bonds of the form $b=((u,s),(v,s+\vep))$.
We write $\bb=(u,s)$ for the starting point of the bond $b$ and
$\tb=(v,s+\vep)$ for its endpoint.  Then, the probability on the
right-hand side of \refeq{tausplit1} equals
\begin{align}
\mP_\vep^\lamb\big((o,0)\db\bb,~b\text{ occupied},~\tb\conn(x,t),~b
 \text{ pivotal for }(o,0)\conn(x,t)\big).
\end{align}
We ignore the intersection with the event that $b$ is pivotal for
$(o,0)\conn(x,t)$, and obtain using the Markov property that
\begin{align}\lbeq{tausplit2}
\tau^\lamb_{t;\vep}(x)=\pi_{t;\vep}^{\sss(0)}(x)+\ddsum_{s=0}^{t-\vep}
 \sum_{u,v\in\Zd}\pi_{s;\vep}^{\sss(0)}(u)\,p_\vep(v-u)\,\tau_{t-s-\vep;
 \vep}(x-v)-R^{\sss(0)}_{t;\vep}(x),
\end{align}
where
\begin{align}\lbeq{R0def}
R^{\sss(0)}_{t;\vep}(x)=\sum_b\mP_\vep^\lamb\big((o,0)\db\bb,~b
 \text{ occupied},~\tb\conn(x,t),~b\text{ not pivotal for }(o,0)
 \conn(x,t)\big).
\end{align}
We will investigate the error term $R^{\smallsup{0}}_{t;\vep}(x)$ further,
again using inclusion-exclusion, by investigating the first pivotal bond
after $\tb$ to arrive at \refeq{pitN}. The term
$\pi^{\smallsup{1}}_{t;\vep}(x)$
is the contribution to $R^{\smallsup{0}}_{t;\vep}(x)$ where such a pivotal
does not exist. Thus, in $\pi^{\smallsup{0}}_{t;\vep}(x)$ for $t\geq \vep$ and
in $\pi^{\smallsup{1}}_{t;\vep}(x)$ for all $t\geq 0$, there is at least
one loop, which, for $L$ large, should yield a small correction only. In
\refeq{pitN}, the contributions from $N\geq 2$ have at least two loops
and are thus again smaller, even though all $N\geq 0$ give essential
contributions to $\pi_{t;\vep}^\lambda(x)$ in \refeq{pitN}.

There are three ways to obtain the lace expansion in \refeq{tautx}
for oriented percolation models. We use the expansion by
Sakai \cite{s00,s01}, as described in \refeq{pi0def}--\refeq{R0def}
above, based on inclusion-exclusion together with the Markov
property for oriented percolation. For unoriented percolation,
Hara and Slade \cite{HS90a} developed an expression for
$\pi_{t;\vep}^{\lamb}(x)$ in terms of sums of nested
expectations, by repeated use of inclusion-exclusion
and using the independence of percolation. This expansion, and
its generalizations to the higher-point functions,
was used in \cite{hs01} to investigate the oriented percolation
$r$-point functions. The original expansion in \cite{HS90a} was for
unoriented percolation, and does not make use of the Markov property.
Nguyen and Yang \cite{ny93, ny95} derived an alternate expression for
$\pi_{t;\vep}^\smallsup{N}(x)$ by adapting the lace expansion of Brydges
and Spencer \cite{BS85} for weakly self-avoiding walk.
In the graphical representation of the Brydges-Spencer expansion,
laces arise which give the ``lace expansion'' its name. Even though in many
of the lace expansions for percolation type models, such as oriented and
unoriented percolation, no laces appear, the name has stuck for historical
reasons.

It is not so hard to see that the Nguyen-Yang expansion is equivalent to
the above expansion using inclusion-exclusion, just as for self-avoiding
walks \cite{ms93}.  Since we find the Sakai expansion simpler, especially
when dealing with the continuum limit, we prefer the Sakai expansion to the
Nguyen-Yang expansion.  In \cite{hs01}, the Hara-Slade expansion was used
to obtain \refeq{pitN} with a {\it different} expression for
$\pi_{t;\vep}^{\sss(N)}(x)$.  In either expansion,
$\pi_{t;\vep}^{\sss(N)}(x)$ is nonnegative for all $t,x,N$, and can be
represented in terms of Feynman-type diagrams.  The Feynman diagrams are
similar for the three expansions and obey similar estimates, even though
the expansion used in this paper produces the simplest diagrams.

In \cite{hs01}, the Nguyen-Yang expansion was also used to deal with the
derivative of the lace expansion coefficients with respect to the
percolation parameter $p$.  In this paper, we use the inclusion-exclusion
expansion also for the derivative of the expansion coefficients with
respect to $\lamb$, rather than on two different expansions as in
\cite{hs01}.

We now comment on the relative merits of the Sakai and the Hara-Slade
expansion.  Clearly, the Hara-Slade expansion is more general, as it also
applies to unoriented percolation.  On the other hand, the Sakai expansion
is somewhat simpler to use, and the bounding diagrams on the arising
Feynman diagrams are simpler.  Finally, the resulting expressions for
$\pi_{t;\vep}^{\sss(N)}(x)$ in the Sakai expansion allow for a continuum
limit, where it is not clear to us how to perform this limit using the
Hara-Slade expansion coefficients.

In \cite{hsa04}, we will adapt the expansion in Section~\ref{s:lace} to
deal with the discretized contact process and oriented percolation
higher-point functions.  For this, we will need ingredients from the
Hara-Slade expansion to compare occupied paths living on a {\it common}
time interval, with independent paths.  This independence does not follow
from the Markov property, and therefore the Hara-Slade expansion, which
does not require the Markov property, will be crucial.  The ``decoupling''
of disjoint paths is crucial in the derivation of the lace expansion for
the higher point functions, and explains the importance of the Hara-Slade
expansion for oriented percolation and the contact process.

To complete this discussion, we note that an alternative route to
the contact process results is via \refeq{tautxcont}. In \cite{BR01},
an approach using a Banach fixed point theorem was used to prove
asymptotics of the two-point function for weakly self-avoiding walk.
The crucial observation is that a lace expansion equation such as
\refeq{tautxcont} can be viewed as a fixed point equation of a certain
operator on sequence spaces. By proving properties of this operator,
Bolthausen and Ritzman were able to deduce properties of the fixed point
sequence, and thus of the weakly self-avoiding walk two-point function.
It would be interesting to investigate whether such an approach may
be used on \refeq{tautxcont} as well.

\subsection{Bounds on the lace expansion}\label{sec-ble}
In order to prove the statements in Proposition~\ref{thm-disc}, we will use
induction in $n$, where $t=n\vep\in\vep\Zp$.  The lace expansion equation
in \refeq{tautx} forms the main ingredient for this induction in time.  We
will explain the inductive method in more detail below.  To advance the
induction hypotheses, we clearly need to have certain bounds on the lace
expansion coefficients.  The form of those bounds will be explained now.
The statement of the bounds involve the small parameter
\begin{align}\lbeq{betadef}
\beta=L^{-d}.
\end{align}
We will use the following set of bounds:
\begin{align}\lbeq{fbd}
|\hat\tau_{s;\vep}(0)|\leq K,&&
|\nabla^2\hat\tau_{s;\vep}(0)|\leq K\sigma^2s,&&
\|\wD^2\,\hat{\tau}_{s;\vep}\|_1\leq\frac{K\beta}{(1+s)^{d/2}},
\end{align}
where we write
$\|\hat f\|_1=\int_{[-\pi,\pi]^d}\frac{d^dk}{(2\pi)^d}\,|\hat f(k)|$
for a function $\hat f:[-\pi,\pi]^d\mapsto\Cbold$.  The bounds on the
lace expansion consist of the following estimates, which will be proved
in Section~\ref{s:pibds}.

\begin{prop}[Bounds on the lace expansion for $d>4$]\label{lem-Pibd.smp}
Assume \refeq{fbd} for some $\lamb_0$ and all $s\leq t$.  Then, there are
$\beta_0=\beta_0(d,K)>0$ and $C=C(d,K)<\infty$ (both independent of
$\vep,L$) such that, for $\lamb\leq\lamb_0$, $\beta<\beta_0$, $s\in\vep\Zp$
with $2\vep\leq s\leq t+\vep$, $q=0,2,4$ and $\Delta'\in[0,1\wedge\Delta]$,
and uniformly in $\vep\in(0,1]$,
\begin{gather}
\sum_{x\in\Zd}|x|^q\,|\pi_{s;\vep}^{\lambda}(x)|\leq\frac{\vep^2C\sigma^q
 \beta}{(1+s)^{(d-q)/2}},\lbeq{pimom}\\
\Big|\wpi^{\lambda}_{s;\vep}(k)-\wpi^{\lambda}_{s;\vep}(0)-\frac{a(k)}
 {\sigma^2}\,\nabla^2\wpi^{\lambda}_{s;\vep}(0)\Big|\leq\frac{\vep^2C
 \beta\,a(k)^{1+\Delta'}}{(1+s)^{(d-2)/2-\Delta'}},\lbeq{pidiff}\\
|\partial_\lamb\wpi^{\lambda}_{s;\vep}(0)|\leq\frac{\vep^2C\beta}
 {(1+s)^{(d-2)/2}}.\lbeq{pider}
\end{gather}
\end{prop}
The main content of Proposition \ref{lem-Pibd.smp} is that
the bounds on $\hat{\tau}_{s;\vep}$ for $s\leq t$ in \refeq{fbd}
imply bounds on $\hat\pi_{s;\vep}$ for all $s\leq t+\vep$. This fact allows
us to use the bounds on $\hat\pi_{s;\vep}$ for all arising $s$ in
\refeq{tautx} in order to advance the appropriate induction hypotheses.
Of course, in order to complete the inductive argument, we need that
the induction statements imply the bounds in \refeq{fbd}.

The proof of Proposition~\ref{lem-Pibd.smp} is deferred to
Section~\ref{s:pibds}. Proposition~\ref{lem-Pibd.smp} is
probably false in dimensions $d\leq 4$. However, when the range
increases with $T$ as in Theorem \ref{thm:2pt-lowdim}, we
are still able to obtain the necessary bounds. In the statement of
the bounds, we recall that $L_{\sT}$ is given in \refeq{Lt-def}.

\begin{prop}[Bounds on the lace expansion for $d\leq4$]\label{lem-Pibd.smp2}
Let $\alpha>0$ in \refeq{alphadef}.  Assume \refeq{fbd}, with $\beta$
replaced by $\beta_{\sT}=L_{\sT}^{-d}$ and $\sigma^2$ by $\sigma_{\sT}^2$,
for some $\lamb_0$ and all $s\leq t$.  Then, there are $L_0=L_0(d,K)<\infty$
(independent of $\vep$) and $C=C(d,K)<\infty$ (independent of $\vep,L$)
such that, for $\lamb\leq\lamb_0$, $L_1\geq L_0$, $s\in\vep\Zp$ with
$2\vep\leq s\leq t+\vep$, $q=0,2,4$ and $\Delta'\in[0,1\wedge\Delta]$,
the bounds in \refeq{pimom}--\refeq{pider} hold for $t\leq T\log T$,
with $\beta$ replaced by $\beta_{\sT}=L_{\sT}^{-d}$ and $\sigma^2$ by
$\sigma_{\sT}^2$.
\end{prop}

The main point in Propositions~\ref{lem-Pibd.smp}--\ref{lem-Pibd.smp2} is
the fact that we need to extract two factors of $\vep$.  One can see that
such factors must be present by investigating, e.g.,
$\pi^{\sss(0)}_{t;\vep}(x)$, which is the probability that $(o,0)$ is
doubly connected to $(x,t)$.  When $t>0$, there must be at least two
spatial bonds, one emanating from $(o,0)$ and one pointing into $(x,t)$.
By \refeq{bprob}, these two spatial bonds give rise to two powers of
$\vep$.  The proof for $N\geq 1$ then follows by induction in $N$.

\subsection{Implementation of the inductive method}\label{sec-indimp}
Our analysis of \refeq{tautx} begins by taking its Fourier
transform, which gives the recursion relation
\begin{align}\lbeq{tauk}
\hat\tau_{t;\vep}^\lamb(k)=\ddsum_{s=0}^{t-\vep}\hat\pi_{s;\vep}^\lamb
 (k)\,\hat p_\vep(k)\,\hat\tau_{t-s-\vep;\vep}^\lamb(k)+\hat\pi_{t;
 \vep}^\lamb(k)\quad\quad(t\in\vep\N).
\end{align}
As already explained in Section~\ref{sec-ble}, it is possible to estimate
$\hat\pi_{s;\vep}^\lamb(k)$, for all $s\leq t$, in terms of $\|\tau_{s;\vep}
^\lamb\|_1\equiv\sum_{x\in\Zd}\tau_{s;\vep}^\lamb(x)=\hat\tau_{s;\vep}^\lamb
(0)$ and $\|\tau_{s;\vep}^\lamb\|_\infty\leq\|\hat\tau_{s;\vep}^\lamb\|_1$
with $s\leq t-\vep$.  Therefore, the right-hand side of \refeq{tauk}
explicitly involves $\hat\tau_{s;\vep}^\lamb(k)$ only for $s\leq t-\vep$.
This opens up the possibility of an inductive analysis of \refeq{tauk}.
A general approach to this type of inductive analysis is given in
\cite{hs02}.  However, here we will need the uniformity in the variable
$\vep$, and therefore we will state a version of the induction in
Section~\ref{s:ind} that is adapted to the uniformity in $\vep$ and thus
the continuum limit.  The advancement of the induction hypotheses is
deferred to Appendix~\ref{s:adv-ind}.

Moreover, we will show that the critical point is given implicitly by the
equation
\begin{align}\lbeq{lambdacdef}
\lambc^{\sss(\vep)}=1-\frac1\vep\ddsum_{s=2\vep}^\infty\hat\pi^{\lambc^{
 \smallsup{\vep}}}_{s;\vep}\!(0)\;\hat p_\vep^{\lambc^{(\vep)}}\!(0),
\end{align}
and that the constants $A^{\sss(\vep)}$ and $v^{\sss(\vep)}$ of
Proposition~\ref{thm-disc} are given by
\begin{align}\lbeq{Avdef}
A^{\sss(\vep)}=\frac{\dpst1+\ddsum_{s=2\vep}^\infty\hat\pi_{s;\vep}^{\lambc
 ^{(\vep)}}\!(0)}{\dpst1+\frac1\vep\ddsum_{s=2\vep}^\infty s\;\hat\pi^{
 \lambc^{\smallsup{\vep}}}_{s;\vep}\!(0)\;\hat p_\vep^{\lambc^{(\vep)}}\!
 (0)},\qquad
v^{\sss(\vep)}=\frac{\dpst\lambc^{\sss(\vep)}-\frac1{\sigma^2\vep}\ddsum_{s
 =2\vep}^\infty\nabla^2\big[\hat\pi^{\lambc^{(\vep)}}_{s;\vep}\!(k)\;\hat
 p_\vep^{\lambc^{(\vep)}}\!(k)\big]_{k=0}}{\dpst1+\frac1\vep\ddsum_{s=2\vep}
 ^\infty s\;\hat\pi^{\lambc^{\smallsup{\vep}}}_{s;\vep}\!(0)\;\hat p_\vep
 ^{\lambc^{(\vep)}}\!(0)},
\end{align}
%
where we have added an argument $\lambc^{\sss(\vep)}$ to emphasize that
$\lamb$ is critical for the evaluation of $\pi_{t;\vep}^\lamb$ on the
right-hand sides.  Convergence of the series on the right-hand sides,
for $d>4$, follows from Proposition~\ref{lem-Pibd.smp}.  For oriented
percolation, i.e., for $\vep=1$, these equations agree with
\cite[(2.11-2.13)]{hs01}.

The result of induction is summarized in the following proposition:
\begin{prop}[Induction]\label{prop-ind}
If Proposition~\ref{lem-Pibd.smp} holds, then \refeq{fbd} holds for
$s\leq t+\vep$.  Therefore, \refeq{fbd} holds for all $s\geq0$ and
\refeq{pimom}--\refeq{pider} hold for all $s\geq2\vep$.  Moreover,
the statements in Proposition~\ref{thm-disc} follow, with the error
terms uniform in $\vep\in(0,1]$.
\end{prop}

There is also a low-dimensional version of Proposition~\ref{prop-ind},
but we refrain from stating it.

\subsection{Continuum limit}\label{sec-outlinedisc}
In this section we state the result necessary to complete the proof of
Theorems~\ref{thm:2pt}--\ref{thm:2pt-lowdim} from
Propositions~\ref{thm-disc}--\ref{thm-disc2}.  In particular,
from now onwards, we specialize to the contact process.

\begin{prop}[Continuum limit]\label{prop-disc}
Suppose that $\lamb^{\sss(\vep)}\to\lamb$ and
$\lamb^{\sss(\vep)}\leq\lambc^{\sss(\vep)}$ for $\vep$ sufficiently small.
Then, for every $t>0$ and $x\in\Zd$, there is a $\pi_t^\lamb(x)$ such that
\begin{align}\lbeq{piconvvep}
\lim_{\vep\daw0}\frac1{\vep^2}\pi_{t;\vep}^{\lamb^{(\vep)}}(x)
 =\pi_t^\lamb(x),&&
\lim_{\vep\daw0}\frac1{\vep^2}[\partial_\ell\pi_{t;\vep}^\ell
 (x)]_{\ell=\lamb^{(\vep)}}=\partial_\lamb\pi_t^\lamb(x).
\end{align}
Consequently, for $\lamb\leq\lambc$ and $q=0,2,4$,
\begin{align}\lbeq{boundspiCP}
\sum_{x\in\Zd}|x|^q\pi_t^\lamb(x)\leq\frac{C\beta}{(1+t)^{(d-q)/2}},&&
\sum_{x\in\Zd}\partial_\lamb\pi_t^\lamb(x)\leq\frac{C\beta}{(1+t)^{(d-2)/2}},
\end{align}
and there exist $A=1+O(L^{-d})$ and $v=1+O(L^{-d})$ such that
\begin{align}\lbeq{convAv}
\lim_{\vep\daw0}A^{\sss(\vep)}=A,&& \lim_{\vep\daw0}v^{\sss(\vep)}=v.
\end{align}
Furthermore, $\partial_\lamb\pi_t^\lamb(x)$ is continuous in $\lamb$.
\end{prop}
In Proposition~\ref{lem-Pibd.smp}, the right-hand sides of
\refeq{pimom}--\refeq{pider} are proportional to $\vep^2$.
The main point in the proof of Proposition~\ref{prop-disc}
is that the lace expansion coefficients, scaled by $\vep^{-2}$,
converge as $\vep\daw0$, using the weak convergence of
$\mP_\vep^\lamb$ to $\mP^\lamb$ \cite[Proposition 2.7]{bg91}.

In Section~\ref{s:continuum}, we will show that
$\frac1{\vep^2}\pi_{t;\vep}^\lamb(x)$ and
$\frac1{\vep^2}\partial_\lamb\pi_{t;\vep}^\lamb(x)$ both converge
pointwise.  We now show that this implies that the limit of
$\frac1{\vep^2}\partial_\lamb\pi_{t;\vep}^\lamb(x)$ equals
$\partial_\lamb\pi_t^\lamb(x)$.
To see this, we use
\begin{align}
\frac1{\vep^2}\pi_{t;\vep}^\lamb(x)=\int_0^\lamb d\lamb'
 ~\frac1{\vep^2}\partial_{\lamb'}\pi_{t;\vep}^{\lamb'}(x).
\end{align}
where we use $\frac1{\vep^2}\pi_{t;\vep}^0(x)=0$ for $t>0$.  By the assumed
pointwise convergence, the left-hand side converges to $\pi_t^\lamb(x)$,
while the right-hand side converges to the integral of the limit of
$\frac1{\vep^2}\partial_{\lamb'}\pi_{t;\vep}^{\lamb'}(x)$, denoted
$f_t^{\lamb'}(x)$ for now, using the dominated convergence theorem.
Therefore, for any $\lamb\leq\lambc$,
\begin{align}
\pi_t^\lamb(x)=\int_0^\lamb d\lamb'~f_t^{\lamb'}(x),
\end{align}
which indeed implies that $f_t^\lamb(x)=\partial_\lamb\pi_t^\lamb(x)$.

\begin{proof}[Proof of Theorems~\ref{thm:2pt}--\ref{thm:2pt-lowdim}
assuming Propositions~\ref{thm-disc}--\ref{thm-disc2} and \ref{prop-disc}]
We only prove Theorem~\ref{thm:2pt}, since the proof of
Theorem~\ref{thm:2pt-lowdim} is identical.  By \cite[Proposition 2.7]{bg91},
we have that, for every $(x,t)$ and $\lamb>0$,
\begin{align}
\lim_{\vep\daw0}\tau_{t;\vep}^\lamb(x)=\tau_t^\lamb(x).
\end{align}
Since 
$\tau_t^\lamb(x)$ is continuous in $\lamb$ (see e.g.,
\cite[pp.38--39]{Ligg99}), we also obtain
$\lim_{\vep\daw0}\tau_{t;\vep}^{\lamb^{(\vep)}}(x)=\tau_t^\lamb(x)$ for
any $\lamb^{\sss(\vep)}\to\lamb$.  Since $\lambc^{\sss(\vep)}\to\lambc$
\cite[Section~2.1]{s01}, $\tau_{t;\vep}^{\lambc^{\sss(\vep)}}(x)$ also
converges to $\tau_t^{\lambc}(x)$.  Using the uniformity in $\vep$ of
the upper and lower bounds in \refeq{tausupvep}, we obtain \refeq{tausup}.


Next, we prove $\lim_{\vep\daw0}\hat\tau_{t;\vep}^{\lambc^{(\vep)}}(k)=
\hat\tau_t^{\lambc}(k)$ for every $k\in[-\pi,\pi]^d$ and $t\geq0$.  Note
that the Fourier transform involves a sum over $\Zd$, such as
\begin{align}\lbeq{tauhatsum}
\hat\tau_{t;\vep}^{\lambc^{(\vep)}}(k)=\sum_{x\in\Zd}
 \tau_{t;\vep}^{\lambc^{(\vep)}}\!(x)\;e^{ik\cdot x}.
\end{align}
To use the pointwise convergence of $\tau_{t;\vep}^{\lambc^{(\vep)}}(x)$,
we first show that the sum over $x\in\Zd$ in \refeq{tauhatsum} can be
approximated by a finite sum.  To see this, we note that
\begin{align}\lbeq{tauhatbd}
\tau_{t;\vep}^\lamb(x)\leq p_\vep^{*t/\vep}(x)=\sum_{n=0}^{t/\vep}
 \binom{t/\vep}{n}(1-\vep)^{t/\vep-n}(\lamb\vep)^n\,D^{*n}(x).
\end{align}
For any fixed $t$, we can choose $\delta_{\sss R}\geq0$, which is
$\vep$-independent and decays to zero as $R\uaw\infty$, such that
\begin{align}\lbeq{tauhatbd-appl}
\sum_{x\in\Zd:\|x\|_\infty>R}\tau_{t;\vep}^\lamb(x)\leq\delta_{\sss R}.
\end{align}
Therefore, the same holds for $\tau_t^\lamb(x)$, and hence we can
approximate both $\hat\tau_{t;\vep}^{\lambc^{(\vep)}}(k)$ and
$\hat\tau_t^{\lambc}(k)$ by sums over $x\in\Zd$ with $\|x\|_\infty\leq R$,
in which we use the pointwise convergence of
$\tau_{t;\vep}^{\lambc^{(\vep)}}(x)$.  Taking $R\uaw\infty$, we obtain $\hat
\tau_t^{\lambc}(k)=\lim_{\vep\daw0}\hat\tau_{t;\vep}^{\lambc^{(\vep)}}(k)$.

Using the above, we obtain
\begin{align}
\hat\tau_t^{\lambc}\big(\tfrac{k}{\sqrt{v\sigma^2t}}\big)&=\lim_{\vep\daw0}
 \hat\tau_{t;\vep}^{\lambc^{(\vep)}}\big(\tfrac{k}{\sqrt{v\sigma^2t}}\big)
=\lim_{\vep\daw0}\hat\tau_{t;\vep}^{\lambc^{(\vep)}}\big(\tfrac{\sqrt{v^{
 (\vep)}}}{\sqrt{v}}\tfrac{k}{\sqrt{v^{(\vep)}\sigma^2t}}\big)\nn\\
&=\lim_{\vep\daw0}A^{\sss(\vep)}\,e^{-\frac{v^{(\vep)}}v\frac{|k|^2}{2d}}\,
 \big[1+O\big(\tfrac{v^{(\vep)}}v|k|^2(1+t)^{-\delta}\big)+O((1+t)^{-(d-4)
 /2})\big]\nn\\
&=A\,e^{-\frac{|k|^2}{2d}}\,\big[1+O(|k|^2(1+t)^{-\delta})+O((1+t)^{-(d-4)
 /2})\big],
\end{align}
which proves \refeq{tauasy}. Similar argument can be used for \refeq{taugyr}.
\end{proof}

\begin{proof}[Proof of \refeq{gammasharp} assuming \refeq{tautxcont} and
Proposition~\ref{prop-disc}]
We now prove that, in the current setting,
$\chi(\lamb)=\int_0^\infty dt\,\hat\tau_t^\lamb(0)$ satisfies the precise
asymptotics in \refeq{gammasharp}, assuming \refeq{tautxcont} and
Proposition~\ref{prop-disc}.

Let $\lamb<\lambc$.  Since $\hat\tau_0^\lamb(0)=1$ and
$\hat\tau_\infty^\lamb(0)=0$, using \refeq{tautxcont} we obtain
\begin{align}
-1=\int_0^\infty dt~\partial_t\hat\tau_t^\lamb(0)&=\int_0^\infty dt~
 \bigg[(\lamb-1)\,\hat\tau_t^\lamb(0)+\int_0^t ds~\hat\pi_s^\lamb(0)\,
 \hat\tau_{t-s}^\lamb(0)\bigg]\nn\\
&=\bigg[\lamb-1+\int_0^\infty ds~\hat\pi_s^\lamb(0)\bigg]\int_0^\infty
 dt~\hat\tau_t^\lamb(0),
\end{align}
so that
\begin{align}\lbeq{tauid}
\chi(\lamb)=\bigg[1-\lamb-\int_0^\infty ds~\hat\pi_s^\lamb(0)\bigg]^{-1}.
\end{align}
By \refeq{lambdacdef} and Proposition~\ref{prop-disc}, $\lambc$ must satisfy
\begin{align}\lbeq{lambcont}
\lambc=1-\int_0^\infty ds~\hat\pi_s^{\lambc}(0),
\end{align}
so that we can rewrite \refeq{tauid} as
\begin{align}\lbeq{tauid2}
\chi(\lamb)=[f(\lambc)-f(\lamb)]^{-1},
\end{align}
where $f(\lamb)=\lamb+\int_0^\infty ds\,\hat\pi_s^\lamb(0)$, since, by
\refeq{lambcont}, $f(\lambc)=1$.  By the mean-value theorem, together
with the fact that $|\partial_\lamb\hat\pi_s^\lamb(0)|$ is integrable
with respect to $s>0$, there is a $\lamb_*\in(\lamb,\lambc)$ such that
\begin{align}\lbeq{tauid3}
\chi(\lamb)=[(\lambc-\lamb)\,f'(\lamb_*)]^{-1}.
\end{align}
By the continuity in $\lamb$ of $\partial_\lamb\pi_t^\lamb(x)$ and its
summability in $(x,t)\in\Zd\times\mR_+$ for $\lamb\leq\lambc$ due to
\refeq{boundspiCP},
$f'(\lamb)=1+\int_0^\infty ds\;\partial_\lamb\hat\pi_t^\lamb(0)$ is also
continuous in $\lamb\leq\lambc$.  Therefore, we obtain \refeq{gammasharp}
with $C=f'(\lambc)^{-1}$.

Finally, we note that the above proof, where the integral is replaced with
a sum over $n\in\Zp$, also shows that the stronger version of $\gamma=1$
holds for oriented percolation.
\end{proof}

The proofs of Theorems~\ref{thm:2pt}--\ref{thm:2pt-lowdim} are now reduced
to the proof of Propositions~\ref{thm-disc}--\ref{thm-disc2} and
\ref{prop-disc}.  Proposition~\ref{prop-disc} will be proved in
Section~\ref{s:continuum}.  The proof of
Propositions~\ref{thm-disc}--\ref{thm-disc2} is reduced to
Propositions~\ref{lem-Pibd.smp}--\ref{prop-ind}, which will be proved in
Sections~\ref{s:pibds}--\ref{s:ind}.  The advancement of the induction
hypotheses is deferred to Appendix~\ref{s:adv-ind}.  We start in
Section~\ref{s:lace} by deriving the lace expansion \refeq{tautx}.


\newcommand{\piv}{{\tt piv}}
\newcommand{\Lb}[2]{\cL_{#1}^{\sss(#2)}}

\section{Lace expansion}\label{s:lace}
In this section, we derive the lace expansion in \refeq{tautx}.
The same type of recursion relation was used for {\it discrete}
models, such as (weakly) self-avoiding walk in $\Zd$
\cite{BS85,HS92b,hhs98,hs02,hs03,Slad87,Slad88,Slad91} and
oriented percolation in $\Zd\times\Zp$ \cite{hs02,hs01,ny93,ny95}.

%

From now on, we will suppress the dependence on $\vep$ and $\lambda$ when no confusion
can arise, and write, e.g., $\pi_t(x)=\pi_{t;\vep}^\lambda(x)$.  In
Section~\ref{ss:exp-two}, we obtain \refeq{tau-2exp},
which is equivalent to the recursion relation in \refeq{tautx}, and
the expression \refeq{pi-def} for $\pi_t(x)$.
In Section~\ref{ss:rep-der},  we obtain the expressions
\refeq{pider-exp}--\refeq{Pider-exp} for $\partial_\lamb\pi_t(x)$.


\subsection{Expansion for the two-point function}\label{ss:exp-two}
In this section, we derive the expansion \refeq{tau-2exp}.
We will also write $\Lambda=\Zd\times\vep\Zp$,
and use bold letters $\ovec, \xvec,\dots$ to represent elements in $\Lambda$,
such as $\ovec=(o,0)$ and $\xvec=(x,t)$, and
write $\tau(\xvec)=\tau_t(x)$, $\pi^{\sss(N)}(\xvec)=\pi_t^{\sss(N)}(x)$,
and so on.

We recall that the two-point function is defined by
    \begin{align}\lbeq{2pt-remind}
    \tau(\xvec)=\mP(\ovec\conn\xvec).
    \end{align}

Before starting with the expansion, we introduce some definition:

\begin{defn}
\label{def-1}
\begin{enumerate}
\item[(i)]
For a bond $b=(\uvec,\vvec)$, we write $\bb=\uvec$ and $\tb=\vvec$.  We
write $b\conn\xvec$ for the event that $b$ is occupied and $\tb\conn\xvec$.
\item[(ii)]
Given a configuration, we say that $\vvec$ is {\em doubly connected to}
$\xvec$, and we write $\vvec\dbc\xvec$, if there are at least two
bond-disjoint paths from $\vvec$ to $\xvec$ consisting of occupied bonds.
By convention, we say that $\xvec\dbc\xvec$ for all $\xvec$.
\item[(iii)]
A bond is said to be {\em pivotal} for $\vvec\conn\xvec$ if
$\vvec\conn\xvec$ in the possibly modified configuration in which that
bond is made occupied, whereas $\vvec$ is not connected to $\xvec$ in
the possibly modified configuration in which that bond is made vacant.
\end{enumerate}
\end{defn}

We split, depending on whether there is a pivotal bond for
$\ovec\conn\xvec$, to obtain
\begin{align}\lbeq{exptau-1}
\tau(\xvec)=\mP(\ovec\dbc\xvec)+\sum_b\mP(\ovec\dbc\bb,~b
 \text{ occupied \& pivotal for }\ovec\conn\xvec).
\end{align}
We denote
\begin{align}\lbeq{pi0-def}
\pi^{\sss(0)}(\xvec)=\mP(\ovec\dbc\xvec),
\end{align}
so that we can rewrite \refeq{exptau-1} as
\begin{align}\lbeq{exptau-2}
\tau(\xvec)=\pi^{\sss(0)}(\xvec)+\sum_b\mP(\ovec\dbc\bb,~b\conn\xvec,~
 b\text{ pivotal for }\ovec\conn\xvec).
\end{align}
Define
\begin{align}\lbeq{R0-def}
R^{\sss(0)}(\xvec)=\sum_b\mP(\ovec\dbc\bb,~b\conn\xvec,~
 b\text{ not pivotal for }\ovec\conn\xvec),
\end{align}
then, by inclusion-exclusion on the event that $b$ is pivotal for
$\ovec\conn\xvec$, we arrive at
\begin{align}\lbeq{exptau-3}
\tau(\xvec)=\pi^{\sss(0)}(\xvec)+\sum_b\mP(\ovec\dbc\bb,~b\conn\xvec)
 -R^{\sss(0)}(\xvec).
\end{align}
The event $\ovec\dbc\bb$ only depends on bonds with time variables less
than or equal to the one of $\bb$, while the event $b\conn\xvec$ only
depends on bonds with time variables larger than or equal to the one of
$\bb$.  Therefore, by the Markov property, we obtain
\begin{align}\lbeq{MPOP}
\mP(\ovec\dbc\bb,~b\conn\xvec)=\mP(\ovec\dbc\bb)\,\mP(b\text{ occupied})
 \,\mP(\tb\conn\xvec)=\pi^{\sss(0)}(\bb)\,p(b)\,\tau(\xvec-\tb),
\end{align}
where we abuse notation to write
\begin{align}
p(b)=p(\tb-\bb).
\end{align}
Therefore, we arrive at
\begin{align}\lbeq{exptau-4}
\tau(\xvec)=\pi^{\sss(0)}(\xvec)+(\pi^{\sss(0)}\sstar p\sstar\tau)(\xvec)
 -R^{\sss(0)}(\xvec),
\end{align}
where we use ``$\sstar$'' to denote convolution in $\Lambda$, i.e.,
\begin{align}
(f\sstar g)(\xvec)=\sum_{\yvec\in\Lambda}f(\yvec)\,g(\xvec-\yvec).
\end{align}
This completes the first step of the expansion, and we are left to
investigate $R^{\sss(0)}(\xvec)$.  For this, we need some further notation.

\begin{defn}
\begin{enumerate}
\item[(i)]
Given a configuration and $\xvec\in\Lambda$, we define $\bC(\xvec)$
to be the set of sites to which $\xvec$ is connected, i.e.,
$\bC(\xvec)=\{\yvec\in\Lambda:\xvec\conn\yvec\}$.  Given a bond $b$,
we also define $\tilde\bC^b(\xvec)$ to be the set of sites to which
$\xvec$ is connected in the (possibly modified) configuration in which
$b$ is made vacant.
\item[(ii)]
Given a site set $\bC$, we say that $\vvec$ is connected to $\xvec$
{\it through} $\bC$, if every occupied path connecting $\vvec$ to $\xvec$
has at least one bond with an endpoint in $\bC$.  This event is written
as $\vvec\ct{\bC}\xvec$.  Similarly, we write $\{b\ct{\bC}\xvec\}=\{b$
occupied$\}\cap\{\tb\ct{\bC}\xvec\}$.
\end{enumerate}
\end{defn}

We then note that
\begin{align}
\{\vvec\conn\bb,~b\conn\xvec,~b\text{ not pivotal for }\vvec\conn\xvec\}
 =\big\{\vvec\conn\bb,~b\ctx{\tilde\bC^b(\vvec)}\xvec\big\}.
\end{align}
Therefore,
\begin{align}\lbeq{R0-rew}
R^{\sss(0)}(\xvec)=\sum_b\mP\big(\ovec\dbc\bb,~b\ctx{\tilde\bC^b(\ovec)}
 \xvec\big).
\end{align}
The event $\{\vvec\ct{\bC}\xvec\}$ can be decomposed into two cases
depending on whether there is or is not a pivotal bond $b$ for
$\vvec\conn\xvec$ such that $\vvec\ct{\bC}\bb$.
Let
\begin{align}
E'(\vvec,\yvec;\bC)&=\{\vvec\ct{\bC}\yvec\}\cap\big\{\nexists\,b
 \text{ pivotal for }\vvec\conn\yvec\text{ s.t. }\vvec\ct{\bC}\bb
 \big\},\lbeq{E'def}\\[5pt]
E(b,\yvec;\bC)&=\{b\text{ occupied}\}\cap E'(\tb,\yvec;\bC).\lbeq{Edef}
\end{align}
See Figure~\ref{fig-E} for a schematic representation of the event
$E(b,\xvec;\bC)$.
\begin{figure}[t]
\begin{center}
\setlength{\unitlength}{0.00067in}
{
\begin{picture}(1000,3000)(0,500)
{
\put(500, 2950){\ellipse{200}{400}}
\path(500,2600)(500,2750)
\put(500, 2400){\ellipse{200}{400}}
\path(500,2050)(500,2200)
\put(500, 1850){\ellipse{200}{400}}
\path(500,1500)(500,1650)
\put(500, 1300){\ellipse{200}{400}}
\path(500,950)(500,1100)
\put(500, 750){\ellipse{200}{400}}

\put(230,450){\makebox(0,0)[lb]{$b$}}
\put(410,3300){\makebox(0,0)[lb]{$\xvec$}}
\put(-300,2000){\makebox(0,0)[lb]{$\bC$}}

\thicklines
\qbezier(-200,1000)(0,2500)(650,3100)
\path(430,550)(570,550)
\path(430,450)(570,450)
\thinlines
}
\end{picture}
}
\end{center}
\caption{\label{fig-E}
Schematic representation of the event $E(b,\xvec;\bC)$.}
\end{figure}
If there are pivotal bonds for $\vvec\conn\xvec$, then we take the
{\it first} such pivotal bond $b$ for which $\vvec\ct{\bC}\bb$.
Therefore, we have the partition
\begin{align}\lbeq{decompE2}
\{\vvec\ct{\bC}\xvec\}=E'(\vvec,\xvec;\bC)\DDcup\BDcup{b}\big\{E'(\vvec,
 \bb;\bC)\cap\{b\text{ occupied \& pivotal for }\vvec\conn\xvec\}\big\}.
\end{align}
Defining
\begin{align}\lbeq{pi1-def}
\pi^{\sss(1)}(\yvec)=\sum_b\mP\big(\{\ovec\dbc\bb\}\cap
 E(b,\xvec;\tilde\bC^b(\ovec))\big),
\end{align}
we obtain
\begin{align}
R^{\sss(0)}(\xvec)=\pi^{\sss(1)}(\xvec)+\sum_{b_1,b_2}\mP\big(\{
 \ovec\dbc\bb_1\}\cap E(b_1,\bb_2;\tilde\bC^{b_1}(\ovec))\cap\{b_2
 \text{ occupied \& pivotal for }\tb_1\conn\xvec\}\big).
\end{align}

To the second term, we apply the inclusion-exclusion relation
\begin{align}\lbeq{incl/exclSa}
\{b\text{ occupied \& pivotal for }\vvec\conn\xvec\}&=\{\vvec\conn\bb,~
 b\conn\xvec\}\setminus\big\{\vvec\conn\bb,~b\ctx{\tilde\bC^{b}(\vvec)}
 \xvec\big\}.
\end{align}
We define
\begin{align}
R^{\sss(1)}(\xvec)=\sum_{b_1,b_2}\mP\big(\{\ovec\dbc\bb_1\}\cap E(b_1,
 \bb_2;\tilde\bC^{b_1}(\ovec))\cap\big\{b_2\ctx{\tilde\bC^{b_2}(\tb_1)}
 \xvec\big\}\big),
\end{align}
so that we obtain
\begin{align}
R^{\sss(0)}(\xvec)&=\pi^{\sss(1)}(\xvec)+\sum_{b_1,b_2}\mP\big(\{\ovec\dbc
 \bb_1\}\cap E(b_1,\bb_2;\tilde\bC^{b_1}(\ovec))\cap\{b_2\conn\xvec\}\big)
 -R^{\sss(1)}(\xvec),
\end{align}
where we use that
\begin{align}
E'(\vvec,\bb;\bC)\cap\{\vvec\conn\bb,~b\conn\xvec\}=E'(\vvec,\bb;\bC)
 \cap\{b\conn\xvec\}.
\end{align}
The event $\{\ovec\dbc\bb_1\}\cap E(b_1,\bb_2;\tilde\bC^{b_1}(\ovec))$
depends only on bonds before $\bb_2$, while $\{b_2\conn\xvec\}$ depends
only on bonds after $\bb_2$.  By the Markov property, we end up with
\begin{align}\lbeq{Rn1rewr2}
R^{\sss(0)}(\xvec)&=\pi^{\sss(1)}(\xvec)+\sum_{b_2}\pi^{\sss(1)}
 (\bb_2)\,p(b_2)\,\tau(\xvec-\tb_2)-R^{\sss(1)}(\xvec)\nn\\
&=\pi^{\sss(1)}(\xvec)+(\pi^{\sss(1)}\sstar p\sstar\tau)(\xvec)
 -R^{\sss(1)}(\xvec),
\end{align}
so that
\begin{align}\lbeq{exptau-5}
\tau(\xvec)=\pi^{\sss(0)}(\xvec)-\pi^{\sss(1)}(\xvec)+\big((\pi^{\sss
 (0)}-\pi^{\sss(1)})\sstar p\sstar\tau\big)(\xvec)+R^{\sss(1)}(\xvec).
\end{align}
This completes the second step of the expansion.

To complete the expansion for $\tau(\xvec)$, we need to investigate
$R^{\sss(1)}(\xvec)$ in more detail.  Note that $R^{\sss(1)}(\xvec)$
involves the probability of a subset of
$\big\{b_2\ctx{\tilde\bC^{b_2}(\bb_1)}\xvec\big\}$.

For this subset, we will use \refeq{decompE2} and \refeq{incl/exclSa}
again, and follow the steps of the above proof.  The expansion is completed
by repeating the above steps indefinitely.  To facilitate the statement and
the proof of the expansion, we make a few more definitions.  For
$\vec b_{\sN}=(b_1,\dots,b_{\sN})$ with $N\geq1$, we define
\begin{align}\lbeq{tildeE-def}
\tilde E_{\vec b_{\sN}}^{\sss(N)}(\xvec)=\{\ovec\dbc\bb_1\}\cap\bigcap_{i
 =1}^{N-1}E\big(b_i,\bb_{i+1};\tilde\bC^{b_i}(\tb_{i-1})\big)\cap E\big(
 b_{\sN},\xvec;\tilde\bC^{b_{\sN}}(\tb_{\sss N-1})\big),
\end{align}
where we use the convention that $\tb_{\sss 0}=\ovec$ and that the empty
intersection, arising when $N=1$, is the whole probability space.
Also, we let
    \begin{align}
    \tilde E_{\vec b_0}^{\sss(0)}(\xvec)
    =\{\ovec \dbc \xvec\}.
    \end{align}
Using this notation, we define
    \begin{align}
    \pi^{\sss(N)}(\xvec)&=\sum_{\vec b_{\sN}}\mP\big(\tilde
    E_{\vec b_{\sN}}^{\sss(N)}(\xvec)\big),\lbeq{pi-def}
    \end{align}
and denote the alternating sum by
    \eq\lbeq{pi-defalt}
    \pi(\xvec)=\sum_{N=0}^\infty(-1)^N \pi^{\sss(N)}(\xvec).
    \en
Note that the sum in \refeq{pi-defalt} is a finite sum, as long
as $t_{\xvec}$ is finite, where $t_{\xvec}$ denotes the time coordinate
of $\xvec$, since each of the bonds $b_{\sss 1},\dots,b_{\sN}$ eats up
at least one time-unit $\vep$, so that $\pi^{\sss(N)}(\xvec)=0$ for
$N\vep>t_{\xvec}$.  The result of the expansion is
summarized as follows.

\begin{prop}[The lace expansion]\label{prop:exp-two}
For any $\lamb\geq0$ and $\xvec\in\Lambda$,
    \begin{align}\lbeq{tau-2exp}
    \tau(\xvec)=\pi(\xvec)+(\pi\sstar p\sstar \tau)(\xvec).
    \end{align}
\end{prop}

\begin{figure}[t]
\begin{center}
\setlength{\unitlength}{0.0006in}
{
\begin{picture}(8000,3500)(-500,0)
{
\put(500, 1700){\ellipse{400}{2600}}
\put(-700,1600){\makebox(0,0)[lb]{$\pi^{\sss(0)}(\xvec)$}}
\put(430,3150){\makebox(0,0)[lb]{$\xvec$}}
\put(430,150){\makebox(0,0)[lb]{$\ovec$}}

\put(3000, 2800){\ellipse{200}{400}}
\path(3000,2450)(3000,2600)
\put(3000, 2250){\ellipse{200}{400}}
\path(3000,1900)(3000,2050)
\put(3000, 1700){\ellipse{200}{400}}
\path(3000,1350)(3000,1500)
\put(3000, 1150){\ellipse{200}{400}}
\path(3000,800)(3000,950)
\put(3000, 600){\ellipse{200}{400}}
\put(1850,1600){\makebox(0,0)[lb]{$\pi^{\sss(1)}(\xvec)$}}
\put(2930,3150){\makebox(0,0)[lb]{$\xvec$}}
\put(2930,150){\makebox(0,0)[lb]{$\ovec$}}
\qbezier(3100,600)(3700,1700)(3100,2800)

\put(5500, 2800){\ellipse{200}{400}}
\path(5500,2450)(5500,2600)
\put(5500, 2250){\ellipse{200}{400}}
\path(5500,1900)(5500,2050)
\put(5500, 1700){\ellipse{200}{400}}
\path(5500,1350)(5500,1500)
\put(5500, 1150){\ellipse{200}{400}}
\path(5500,800)(5500,950)
\put(5500, 600){\ellipse{200}{400}}
\qbezier(5600,600)(6000,1150)(5580,1600)
\qbezier(5580,1800)(6000,2250)(5600,2800)
\put(4350,1600){\makebox(0,0)[lb]{$\pi^{\sss(2)}(\xvec)$}}
\put(5430,3150){\makebox(0,0)[lb]{$\xvec$}}
\put(5430,150){\makebox(0,0)[lb]{$\ovec$}}
\put(6200,1500){\makebox(0,0)[lb]{$\bigcup$}}

\put(6900, 2800){\ellipse{200}{400}}
\path(6900,2450)(6900,2600)
\put(6900, 2250){\ellipse{200}{400}}
\path(6900,1900)(6900,2050)
\put(6900, 1700){\ellipse{200}{400}}
\path(6900,1350)(6900,1500)
\put(6900, 1150){\ellipse{200}{400}}
\path(6900,800)(6900,950)
\put(6900, 600){\ellipse{200}{400}}
\qbezier(7000,600)(7400,1150)(6980,1600)
\qbezier(7000,1150)(7400,2250)(7000,2800)
\put(6830,3150){\makebox(0,0)[lb]{$\xvec$}}
\put(6830,150){\makebox(0,0)[lb]{$\ovec$}}
}
\end{picture}
}
\end{center}
\caption{\label{fig-pi}
Schematic representations of $\pi^{\sss(0)}(\xvec)$,
$\pi^{\sss(1)}(\xvec)$ and $\pi^{\sss(2)}(\xvec)$.}
\end{figure}

\begin{proof}
By \refeq{Rn1rewr2}, we are left to identify $R^{\sss(1)}(\xvec)$.
For $N\geq1$, we define
\begin{align}\lbeq{Qn-deffin}
R^{\sss(N)}(\xvec)=\sum_{\vec b_N}\mP\big(\tilde E_{\vec b_{N-1}}^{\sss
 (N-1)}(\bb_{\sN})\cap\big\{b_{\sN}\ctx{\tilde\bC^{b_N}(\tb_{N-1})}\xvec
 \big\}\big).
\end{align}
We prove below
\begin{align}\lbeq{rem-ind}
R^{\sss(N)}(\xvec)=\pi^{\sss(N)}(\xvec)+(\pi^{\sss(N)}\sstar p\sstar\tau)
 (\xvec)-R^{\sss(N+1)}(\xvec).
\end{align}
The equation \refeq{tau-2exp} follows by repeated use of \refeq{rem-ind}
until the remainder $R^{\sss(N+1)}(\xvec)$ vanishes, which must happen
at least when $N\vep>t_{\xvec}$.  To complete the proof of
Proposition~\ref{prop:exp-two}, we are left to prove \refeq{rem-ind},
which is a generalization of \refeq{Rn1rewr2}.

First we rewrite $b_{\sN}\ctx{\tilde\bC^{b_N}(\tb_{N-1})}\xvec$ in
\refeq{Qn-deffin}.  As in \refeq{decompE2}, this event can be
decomposed into two cases, depending on whether there is or is not
a pivotal bond $b$ for $\tb_{\sN}\conn\xvec$ such that
$\tb_{\sN}\ctx{\tilde\bC^{b_N}(\tb_{N-1})}\bb$.  The contribution
where there is no such a bond equals
$E(b_{\sN},\xvec;\tilde\bC^{b_N}(\tb_{\sss N-1}))$.  If there are
such pivotal bonds, then we take the {\it first} bond $b$ among these
bonds and obtain (cf., \refeq{decompE2})
\begin{align}\lbeq{decompE}
\big\{b_{\sN}\ctx{\tilde\bC^{b_N}(\tb_{N-1})}\xvec\big\}&=E(b_{\sN},\xvec;
 \tilde\bC^{b_N}(\tb_{\sss N-1}))\nn\\
&\quad\DDcup\BDcup{b}\big\{E(b_{\sN},\bb;\tilde\bC^{b_N}(\tb_{\sss N-1}))
 \cap\{b\text{ occupied \& pivotal for }\tb_{\sN}\conn\xvec\}\big\}.
\end{align}
By \refeq{pi-def}, the contribution from
$E(b_{\sN},\xvec;\tilde\bC^{b_N}(\tb_{\sss N-1}))$ in the right-hand side
is $\pi^{\sss(N)}(\xvec)$, which is the first term in the right-hand side
of \refeq{rem-ind}.  For the contribution from the remaining event in
\refeq{decompE}, we use \refeq{incl/exclSa} to arrive at
\begin{gather}
\sum_{\vec b_N,b}\mP\big(\tilde E_{\vec b_N}^{\sss(N)}(\bb)\cap\{b
 \text{ occupied \& pivotal for }\tb_{\sN}\conn\xvec\}\big)
=\sum_{\vec b_N,b}\mP\big(\tilde E_{\vec b_N}^{\sss(N)}(\bb)\cap\{
 b\conn\xvec\}\big)-R^{\sss(N+1)}(\xvec).\lbeq{term2}
\end{gather}
The last term in the above expression is the last term in the right-hand
side of \refeq{rem-ind}.  Again by the Markov property and \refeq{pi-def},
the first term in the right-hand side of \refeq{term2} equals the second
term in the right-hand side of \refeq{rem-ind}.
This completes the proof of \refeq{rem-ind} and thus the proof of
Proposition~\ref{prop:exp-two}.
\end{proof}

\subsection{Representation for the derivative}\label{ss:rep-der}
In this section, we derive a formula for $\partial_\lamb\pi(\xvec)$.
To state the result below, we define
\begin{align}
\piv[\vvec,\xvec]=\{b:b\text{ pivotal for }\vvec\conn\xvec\}.
\end{align}

\begin{prop}\label{prop:exp-der}
For $\lamb>0$ and $\xvec\in\Lambda$,
\begin{align}\lbeq{pider-exp}    \partial_\lamb\pi(\xvec)=\frac1\lamb\sum_{N=1}^\infty(-1)^N\,
 \Pi^{\sss(N)}(\xvec),
\end{align}
where $\Pi^{\sss(N)}(\xvec)=\sum_{n=1}^N\Pi^{\sss(N;n)}(\xvec)$ with
\begin{align}\lbeq{Pider-exp}
\Pi^{\sss(N;n)}(\xvec)=\sum_{\substack{\vec b_N,b:\\ b\text{ spatial}}}
 \mP\big(\tilde E_{\vec b_N}^{\sss(N)}(\xvec)\cap\big\{b\in\{b_n\}
 \Dcup{}{}\piv[\tb_n,\bb_{n+1}]\big\}\big),
\end{align}
and where $\bb_{\sss N+1}$ is defined to be $\xvec$.
\end{prop}

By the same reason as in \refeq{pi-defalt}, the sum in \refeq{pider-exp} is
a finite sum.  We prove \refeq{pider-exp} by differentiating the expression
\refeq{tau-2exp} for $\tau(\xvec)$ and comparing it with the expression for
$\partial_\lamb\tau(\xvec)$ obtained by using Russo's formula, rather than
differentiating $\partial_\lamb\pi_t(x)$ directly.  Possibly, one can
also use direct differentiation of the expressions \refeq{pi-def} for
$\pi(\xvec)$, but this is cumbersome because of the complex combination
of {\it increasing} and {\it decreasing} events consisting of
$\pi^{\sss(N)}(\xvec)$, where an increasing (respectively, decreasing)
event is an event that is more (respectively, less) likely to occur as
$\lamb$ increases.  We note that, instead of a difference of two terms due
to the pivotals for the increasing and decreasing events, we only obtain a
single sum over pivotals.  Thus, an intricate cancellation takes place.
This is further demonstrated by the fact that there is no contribution from
$N=0$.  In particular, it is {\it not} true that
$\partial_\lamb\pi^{\sss(N)}(\xvec)=\frac1\lamb\Pi^{\sss(N)}(\xvec)$.

\begin{proof}
In the proof it will be convenient to split
\begin{align}
\Pi^{\sss(N,n)}(\xvec)=\Pi^{\sss(N;n,1)}(\xvec)+\Pi^{\sss(N;n,2)}(\xvec),
\end{align}
where $\Pi^{\sss(N;n,1)}(\xvec)$ is the contribution from $b=b_n$ in
\refeq{Pider-exp}, whereas $\Pi^{\sss(N;n,2)}(\xvec)$ is the contribution
from $b\in\piv[\tb_n,\bb_{n+1}]$.

To obtain an expression for $\partial_\lamb\tau(\xvec)$, we use
{\it Russo's formula} \cite{bg91,g99}.  Let $E$ be an increasing
event that depends only on finitely many spatial bonds.  Then
\begin{align}\lbeq{Russo}
\partial_\lamb\mP(E)=\frac1\lamb\sum_{b\text{ spatial}}
 \mP(E\text{ occurs},~b\text{ is pivotal for }E),
\end{align}
where we use the fact that only spatial pivotal bonds for $E$ are
responsible to the differentiation with respect to $\lamb$.  Let
$\Box_{\sR}=[-R,R]^d\cap\Zd$.  We apply \refeq{Russo} to $E=E_{\sR}
(\xvec)\equiv\{\ovec\conn\xvec$ in $\Box_{\sR}\times[0,t_{\xvec}]\}$,
which is the set of bond configurations whose restriction on bonds
$(\uvec,\vvec)\subset\Box_{\sR}\times[0,t_{\xvec}]$ are in
$\{\ovec\conn\xvec\}$.  Note that
$\lim_{R\to\infty}\mP(E_{\sR}(\xvec))=\tau(\xvec)$, and that,
for any $\lamb_0\in[0,\infty)$,
\begin{align}
\lim_{R\to\infty}\partial_\lamb\mP(E_{\sR}(\xvec))=\frac1\lamb\sum_{
 b\text{ spatial}}\mP(\ovec\conn\xvec,~b\text{ is pivotal for }\ovec
 \conn\xvec),\lbeq{uniformconv}
\end{align}
uniformly in $\lamb\in[0,\lamb_0]$, which we will show at the end
of this section.  Therefore, we can exchange the order of
$\lim_{R\to\infty}$ and $\partial_\lamb$, and obtain
\begin{align}
\partial_\lamb\tau(\xvec)=\frac1\lamb\sum_{b\text{ spatial}}\mP(\ovec
 \conn\xvec,~b\text{ is pivotal for }\ovec\conn\xvec).\lbeq{russo-appl}
\end{align}

We follow the same strategy as in Section~\ref{ss:exp-two} to obtain
a recursion relation, now for $\partial_\lamb\tau(\xvec)$ rather than
for $\tau(\xvec)$.  Then, \refeq{russo-appl} equals
\begin{align}
\partial_\lamb\tau(\xvec)=\frac1\lamb\sum_{b\text{ spatial}}\bigg[
&\mP(\ovec\db\bb,~b\text{ occupied \& pivotal for }\ovec\conn\xvec)\nn\\
&+\sum_{b_1<b}\mP(\ovec\db\bb_1,~b_1\text{ and }b
 \text{ occupied \& pivotal for }\ovec\conn\xvec)\bigg],\lbeq{russo}
\end{align}
where $\sum_{b_1<b}$ is the sum over bonds $b_1$ with $t_{\bb_1}<t_{\bb}$.
The first and second terms in the brackets of the right-hand side
correspond respectively to when $b$ is or is not the first element of
$\piv[\ovec,\xvec]$.  The contribution from the first term is the same as
\refeq{exptau-1}, apart from the factor $\frac1\lamb$ and the restriction
that $b$ has to be a spatial bond.  Thus, the first term equals
\begin{align}
(\pi^{\sss(0)}\sstar\vep D\sstar\tau)(\xvec)+\frac1\lamb\sum_{N=1}^\infty
 (-1)^N\big[\Pi^{\sss(N;1,1)}(\xvec)+(\Pi^{\sss(N;1,1)}\sstar p\sstar\tau)
 (\xvec)\big],\lbeq{der-exp00}
\end{align}
where we abuse notation to write
\begin{align}
D((y,s))=D(y)\,\delta_{s,\vep}.
\end{align}

For the second term in \refeq{russo}, we use
\begin{align}\lbeq{diffpivs}
&\{b_1\text{ and }b\text{ occupied \& pivotal for }\ovec\conn\xvec\}\nn\\
&\quad=\{b_1\text{ occupied \& pivotal for }\ovec\conn\xvec\}\cap\{b\in
 \piv[\tb_1,\xvec]\}.
\end{align}
We ignore the condition that $b_1$ is pivotal for $\ovec\conn\xvec$ and use
inclusion-exclusion in the form \refeq{incl/exclSa} to make up for the
arising error.  Using the Markov property, the contribution from the second
term in \refeq{russo} is
\begin{align}\lbeq{der-exp01}
\frac1\lamb\sum_{\substack{b_1,b:\\ b\text{ spatial}}}\mP\big(\ovec\db
 \bb_1,~b_1\conn\xvec,~b\in\piv[\tb_1,\xvec]\big)-Q^{\sss(1)}(\xvec)
 =(\pi^{\sss(0)}\sstar p\sstar\partial_\lamb\tau)(\xvec)-Q^{\sss(1)}(\xvec),
\end{align}
where we define $Q^{\sss(n)}(\xvec)$ by
\begin{align}\lbeq{Rn-def}
Q^{\sss(n)}(\xvec)&=\frac1\lamb\sum_{\substack{\vec b_n,b:\\ b
 \text{ spatial}}}\mP\big(\tilde E_{\vec b_{n-1}}^{\sss(n-1)}(\bb_n)\cap
 \big\{b_n\ctx{\tilde\bC^{b_n}(\tb_{n-1})}\xvec\big\}\cap\{b\in\piv[\tb_n,
 \xvec]\}\big),
    \end{align}
and we recall that $\tilde E^{\sss(0)}(\bb_1)=\{\ovec\db\bb_1\}$ and
$\tb_0=\ovec$.  Note that $Q^{\sss(n)}(\xvec)$ is similar to
$R^{\sss(n)}(\xvec)$ in \refeq{Qn-deffin}, except for the factor
$\frac1\lamb$, the sum over spatial bonds $b$, and the extra condition
$b\in\piv[\tb_n,\xvec]$.  Therefore, by \refeq{russo}--\refeq{der-exp00}
and \refeq{der-exp01}, we have
\begin{align}\lbeq{russo-result}
\partial_\lamb\tau(\xvec)&=(\pi^{\sss(0)}\sstar\vep D\sstar\tau)(\xvec)
 +(\pi^{\sss(0)}\sstar p\sstar\partial_\lamb\tau)(\xvec)-Q^{\sss(1)}
 (\xvec)\nn\\
&\qquad+\frac1\lamb\sum_{N=1}^\infty(-1)^N\big[\Pi^{\sss(N;1,1)}(\xvec)
 +(\Pi^{\sss(N;1,1)}\sstar p\sstar\tau)(\xvec)\big].
\end{align}

Below, we will use inclusion-exclusion to prove that, for $n\geq1$,
\begin{align}\lbeq{derem-ind}
Q^{\sss(n)}(\xvec)&=(\pi^{\sss(n)}\sstar\vep D\sstar\tau)(\xvec)+(\pi^{\sss(n)}
 \sstar p\sstar\partial_\lamb\tau)(\xvec)-Q^{\sss(n+1)}(\xvec)\nn\\
&\qquad+\frac1\lamb\sum_{N=n+1}^\infty(-1)^{N-n}\big[\Pi^{\sss(N;n+1,1)}
 (\xvec)+(\Pi^{\sss(N;n+1,1)}\sstar p\sstar\tau)(\xvec)\big]\nn\\
&\qquad+\frac1\lamb\sum_{N=n}^\infty(-1)^{N-n}\big[\Pi^{\sss(N;n,2)}
 (\xvec)+(\Pi^{\sss(N;n,2)}\sstar p\sstar\tau)(\xvec)\big].
\end{align}
Before proving \refeq{derem-ind}, we complete the proof of \refeq{pider-exp}
assuming \refeq{derem-ind}.  By repeated applications of \refeq{derem-ind}
to \refeq{russo-result} until the remainder term $Q^{\sss(n)}(\xvec)$
vanishes, we obtain
\begin{align}
\partial_\lamb\tau(\xvec)&-(\pi\sstar\vep D\sstar\tau)(\xvec)
 +(\pi\sstar p\sstar\partial_\lamb\tau)(\xvec)\nn\\
&\quad=\frac1\lamb\sum_{N=1}^\infty(-1)^N\big[\Pi^{\sss(N;1,1)}(\xvec)
 +(\Pi^{\sss(N;1,1)}\sstar p\sstar\tau)(\xvec)\big]\nn\\
&\qquad+\frac1\lamb\sum_{n=1}^\infty(-1)^n\sum_{N=n+1}^\infty(-1)^{N-n}
 \big[\Pi^{\sss(N;n+1,1)}(\xvec)+(\Pi^{\sss(N;n+1,1)}\sstar p\sstar\tau)
 (\xvec)\big]\nn\\
&\qquad+\frac1\lamb\sum_{n=1}^\infty(-1)^n\sum_{N=n}^\infty(-1)^{N-n}
 \big[\Pi^{\sss(N;n,2)}(\xvec)+(\Pi^{\sss(N;n,2)}\sstar p\sstar\tau)
 (\xvec)\big]\nnmb\\
&\quad=\frac1\lamb\sum_{N=1}^\infty(-1)^N\big[\Pi^{\sss(N)}(\xvec)
 +(\Pi^{\sss(N)}\sstar p\sstar\tau)(\xvec)\big].
\end{align}
Differentiating both sides of \refeq{tau-2exp} and comparing with the
above expression, we obtain
\begin{align}\lbeq{pider-induction}
\partial_\lamb\pi(\xvec)+(\partial_\lamb\pi\sstar p\sstar\tau)(\xvec)=\frac1
 \lamb\sum_{N=1}^\infty(-1)^N\big[\Pi^{\sss(N)}(\xvec)+(\Pi^{\sss(N)}\sstar
 p\sstar\tau)(\xvec)\big].
\end{align}
Using this identity, we prove \refeq{pider-exp} by induction on
$t_{\xvec}/\vep$.  Since $\pi((x,0))=\delta_{o,x}$ and
$\Pi^{\sss(N)}((x,0))=0$ for all $N\geq1$, we obtain \refeq{pider-exp}
for $t_{\xvec}/\vep=0$.  Suppose that \refeq{pider-exp} holds for all
$t_{\xvec}/\vep\leq m$.  Then the contribution from the second term in
the brackets of \refeq{pider-induction} equals the second term on the
left-hand side of \refeq{pider-induction}, and thus \refeq{pider-exp}
for $t_{\xvec}/\vep=m+1$ holds.  This completes the inductive proof of
\refeq{pider-exp}.

In order to complete the proof of Proposition \ref{prop:exp-der},
we prove \refeq{derem-ind}.  Because of the condition
$b_n\ctx{\tilde\bC^{b_n}(\tb_{n-1})}\xvec$ in \refeq{Rn-def}, either
the event $E(b_n,\xvec;\tilde\bC^{b_n}(\tb_{n-1}))$ occurs or there
is an occupied bond $b_{n+1}\in\piv[\tb_n,\xvec]$ for which the event
$E(b_n,\bb_{n+1};\tilde\bC^{b_n}(\tb_{n-1}))$ occurs.  The contribution
from the former case to $Q^{\sss(n)}(\xvec)$ is
\begin{align}\lbeq{der-exp10-}
\frac1\lamb\sum_{\substack{\vec b_n,b:\\ b\text{ spatial}}}\mP\big(
 \tilde E_{\vec b_{n-1}}^{\sss(n-1)}(\bb_n)\cap E(b_n,\xvec;\tilde
 \bC^{b_n}(\tb_{n-1}))\cap\{b\in\piv[\tb_n,\xvec]\}\big)
=\frac1\lamb\Pi^{\sss(n;n,2)}(\xvec).
\end{align}
The contribution from the latter case is, as in \refeq{russo},
\begin{align}\lbeq{IE3}
\frac1\lamb\sum_{\substack{\vec b_n,b:\\ b\text{ spatial}}}\bigg[&\mP\big(
 \tilde E_{\vec b_n}^{\sss(n)}(\bb)\cap\{b\text{ occupied \& pivotal for }
 \tb_n\conn\xvec\big)\nn\\
&+\sum_{b_{n+1}<b}\mP\big(\tilde E_{\vec b_n}^{\sss(n)}(\bb_{n+1})\cap
 \{b_{n+1}\text{ and }b\text{ occupied \& pivotal for }\tb_n\conn\xvec\}
 \big)\nn\\
&+\sum_{b_{n+1}>b}\mP\big(\tilde E_{\vec b_n}^{\sss(n)}(\bb_{n+1})\cap
 \{b\text{ and }b_{n+1}\text{ occupied \& pivotal for }\tb_n\conn\xvec\}
 \big)\bigg],
\end{align}
where the first, second and third terms in the brackets correspond
respectively to when $b_{n+1}=b$, when $b_{n+1}$ is between $\tb_n$ and
$\bb$, and when $b_{n+1}$ is between $\tb$ and $\xvec$.  The first term
is similar to that in \refeq{russo}, and its contribution equals, as in
\refeq{der-exp00},
\begin{align}\lbeq{der-exp10}
(\pi^{\sss(n)}\sstar\vep D\sstar\tau)(\xvec)+\frac1\lamb
 \sum_{N=n+1}^\infty(-1)^{N-n}\big[\Pi^{\sss(N;n+1,1)}(\xvec)
 +(\Pi^{\sss(N;n+1,1)}\sstar p\sstar\tau)(\xvec)\big].
\end{align}

For the second term in \refeq{IE3}, we apply \refeq{diffpivs}, with $b_1$
and $\ovec$ being replaced respectively by $b_{n+1}$ and $\tb_n$, and use
the inclusion-exclusion relation \refeq{incl/exclSa} and the Markov
property.  Then, the contribution from the second term equals, as in
\refeq{der-exp01},
\begin{align}
&\frac1\lamb\sum_{\substack{\vec b_{n+1},b:\\ b\text{ spatial}}}\mP\big(
 \tilde E_{\vec b_n}^{\sss(n)}(\bb_{n+1})\cap\{b_{n+1}
 \text{ occupied \& pivotal for }\tb_n\conn\xvec\}\cap\{b\in\piv
 [\tb_{n+1},\xvec]\}\big)\nn\\
&\qquad=\frac1\lamb\sum_{\substack{\vec b_{n+1},b:\\ b\text{ spatial}}}
 \mP\big(\tilde E_{\vec b_n}^{\sss(n)}(\bb_{n+1})\cap\{b_{n+1}\conn\xvec\}
 \cap\{b\in\piv[\tb_{n+1},\xvec]\}\big)-Q^{\sss(n+1)}(\xvec)\nn\\
&\qquad=(\pi^{\sss(n)}\sstar p\sstar\partial_\lamb\tau)(\xvec)-Q^{\sss(n+1)}
 (\xvec).\lbeq{der-exp13}
\end{align}

For the third term in \refeq{IE3}, we use
\begin{align}
&\{b\text{ and }b_{n+1}\text{ occupied \& pivotal for }
 \tb_n\conn\xvec\}\nn\\
&\quad=\{b\in\piv[\tb_n,\bb_{n+1}]\}\cap\{b_{n+1}
 \text{ occupied \& pivotal for }\tb_n\conn\xvec\}.
\end{align}
By the inclusion-exclusion relation \refeq{incl/exclSa}, the contribution
from the third term equals
\begin{align}
&\frac1\lamb\sum_{\substack{\vec b_{n+1},b:\\ b\text{ spatial}}}\mP\big(
 \tilde E_{\vec b_n}^{\sss(n)}(\bb_{n+1})\cap\{b\in\piv[\tb_n,\bb_{n+1}]\}
 \cap\{b_{n+1}\text{ occupied \& pivotal for }\tb_n\conn\xvec\}\big)\nn\\
&\quad=\frac1\lamb\sum_{\substack{\vec b_{n+1},b:\\ b\text{ spatial}}}
 \Big[\mP\big(\tilde E_{\vec b_n}^{\sss(n)}(\bb_{n+1})\cap\{b\in\piv[
 \tb_n,\bb_{n+1}]\}\cap\{b_{n+1}\conn\xvec\}\big)\nn\\
&\hspace{5pc}-\mP\big(\tilde E_{\vec b_n}^{\sss(n)}(\bb_{n+1})\cap\{b\in
 \piv[\tb_n,\bb_{n+1}]\}\cap\{b_{n+1}\ctx{\tilde \bC^{b_{n+1}}(\tb_n)}
 \xvec\}\big)\Big],\lbeq{der-exp110}
\end{align}
where the first term equals, by the Markov property,
\begin{align}\lbeq{der-exp11}
\frac1\lamb(\Pi^{\sss(n;n,2)}\sstar p\sstar\tau)(\xvec).
\end{align}
For the second term in \refeq{der-exp110}, we use the same argument
as above \refeq{decompE}.  Because of the condition
$b_{n+1}\ctx{\tilde\bC^{b_{n+1}}(\tb_n)}\xvec$, either the event
$E(b_{n+1},\xvec;\tilde\bC^{b_{n+1}}(\tb_n))$ occurs or there
is an occupied bond $b_{n+2}\in\piv[\tb_{n+1},\xvec]$ such that
$E(b_{n+1},\bb_{n+2};\tilde\bC^{b_{n+1}}(\tb_n))$ occurs.
By repeated use of inclusion-exclusion and the Markov property,
as above \refeq{decompE}, the contribution from the second term in
\refeq{der-exp110} equals
\begin{align}\lbeq{der-exp12}
\frac1\lamb\sum_{N=n+1}^{\infty}(-1)^{N-n}\big[\Pi^{\sss(N;n,2)}(\xvec)
 +(\Pi^{\sss(N;n,2)}\sstar p\sstar \tau)(\xvec)\big].
\end{align}

Combining \refeq{der-exp10-}, \refeq{der-exp10}--\refeq{der-exp13} and
\refeq{der-exp11}--\refeq{der-exp12}, we obtain \refeq{derem-ind}.
This completes the proof of Proposition~\ref{prop:exp-der}, assuming
the uniformity of \refeq{uniformconv}.
\end{proof}

\begin{proof}[Proof of the uniformity of \refeq{uniformconv}]
Given $\lamb_0\in[0,\infty)$, we prove that
$\partial_\lamb\mP(E_{\sR}(\xvec))$ converges to the right-hand
side of \refeq{uniformconv}, uniformly in $\lamb\in[0,\lamb_0]$.

Recall that $E_{\sR}(\xvec)=\{\ovec\conn\xvec$ in
$\Box_{\sR}\times[0,t_{\xvec}]\}$. The difference between
$\partial_\lamb\mP(E_{\sR}(\xvec))$ and the right-hand side of
\refeq{uniformconv} is bounded by
\begin{gather}
\frac1\lamb\sum_{\substack{b\text{ spatial}\\ b\subset\Box_R\times[0,
 t_{\xvec}]}}\mP\big(b\text{ occupied \& pivotal for }E_{\sR}(\xvec),
 \text{ but not pivotal for }\ovec\conn\xvec\big)\nn\\
+\frac1\lamb\sum_{\substack{b\text{ spatial}\\ b\not\subset\Box_R\times
 [0,t_{\xvec}]}}\mP(b\text{ occupied \& pivotal for }\ovec\conn\xvec).
 \lbeq{russobd}
\end{gather}
First, we bound the second term, using $\{b$ occupied \& pivotal for
$\ovec\conn\xvec\}\subset\{\ovec\conn\bb\}\cap\{b\conn\xvec\}$ as well
as the Markov property and \refeq{tauhatbd}, by
\begin{align}
\vep\!\!\sum_{b\not\subset\Box_R\times[0,t_{\xvec}]}\shift p^{\sstar t_{\bb}
 /\vep}(\bb)\;D(\tb-\bb)\;p^{\sstar(t_{\xvec}-t_{\tb})/\vep}(\xvec-\tb)\leq
 \vep\sum_{j=1}^{t_{\xvec}/\vep}(1-\vep+\lamb\vep)^{t_{\xvec}/\vep-j}\shift
 \sum_{y\in\Zd:\|y\|_\infty\geq R}\shift p_\vep^{*(j-1)}(y),\lbeq{russobd2}
\end{align}
where we take the sum over the spatial component of $\xvec$ to obtain
the bound.  Similarly to \refeq{tauhatbd-appl}, this is further bounded,
uniformly in $\lamb$ and $\vep$, by $c\delta_{\sR}'$
where $c=c(\lamb_0,t_{\xvec})$ and $\delta_{\sR}'=\delta_{\sR}'(\lamb_0)$
are some finite constants satisfying $\lim_{R\to\infty}\delta_{\sR}'=0$.

Next, we consider the first term in \refeq{russobd}.  Note that, if $b$ is
pivotal for $E_{\sR}(\xvec)$, but not pivotal for $\ovec\conn\xvec$, then
there must be a detour from some $\yvec\in\Box_{\sR}\times[0,t_{\bb}]$ to
another $\zvec\in\Box_{\sR}\times[t_{\tb},t_{\xvec}]$ that passes through
$(\Zd\setminus\Box_{\sR})\times[0,t_{\xvec}]$ without traversing $b$.
Therefore, the event in the first term of \refeq{russobd} is a subset of
\begin{align}
\bigcup_{\substack{\yvec,\zvec,\uvec\in\Box_R\times[0,t_{\xvec}]\\ \vvec
 \notin\Box_R\times[0,t_{\xvec}]}}\{\ovec\conn\yvec\}\cap\big\{\{\yvec
 \conn\bb,~b\conn\zvec\}\circ\{\yvec\conn\uvec,~(\uvec,\vvec)\conn\zvec\}
 \big\}\cap\{\zvec\conn\xvec\},\lbeq{russobd3}
\end{align}
where $(\uvec,\vvec)$ is the first bond along the detour that {\it crosses}
the boundary of $\Box_{\sR}\times[0,t_{\xvec}]$, so that it is a spatial
bond, and $E_1\circ E_2$ is the event that $E_1$ and $E_2$ occur
{\it disjointly}, i.e., there is a bond set $B$ such that $E_1$ occurs on
$B$ and $E_2$ occurs on the complement of $B$.  By the Markov property,
the three events joined by $\cap$ are independent of each other.  For the
middle event, we use the {\it van~den~Berg-Kesten (BK) inequality} \cite{g99},
which asserts that $\mP(E_1\circ E_2)\leq\mP(E_1)\,\mP(E_2)$
when both $E_1$ and $E_2$ are {\it increasing}
events, i.e., $E_1$ and $E_2$ are more likely to occur as $\lamb$
increases, as in \refeq{russobd3}. Sometimes, we also make use of the
{\it van~den~Berg-Kesten-Reimer (BKR) inequality} \cite{bcr99}, which
proves $\mP(E_1\circ E_2)\leq\mP(E_1)\,\mP(E_2)$ for {\it any} events
$E_1$ and $E_2$.
Then, we use \refeq{tauhatbd} as in
\refeq{russobd2} and obtain that the first term in \refeq{russobd} (even
when we sum over the spatial component of $\xvec$) is bounded by
$c'\delta_{\sR}''$, where $c'=c'(\lamb_0,t_{\xvec},\vep)$ and
$\delta''_{\sR}=\delta_{\sR}''(\lamb_0)$ are some finite constants satisfying
$\lim_{R\to\infty}\delta''_{\sR}=0$.  This completes the proof of the
uniformity in $\lamb$ of \refeq{uniformconv}.

In the above proof, we did not care about the uniformity in $\vep$,
since it has been fixed and positive in this section.  In fact, the
above constant $c'$ is of order $O(\vep^{-2})$ and diverges as $\vep\to0$.
This is because the contribution from \refeq{russobd3} involves the sums
over $t_{\tb}(=t_{\bb}+\vep),t_{\yvec},t_{\zvec},t_{\vvec}(=t_{\uvec}+\vep)
\in\vep\Zp$ that give rise to the factor $\vep^{-4}$.  However, the factor
$\vep^{-2}$ is cancelled by the bond occupation probabilities of the
{\it spatial} bonds $b$ and $(\uvec,\vvec)$, and therefore
$c'=O(\vep^{-2})$.  We could improve this to $c'=O(1)$ by using the ideas
in Section~\ref{ss:diag-bds}, and hence obtain the uniformity in $\vep$
as well, though this is not necessary here.
\end{proof}


\newcommand{\bspat}{B_{\tt spat}}
\newcommand{\btemp}{B_{\tt temp}}
\newcommand{\vhi}{\varphi}

\section{Bounds on the lace expansion}\label{s:pibds}
In this section, we prove
Propositions~\ref{lem-Pibd.smp}--\ref{lem-Pibd.smp2}.  By
\refeq{pi-def}--\refeq{pi-defalt} and \refeq{pider-exp}--\refeq{Pider-exp},
it suffices to prove the following bounds on $\pi_t^{\sss(N)}(x)$ and
$\Pi_t^{\sss(N;n)}(x)$ in order to prove these propositions.

\begin{lem}\label{lem:piPi-bds}
Suppose that \refeq{fbd} holds for some $\lamb_0$ and all $s\leq t$.
\begin{enumerate}
\item[(i)]
Let $d>4$.  Then, there are $\beta_0>0$ and $\Cg<\infty$ such that,
for $\lamb\leq\lamb_0$, $\beta<\beta_0$, $s\in\vep\Zp$ with
$2\vep\leq s\leq t+\vep$, and $q=0,2,4$,
\begin{align}
\sum_x|x|^q\,\pi_s^{\sss(N)}(x)&\leq\frac{\vep^2(\Cg\beta)^{1\vee N}\sigma^q
 N^{q/2}}{(1+s)^{(d-q)/2}},\qquad\text{for }N\geq0,\lbeq{piN-bd}\\
\sum_x\Pi_s^{\sss(N;n)}(x)&\leq\frac{\vep^2(\Cg\beta)^N}{(1+s)^{(d-2)/2}},
 \qquad\qquad\text{for }N\geq n\geq1.\lbeq{PiNn-bd}
\end{align}
\item[(ii)]
Let $d\leq4$ with $\alpha=bd-\frac{4-d}2>0$, $\mu\in(0,\alpha)$ and
$t\leq T\log T$.  Then, there are $\beta_0>0$ and $\Cg<\infty$ such that,
for $\lamb\leq\lamb_0$, $\beta_1<\beta_0$, $s\in\vep\Zp$ with
$2\vep\leq s\leq t+\vep$, and $q=0,2,4$,
\begin{align}
\sum_x|x|^q\,\pi_s^{\sss(N)}(x)&\leq\frac{\vep^2\Cg\beta_{\sT}(\Cg\hat\beta
 _{\sT})^{0\vee(N-1)}\sigma_{\sT}^qN^{q/2}}{(1+s)^{(d-q)/2}},\qquad
 \text{for }N\geq0,\lbeq{piN-bd-lower}\\
\sum_x\Pi_s^{\sss(N;n)}(x)&\leq\frac{\vep^2\Cg\beta_{\sT}(\Cg\hat\beta_{\sT})
 ^{N-1}}{(1+s)^{(d-2)/2}},\qquad\qquad\qquad\text{ for }N\geq n\geq1,
 \lbeq{PiNn-bd-lower}
\end{align}
where $\beta_{\sT}=\beta_1T^{-bd}$ and $\hat\beta_{\sT}=\beta_1T^{-\mu}$.
\end{enumerate}
\end{lem}

\begin{proof}[Proof of Propositions~\ref{lem-Pibd.smp}--\ref{lem-Pibd.smp2}
assuming Lemma~\ref{lem:piPi-bds}]
The inequalities \refeq{pimom} and \refeq{pider} follow from
\refeq{pi-def}--\refeq{pi-defalt}, \refeq{pider-exp}--\refeq{Pider-exp}
and \refeq{piN-bd}--\refeq{PiNn-bd}, if $\beta$ is sufficiently small.
The proof of \refeq{pidiff} is the same as that of Proposition~2.2(ii)
in \cite[Section~4.3]{hs01}, together with \refeq{pimom}.
This completes the proof of Proposition~\ref{lem-Pibd.smp}.

Proposition~\ref{lem-Pibd.smp2} is proved by using
\refeq{piN-bd-lower}--\refeq{PiNn-bd-lower} instead of
\refeq{piN-bd}--\refeq{PiNn-bd}, if $\hat\beta_{\sT}$ is sufficiently small.
\end{proof}

Lemma~\ref{lem:piPi-bds} is proved in
Sections~\ref{ss:diag-bds}--\ref{ss:diag-est-lower}.  In
Section~\ref{ss:diag-bds}, we first introduce certain diagram functions
$P_t^{\sss(N)}(x)$ and $\tilde P_t^{\sss(N;n)}(x)$ that are defined in
terms of two-point functions, and prove that these diagram functions are
upper bounds on $\pi_t^{\sss(N)}(x)$ and $\Pi_t^{\sss(N;n)}(x)$,
respectively.  Then, we bound these diagram functions assuming the bounds
in \refeq{fbd} on the two-point function, for $d>4$ in
Section~\ref{ss:diag-est} and for $d\leq4$ in
Section~\ref{ss:diag-est-lower}.  Finally, in
Section~\ref{ss:fc}, we use these diagram functions to obtain finite-volume
approximations of $\pi_t^{\sss(N)}(x)$ and $\Pi_t^{\sss(N;n)}(x)$, which
will be used in Section~\ref{s:continuum} to prove the continuum limit
$\vep\daw0$.

\subsection{Bounds in terms of the diagram functions}\label{ss:diag-bds}
In this section, we prove that $\pi_t^{\sss(N)}(x)$ and
$\Pi_t^{\sss(N;n)}(x)$ are bounded by certain diagram functions
$P_t^{\sss(N)}(x)$ and $\tilde P_t^{\sss(N;n)}(x)$ that are defined below
in terms of two-point functions.

The strategy in this section is similar to \cite[Section~4.1]{hs01} for
oriented percolation in $\Zd\times\Zp$, where bounds on
$\pi_t^{\sss(N)}(x)$ were proved by using some diagram functions arising
from the Hara-Slade lace expansion.  Since the expansion used in this paper
is somewhat simpler, we can use simpler diagram functions.  However, to
consider the case $\vep\ll1$ as in \cite{s01}, extra care is needed to
obtain the factor $\vep^2$ in \refeq{PN-bd}--\refeq{tPNn-bd-lower}.

\subsubsection{Preliminaries}
Before defining the diagram functions, we start by some preliminaries.
For $\vvec=(v,s)\in\Lambda$ and a bond $b$, we write $\vvec_+=(v,s+\vep)$
and $\{\vvec\conn b\}=\{\vvec\conn\bb\}\cap\{b\text{ occupied}\}$ (cf.,
Definition~\ref{def-1}(i)).  For convenience, we will also use
abbreviations, such as
\begin{align}
\{\vvec\conn b\conn\xvec\}&=\{\vvec\conn b\}\cap\{\tb\conn\xvec\}.
\end{align}
Let $I'(\vvec,\xvec,\xvec)=\{\vvec\conn\xvec\}$, and define,
for $\yvec\ne\xvec$,
\begin{align}\lbeq{I-def}
I'(\vvec,\yvec,\xvec)=\Bigg\{\bigcup_{\substack{b\text{ spatial}:\\ \bb=
 \yvec}}\{\vvec\conn b\conn\xvec\}\Bigg\}\cup\Bigg\{\bigcup_{\substack{b
 \text{ spatial}:\\ \tb=\yvec}}\{\vvec\conn(\bb,\bb_+)\conn\xvec,~b
 \text{ occupied}\}\Bigg\}.
\end{align}
We note that $I'(\vvec,\yvec,\xvec)$ for $\yvec\ne\xvec$ equals
$I'(\vvec,\xvec,\xvec)$ with an extra spatial bond $b$ being embedded
(or added) along the connection from $\vvec$ to $\xvec$.  Denoting
\begin{align}
I(b,\yvec,\xvec)=\{b\text{ occupied}\}\cap I'(\tb,\yvec,\xvec),
\end{align}
we define
\begin{align}
M(b,\vvec;\xvec,\yvec)=\big\{I(b,\yvec,\xvec)\circ\{\vvec\conn\xvec\}\big\}
 \cup\big\{\{b\conn\xvec\}\circ I'(\vvec,\yvec,\xvec)\big\}.\lbeq{M2-def}
\end{align}
Note that, when neither $\tb$ nor $\vvec$ is $\xvec$, the event
$M(b,\vvec;\xvec,\yvec)$ equals $\{b\conn\xvec\}\circ\{\vvec\conn\xvec\}$
with an extra spatial bond being embedded either between $\tb$ and
$\xvec$ due to $I(b,\yvec,\xvec)$, or between $\vvec$ and $\xvec$
due to $I'(\vvec,\yvec,\xvec)$.  In addition, we define
$M^{\sss+}(b,b',\vvec;\xvec,\yvec)$ to be $M(b,\vvec;\xvec,\yvec)$
with the connection from $\tb$ to $\xvec$ being replaced by
$\tb\conn b'\conn\xvec$.  For example, the second event
$\{b\conn\xvec\}\circ I'(\vvec,\yvec,\xvec)$ in \refeq{M2-def} is simply
replaced by $\{\tb\conn b'\conn\xvec\}\circ I'(\vvec,\yvec,\xvec)$ in the
definition of $M^{\sss+}(b,b',\vvec;\xvec,\yvec)$.  Replacing the first
event $I(b,\yvec,\xvec)\circ\{\vvec\conn\xvec\}$ in \refeq{M2-def} is more
complicated, due to the three possibilities of embedding $b'$ into
$\vvec\conn b\conn\xvec$ in \refeq{I-def} and the other three possibilities
of embedding $b'$ into $\vvec\conn(\bb,\bb_+)\conn\xvec$ in \refeq{I-def},
and therefore we refrain from giving a formula for
$M^{\sss+}(b,b',\vvec;\xvec,\yvec)$.

Recall \refeq{pi-def} and \refeq{Pider-exp} for the definitions of
$\pi^{\sss(N)}(\xvec)$ and $\Pi^{\sss(N;n)}(\xvec)$ that involve the event
$\tilde E_{\vec b_N}^{\sss(N)}(\xvec)$.  Our first claim is that
$\tilde E_{\vec b_N}^{\sss(N)}(\xvec)$ satisfies the following successive
relations:

\begin{lem}\label{lem:supevents}
For $N\geq1$,
\begin{align}
\tilde E_{\vec b_N}^{\sss(N)}(\xvec)&\subset\tilde E_{\vec b_{N-1}}
 ^{\sss(N-1)}(\bb_{\sN})\cap M(b_{\sN},\tb_{\sss N-1};\xvec,\xvec),
 \lbeq{tE0-supset}\\
\tilde E_{\vec b_{N-1}}^{\sss(N-1)}(\bb_{\sN})&\cap M(b_{\sN},\tb_{\sss
 N-1};\xvec,\yvec)\nn\\
&\subset\bigcup_{\vvec\in\Lambda}\Big\{\big\{\tilde E_{\vec b_{N-2}}
 ^{\sss(N-2)}(\bb_{\sss N-1})\cap M(b_{\sss N-1},\tb_{\sss N-2};\bb_{\sN},
 \vvec)\big\}\circ M(b_{\sN},\vvec;\xvec,\yvec)\Big\},\lbeq{tE-supset}\\
\tilde E_{\vec b_{N-1}}^{\sss(N-1)}(\bb_{\sN})&\cap M(b_{\sN},\tb_{\sss
 N-1};\xvec,\yvec)\cap\{b\in\piv[\tb_{\sN},\bb_{\sss N+1}]\}\nnmb\\
&\subset\bigcup_{\vvec\in\Lambda}\Big\{\big\{\tilde E_{\vec b_{\sss N-2}}
 ^{\sss(N-2)}(\bb_{\sss N-1})\cap M(b_{\sss N-1},\tb_{\sss N-2};\bb_{\sN},
 \vvec)\big\}\circ M^{\sss+}(b_{\sN},b,\vvec;\xvec,\yvec)\Big\},
 \lbeq{tEb-supset=}
\end{align}
and for $N>n\geq1$,
\begin{align}
&\tilde E_{\vec b_{N-1}}^{\sss(N-1)}(\bb_{\sN})\cap M(b_{\sN},\tb_{\sss
 N-1};\xvec,\yvec)\cap\{b\in\piv[\tb_n,\bb_{n+1}]\}\nn\\
&\quad\subset\bigcup_{\vvec\in\Lambda}\Big\{\big\{\tilde E_{\vec b_{\sss
 N-2}}^{\sss(N-2)}(\bb_{\sss N-1})\cap M(b_{\sss N-1},\tb_{\sss N-2};
 \bb_{\sN},\vvec)\cap\{b\in\piv[\tb_n,\bb_{n+1}]\}\big\}\circ M(b_{\sN},
 \vvec;\xvec,\yvec)\Big\},\lbeq{tEb-supset<}
\end{align}
where $\tb_0=\tb_{-1}=\ovec$,
$M(b_0,\tb_{-1};\bb_1,\vvec)=\{\ovec\db\bb_1\}\cap I'(\ovec,\vvec,\bb_1)$,
$\tilde E_{\vec b_0}^{\sss(0)}(\bb_1)=\{\ovec\db\bb_1\}$ and
$\tilde E_{\vec b_{-1}}^{\sss(-1)}(\bb_0)$ equals the whole probability
space.
\end{lem}

We note that the left-hand side of \refeq{tEb-supset<} is the same as that
of \refeq{tEb-supset=}, except that $b$ is pivotal for $\tb_n\conn\bb_{n+1}$
with $n<N$.

\begin{proof}
The relation \refeq{tE0-supset} follows immediately from
\refeq{tildeE-def} and
\begin{align}
E(b_{\sN},\xvec;\tilde \bC^{b_N}(\tb_{\sss N-1}))\subset\{b_{\sN}\conn\xvec
 \}\circ\{\tb_{\sss N-1}\conn\xvec\}=M(b_{\sN},\tb_{\sss N-1};\xvec,\xvec).
\end{align}

We only prove \refeq{tE-supset}, since
\refeq{tEb-supset=}--\refeq{tEb-supset<} can be proved similarly.  First, we
use \refeq{tE0-supset} to obtain $\tilde E_{\vec b_{\sss N-1}}^{\sss(N-1)}
(\bb_{\sN})\subset\tilde E_{\vec b_{\sss N-2}}^{\sss(N-2)}(\bb_{\sss N-1})
\cap M(b_{\sss N-1},\tb_{\sss N-2};\bb_{\sN},\bb_{\sN})$.
Since $\tilde E_{\vec b_{\sss N-2}}^{\sss(N-2)}(\bb_{\sss N-1})$ depends
only on bonds before $\bb_{\sss N-1}$, it suffices to prove
\begin{align}\lbeq{MM-supset}
M(b_{\sss N-1},\tb_{\sss N-2};\bb_{\sN},\bb_{\sN})\cap M(b_{\sN},\tb_{\sss
 N-1};\xvec,\yvec)\subset\bigcup_{\vvec\in\Lambda}\big\{M(b_{\sss N-1},
 \tb_{\sss N-2};\bb_{\sN},\vvec)\circ M(b_{\sN},\vvec;\xvec,\yvec)\big\}.
\end{align}

Recall that $M(b_{\sss N-1},\tb_{\sss N-2};\bb_{\sN},\bb_{\sN})=\{b_{\sss
N-1}\conn\bb_{\sN}\}\circ\{\tb_{\sss N-2}\conn\bb_{\sN}\}$.  The event in
the left-hand side of \refeq{MM-supset} implies existence of
$\vvec\in\bC(\tb_{\sss N-1})$ such that $\vvec\conn\bb_{\sN}$ and
$M(b_{\sN},\vvec;\xvec,\yvec)$ occur disjointly.  Therefore,
\begin{align}\lbeq{MM-presupset}
M(b_{\sss N-1}&,\tb_{\sss N-2};\bb_{\sN},\bb_{\sN})\cap M(b_{\sN},\tb_{\sss
 N-1};\xvec,\yvec)\nn\\
&\subset\bigcup_{\vvec\in\Lambda}\Big\{\big\{\{b_{\sss N-1}\conn\vvec\conn
 \bb_{\sN}\}\circ\{\tb_{\sss N-2}\conn\bb_{\sN}\}\big\}\nn\\
&\qquad\qquad\cup\big\{\{b_{\sss N-1}\conn\bb_{\sN}\}\circ\{\tb_{\sss N-2}
 \conn\vvec\conn\bb_{\sN}\}\big\}\Big\}\circ M(b_{\sN},\vvec;\xvec,\yvec).
\end{align}
We investigate the vicinity of $\vvec\in\Lambda$ in \refeq{MM-presupset},
where there are two disjoint connections, $\vvec\conn\bb_{\sN}$ and
$\vvec\conn\xvec$.  Since there is at most one temporal bond growing out of
each vertex in $\Lambda$, at least one of the two connections has to use a
spatial bond at $\vvec$.  Therefore,
\begin{align}
\{\vvec\conn\bb_{\sN}\}\circ\{\vvec\conn\xvec\}\subset\bigcup_{\substack{
 b\text{ spatial}:\\ \bb=\vvec}}\Big\{\big\{\{b\conn\bb_{\sN}\}\circ\{\bb
 \conn\xvec\}\big\}\cup\big\{\{(\bb,\bb_+)\conn\bb_{\sN}\}\circ\{b\conn
 \xvec\}\big\}\Big\}.
\end{align}
Substituting this relation into \refeq{MM-presupset} and relabelling
$\tb=\vvec$ in the latter event
$\{(\bb,\bb_+)\conn\bb_{\sN}\}\circ\{b\conn\xvec\}$, we obtain
\refeq{MM-supset}, and thus \refeq{tE-supset}.  This completes the proof.
\end{proof}

\subsubsection{Diagrammatic bounds}
Inspired by the successive relations \refeq{tE0-supset}--\refeq{tEb-supset<},
we inductively construct the diagram functions $P_t^{\sss(N)}(x)$ and
$\tilde P_t^{\sss(N;n)}(x)$ as follows.  For $b=(\uvec,\vvec)$ with
$\uvec=(u,s)$ and $\vvec=(v,s+\vep)$, we abuse notation to write $p(b)$ or
$p(\vvec-\uvec)$ for $p_\vep(v-u)$, and $D(b)$ or $D(\vvec-\uvec)$ for
$D(v-u)$.  Let
\begin{align}\lbeq{vhi-def}
\vhi(\xvec-\uvec)&=\delta_{\uvec,\xvec}+(p\sstar\tau)(\xvec-\uvec),
\end{align}
and
\begin{align}\lbeq{Lne-def}
L(\uvec,\vvec;\xvec)&=\begin{cases}
 \vhi(\xvec-\uvec)~(\tau\sstar\lamb\vep D)(\xvec-\vvec)+(\vhi\sstar\lamb\vep
 D)(\xvec-\uvec)~\tau(\xvec-\vvec),&\text{if }\uvec\ne\vvec,\\
(D\sstar\tau)(\xvec-\uvec)~(\tau\sstar\lamb\vep D)(\xvec-\uvec)+(D\sstar\tau
 \sstar\lamb\vep D)(\xvec-\uvec)~\tau(\xvec-\uvec),&\text{if }\uvec=\vvec.
\end{cases}
\end{align}
We define
\begin{align}\lbeq{Pzero-def}
P^{\sss(0)}(\xvec)=\delta_{\ovec,\xvec}+\lamb\vep L(\ovec,\ovec;\xvec),
\end{align}
and define the {\it zeroth admissible lines} to be the two lines from
$\ovec$ to $\xvec$ in each diagram of $\lamb\vep L(\ovec,\ovec;\xvec)$.
With lines, we mean here $(\lamb\vep D\sstar\tau)(\xvec)$ and
$(\tau\sstar\lamb\vep D)(\xvec)$ for the contribution from the
first term in \refeq{Lne-def} with $\uvec=\vvec=\ovec$, and
$(\lamb\vep D\sstar\tau\sstar\lamb\vep D)(\xvec)$ and $\tau(\xvec)$
for the contribution from the second term in \refeq{Lne-def} with
$\uvec=\vvec=\ovec$.

Given an admissible line $\ell$ from $\vvec$ to $\xvec$ of a diagram
function, say $\tau(\xvec-\vvec)$ for simplicity, and given
$\yvec\ne\xvec$, Construction~$\bspat^\ell(\yvec)$ is defined to be the
operation in which $\tau(\xvec-\vvec)$ is replaced by
\begin{align}\lbeq{bspat}
\tau(\yvec-\vvec)~(\lamb\vep D\sstar\tau)(\xvec-\yvec),
\end{align}
and Construction~$\btemp^\ell(\yvec)$ is defined to be the operation in
which $\tau(\xvec-\vvec)$ is replaced by
\begin{align}\lbeq{btemp}
\sum_{b:\tb=\yvec}\tau(\bb-\vvec)~\lamb\vep D(b)~\mP((\bb,\bb_+)\conn\xvec).
\end{align}
We note that \refeq{bspat}--\refeq{btemp} are inspired by \refeq{I-def}.
The sum of the results of Construction~$\bspat^\ell(\yvec)$ and
Construction~$\btemp^\ell(\yvec)$ is simply said to be the result of
Construction~$B^\ell(\yvec)$.  We define Construction~$B^\ell(s)$ to be
the operation in which Construction~$B^\ell(y,s)$ is performed and then
followed by summation over $y\in\Zd$.  Construction~$\bspat^\ell(s)$ and
Construction~$\btemp^\ell(s)$ are defined similarly.

We denote the result of applying Construction~$B^\ell(\yvec)$ to a diagram
function $f(\xvec)$ by $f(\xvec,B^\ell(\yvec))$, and define
$f(\xvec,\bspat^\ell(\yvec))$ and $f(\xvec,\btemp^\ell(\yvec))$ similarly.
We construct $P^{\sss(N)}(\xvec)$ from $P^{\sss(N-1)}(\xvec)$ by
\begin{align}\lbeq{PNPN-1}
P^{\sss(N)}(\xvec)=2\lamb\vep\sum_{\vvec\in\Lambda}P^{\sss(N-1)}(\vvec)~
 L(\vvec,\vvec;\xvec)+\sum_\ell\sum_{\substack{\vvec,\yvec\in\Lambda\\
 \vvec\ne\yvec}}P^{\sss(N-1)}(\vvec,B^\ell(\yvec))~L(\vvec,\yvec;\xvec),
\end{align}
where $\sum_\ell$ is the sum over the $(N-1)^{\rm st}$ admissible
lines in each diagram.  We define
\begin{align}
\sum_\ell P^{\sss(N-1)}(\vvec,B^\ell(\vvec))=2\lamb\vep P^{\sss(N-1)}
 (\vvec).
\end{align}
Then, \refeq{PNPN-1} equals
\begin{align}\lbeq{PN-def}
P^{\sss(N)}(\xvec)=\sum_\ell\sum_{\vvec,\yvec\in\Lambda}P^{\sss(N-1)}
 (\vvec,B^\ell(\yvec))~L(\vvec,\yvec;\xvec).
\end{align}
We call the newly added lines, contained in $L(\vvec,\yvec;\xvec)$,
the {\it $N^{\rm th}$ admissible lines}.

\begin{figure}[t]
\begin{center}
\setlength{\unitlength}{0.0006in}
{
\begin{picture}(9000,3500)(-1000,0)
{
\qbezier(500,400)(-500,1700)(500,3000)
\qbezier(500,400)(1500,1700)(500,3000)
\put(-1200,1500){\makebox(0,0)[lb]{$P_t^{\sss(0)}(x)=$}}
\put(220,3050){\makebox(0,0)[lb]{$(x,t)$}}
\put(240,50){\makebox(0,0)[lb]{$(o,0)$}}

\qbezier(3000,400)(2000,1700)(3000,3000)
\qbezier(3000,400)(4000,1700)(3000,3000)
\put(1300,1500){\makebox(0,0)[lb]{$P_t^{\sss(1)}(x)=$}}
\put(2720,3050){\makebox(0,0)[lb]{$(x,t)$}}
\put(2740,50){\makebox(0,0)[lb]{$(o,0)$}}
\path(2555,1255)(3450,2150)
\thicklines
\path(3350,2150)(3550,2150)
\path(3350,2220)(3550,2220)
\thinlines

\qbezier(5600,400)(4600,1700)(5600,3000)
\qbezier(5600,400)(6600,1700)(5600,3000)
\path(5325,855)(6120,1650)
\thicklines
\path(6010,1650)(6210,1650)
\path(6010,1720)(6210,1720)
\thinlines
\path(5110,1710)(5915,2515)
\thicklines
\path(5815,2515)(6015,2515)
\path(5815,2585)(6015,2585)
\thinlines
\put(3900,1500){\makebox(0,0)[lb]{$P_t^{\sss(2)}(x)=$}}
\put(5320,3050){\makebox(0,0)[lb]{$(x,t)$}}
\put(5340,50){\makebox(0,0)[lb]{$(o,0)$}}
\put(6420,1500){\makebox(0,0)[lb]{$+$}}

\path(7165,765)(7900,1500)
\thicklines
\path(7800,1500)(8000,1500)
\path(7800,1570)(8000,1570)
\thinlines

\thicklines
\path(6980,2450)(7180,2450)
\path(6980,2520)(7180,2520)
\thinlines
\path(7900,1800)(7080,2450)

\qbezier(7400,400)(6400,1700)(7400,3000)
\qbezier(7400,400)(8400,1700)(7400,3000)
\put(7120,3050){\makebox(0,0)[lb]{$(x,t)$}}
\put(7140,50){\makebox(0,0)[lb]{$(o,0)$}}
}
\end{picture}
}
\end{center}
\caption{Graphical representations of $P_t^{\sss(0)}(x), P_t^{\sss(1)}(x)$
and $P_t^{\sss(2)}(x)$. Lines indicate two-point functions, and small bars
indicate a convolution with $p_{\vep}$.  Spatial bonds that are present at
all vertices in the diagrams are left implicit.}\label{fig-P}
\end{figure}
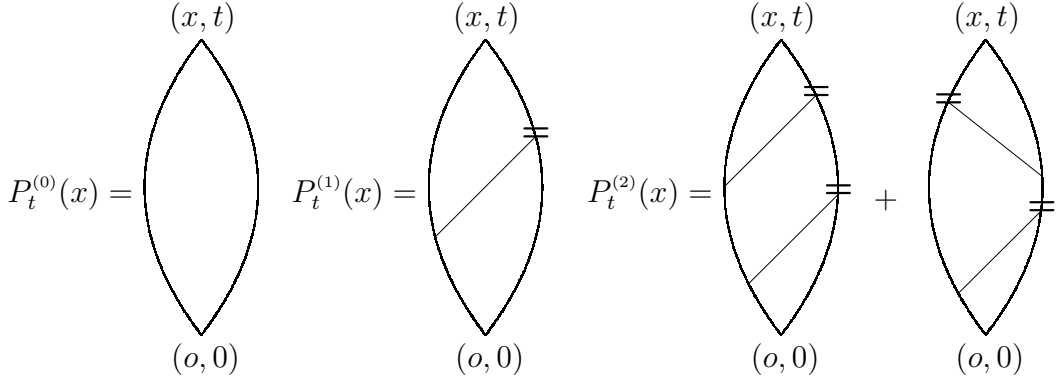

Finally, for $N\geq n\geq1$, we define
\begin{align}\lbeq{tPNn-def}
\tilde P^{\sss(N;n)}(\xvec)=\sum_\ell\sum_{\yvec\in\Lambda}P^{\sss(N)}
 (\xvec,\bspat^\ell(\yvec)),
\end{align}
where $\sum_\ell$ is the sum over the $n^{\rm th}$ admissible lines.

Thanks to the construction in terms of two-point functions, the diagram
functions can be estimated by using \refeq{fbd}, and this will be done in
Sections~\ref{ss:diag-est}--\ref{ss:diag-est-lower}.  The following is
the main statement of this section:

\begin{lem}\label{lem:BKbds}
For $\lamb\geq0$ and $N\geq n\geq1$,
\begin{align}\lbeq{BKbds}
\pi_t^{\sss(N-1)}(x)\leq P_t^{\sss(N-1)}(x),&&
\Pi_t^{\sss(N;n)}(x)\leq\tilde P_t^{\sss(N;n)}(x).
\end{align}
\end{lem}

\begin{proof}
We begin with proving $\pi^{\sss(0)}(\xvec)\leq P^{\sss(0)}(\xvec)$.
The first term in \refeq{Pzero-def} is the contribution from the case
of $\ovec=\xvec$.  If $\ovec\ne\xvec$, there are at least two nonzero
disjoint occupied paths from $\ovec$ to $\xvec$.  As explained below
\refeq{MM-presupset}, at least one of two nonzero disjoint occupied
paths from $\ovec$ has to use a spatial bond at $\ovec$.  That is,
\begin{align}\lbeq{pi0-prebd}
\pi^{\sss(0)}(\xvec)\leq\sum_{\substack{b\text{ spatial}:\\ \bb=\ovec}}
 \mP(\{b\conn\xvec\}\circ\{\ovec\conn\xvec\})=\sum_{\substack{b
 \text{ spatial}:\\ \bb=\ovec}}\mP(M(b,\ovec;\xvec,\xvec)).
\end{align}
We use the same observation at $\xvec$: at least one of the two nonzero
disjoint connections, $\tb\conn\xvec$ and $\ovec\conn\xvec$, has to use
another spatial bond at $\xvec$.  Therefore, we can bound the right-hand
side of \refeq{pi0-prebd} by $\lamb\vep L(\ovec,\ovec;\xvec)$ using the
BK inequality.  This completes the proof of
$\pi^{\sss(0)}(\xvec)\leq P^{\sss(0)}(\xvec)$.

Next, we consider $\pi^{\sss(N)}(\xvec)$ for $N\geq1$.  Let
\begin{align}\lbeq{piN-gendef}
\pi^{\sss(n)}(\xvec,\yvec)=\sum_{\vec b_n}\mP\big(\tilde E_{\vec
 b_{n-1}}^{\sss(n-1)}(\bb_n)\cap M(b_n,\tb_{n-1};\xvec,\yvec)\big).
\end{align}
By the convention in Lemma~\ref{lem:supevents}, $\pi^{\sss(0)}(\xvec,\yvec)
=\mP(\{\ovec\db\xvec\}\cap I'(\ovec,\yvec,\xvec))$.  We prove below by
induction that
\begin{gather}
\pi^{\sss(n)}(\xvec,\yvec)\leq(2\lamb\vep)^{-\delta_{\xvec,\yvec}}
 \sum_\ell P^{\sss(n)}(\xvec,B^\ell(\yvec))\lbeq{PN-ind}
\end{gather}
holds for all $n\geq0$, where $\sum_\ell$ is the sum over the $n^{\rm th}$
admissible lines.  The inequality
$\pi^{\sss(N)}(\xvec)\leq P^{\sss(N)}(\xvec)$ for $N\geq1$ follows from
\refeq{tE0-supset} and \refeq{PN-ind} for $\yvec=\xvec$, together with
the convention in \refeq{PN-def}, i.e.,
$2\lamb\vep P^{\sss(N)}(\xvec)=\sum_\ell P^{\sss(N)}(\xvec,B^\ell(\xvec))$.

For $n=0$, we can assume $\yvec\ne\xvec$, since
$\pi^{\sss(0)}(\xvec)=\pi^{\sss(0)}(\xvec,\xvec)\leq P^{\sss(0)}(\xvec)$
has already been proved.  By the equivalence $\{\ovec\db\xvec\}\cap
I'(\ovec,\yvec,\xvec)=\{\ovec\conn\xvec\}\circ I'(\ovec,\yvec,\xvec)$
and by \refeq{I-def}, we obtain
\begin{align}\lbeq{tM1unev-bd}
\pi^{\sss(0)}(\xvec,\yvec)&\leq\sum_{\substack{b\text{ spatial}:\\
 \bb=\yvec}}\mP(\{\ovec\conn\xvec\}\circ\{\ovec\conn b\conn\xvec\})\nn\\
&\qquad+\sum_{b:\tb=\yvec}\lamb\vep D(b)~\mP(\{\ovec\conn\xvec\}\circ\{
 \ovec\conn(\bb,\bb_+)\conn\xvec\}),
\end{align}
where we use the BK inequality to derive $\lamb\vep D(b)$ in the second sum.

We only prove that the first sum in \refeq{tM1unev-bd} is bounded by
$\sum_\ell P^{\sss(0)}(\xvec,\bspat^\ell(\yvec))$, by investigating the
vicinity of $\ovec$ and $\xvec$ in the diagram functions, as in the proof
of $\pi^{\sss(0)}(\xvec)\leq P^{\sss(0)}(\xvec)$.  The second sum in
\refeq{tM1unev-bd} can be proved similarly to be bounded by
$\sum_\ell P^{\sss(0)}(\xvec,\btemp^\ell(\yvec))$.  In the first sum in
\refeq{tM1unev-bd}, there are three contributions: (i) $\yvec=\ovec$,
(ii) $b=(\yvec,\xvec)$ and (iii) $\yvec\ne\ovec$ and $\tb\ne\xvec$.
The contribution due to $\yvec=\ovec$ is bounded by
\begin{align}
\lamb\vep L(\ovec,\ovec;\xvec)=\sum_\ell P^{\sss(0)}(\xvec,\bspat^\ell
 (\ovec))-2(\lamb\vep D\sstar\tau)(\xvec)~(\lamb\vep D\sstar\tau\sstar
 \lamb\vep D)(\xvec),
\end{align}
while the contribution due to $b=(\yvec,\xvec)$ is bounded by
\begin{align}
&\mP(\{\ovec\conn\xvec\}\circ\{\ovec\conn\yvec\})~\lamb\vep D(\xvec-\yvec)
 \leq[(\lamb\vep D\sstar\tau)(\xvec)~\tau(\yvec)+(\lamb\vep D\sstar\tau)
 (\yvec)~\tau(\xvec)]~\lamb\vep D(\xvec-\yvec)\nn\\
&\quad=\sum_\ell P^{\sss(0)}(\xvec,\bspat^\ell(\yvec))-[(\lamb\vep D\sstar
 \tau)(\yvec)~(\tau\sstar\lamb\vep D)(\xvec)+(\lamb\vep D\sstar\tau\sstar
 \lamb\vep D)(\xvec)~\tau(\yvec)]~\lamb\vep D(\xvec-\yvec).
\end{align}
We can estimate the case (iii) similarly, and obtain a bound which is
$\lamb\vep L(\ovec,\yvec;\xvec)$ with one of the two $\tau$'s in each
product $\tau\cdot\tau$ in \refeq{Lne-def} replaced by \refeq{bspat}.
Summarizing these bounds, we conclude that the first sum in
\refeq{tM1unev-bd} is bounded by
$\sum_\ell P^{\sss(0)}(\xvec,\bspat^\ell(\yvec))$.  This completes the
proof of \refeq{PN-ind} for $n=0$, and initializes the inductive proof
for $n\geq1$.

To advance the induction hypotheses, we assume that \refeq{PN-ind} holds
for $n=N-1$, and that
\begin{align}
\sum_{b:\bb=\uvec}\mP(M(b,\vvec;\xvec,\yvec))\leq(2\lamb\vep)^{\delta_{
 \uvec,\vvec}-\delta_{\xvec,\yvec}}\sum_\ell L(\uvec,\vvec;\xvec,B^\ell
 (\yvec))\lbeq{M2-bd}
\end{align}
holds, where we write $\sum_\ell L(\uvec,\vvec;\xvec,B^\ell(\xvec))
=2\lamb\vep L(\uvec,\vvec;\xvec)$, similarly to the convention used
in \refeq{PN-def}.  We will prove \refeq{M2-bd} below.  By
\refeq{tE-supset} and the BKR inequality, together with \refeq{PN-def},
\refeq{PN-ind} and \refeq{M2-bd}, we obtain
\begin{align}
\pi^{\sss(N)}(\xvec,\yvec)&\leq\sum_{\uvec,\vvec}\pi^{\sss(N-1)}
 (\uvec,\vvec)\sum_{b_{\sN}:\bb_{\sN}=\uvec}\mP(M(b_{\sN},\vvec;
 \xvec,\yvec))\nn\\
&\leq(2\lamb\vep)^{-\delta_{\xvec,\yvec}}\sum_{\ell,\ell'}\sum_{
 \uvec,\vvec}P^{\sss(N-1)}(\uvec,B^{\ell'}(\vvec))~L(\uvec,\vvec;
 \xvec,B^\ell(\yvec))\nn\\
&=(2\lamb\vep)^{-\delta_{\xvec,\yvec}}\sum_\ell P^{\sss(N)}(\xvec,
 B^\ell(\yvec)).
\end{align}
This advances the induction hypotheses, and hence completes the proof of
\refeq{PN-ind}, assuming that \refeq{M2-bd} holds.

It thus remains to prove \refeq{M2-bd}.  We only consider the case of
$\uvec=\vvec$ and $\xvec=\yvec$, since it explains why the factor
$2\lamb\vep$ is in the definition of the diagram functions.  The case of
$\uvec=\vvec$ and $\xvec\ne\yvec$ can be proved similarly to \refeq{PN-ind}
for $n=0$, and the proof of the remaining case is simpler because extracting
the factors of $\vep$ at $\uvec\ne\vvec$ is unnecessary.  Let $\uvec=\ovec$
in \refeq{M2-bd}, by translation invariance.  Then, the left-hand side of
\refeq{M2-bd} equals the rightmost expression in \refeq{pi0-prebd}, except
for the condition that $b$ is a spatial bond.  By the same observation at
$\ovec$ as in \refeq{pi0-prebd}, we obtain
\begin{align}
\sum_{b:\bb=\ovec}\mP(M(b,\ovec;\xvec,\xvec))&=\sum_{\substack{b
 \text{ spatial}:\\ \bb=\ovec}}\mP(M(b,\ovec;\xvec,\xvec))+\mP
 (M((\ovec,\ovec_+),\ovec;\xvec,\xvec))\nnmb\\
&\leq\sum_{\substack{b\text{ spatial}:\\ \bb=\ovec}}\mP(M(b,\ovec;
 \xvec,\xvec))+\sum_{\substack{b\text{ spatial}:\\ \bb=\ovec}}\mP(\{
 (\ovec,\ovec_+)\conn\xvec\}\circ\{b\conn\xvec\})\nnmb\\
&\leq2\sum_{\substack{b\text{ spatial}:\\ \bb=\ovec}}\mP(M(b,\ovec;
 \xvec,\xvec)) \leq 2\lamb\vep L(\ovec,\ovec;\xvec).
\end{align}
This completes the proof of \refeq{M2-bd}, and hence the proof
of the first inequality in \refeq{BKbds}.

To prove the second inequality in \refeq{BKbds}, we recall that
$\Pi^{\sss(N,n)}(\xvec)=\Pi^{\sss(N;n,1)}(\xvec)+\Pi^{\sss(N;n,2)}(\xvec)$,
where the first and second terms are the contributions to \refeq{Pider-exp}
from $b=b_n$ and from $b\in\piv[\tb_n,\bb_{n+1}]$, respectively.  There are
two $n^{\rm th}$ admissible lines terminating at $\bb_{n+1}$, one from
$\bb_n$ and the other from some vertex $\wvec$.  We can bound
$\Pi^{\sss(N;n,1)}(\xvec)$ by $P^{\sss(N)}(\xvec)$ with the nonzero
admissible line from $\bb_n$, say $(p\sstar \tau)(\bb_{n+1}-\bb_n)$,
replaced by $(\lamb\vep D\sstar \tau)(\bb_{n+1}-\bb_n)$; if $\wvec=\bb_n$,
we replace the factor $p$ in one of the two admissible lines by
$\lamb\vep D$ as above, and add both contributions.  For
$\Pi^{\sss(N;n,2)}(\xvec)$, we use \refeq{tEb-supset=}--\refeq{tEb-supset<}
to obtain the bound
$\sum_\ell\sum_{b\ne b_n}P^{\sss(N)}(\xvec,\bspat^\ell(\bb))$, where
$\sum_\ell$ is the sum over the $n^{\rm th}$ admissible lines.  Together
with the bound on $\Pi^{\sss(N;n,1)}(\xvec)$, we obtain \refeq{tPNn-def}.
This completes the proof the second inequality in \refeq{BKbds}, and hence
the proof of Lemma~\ref{lem:BKbds}.
\end{proof}

\subsection{Estimate of the diagram functions above four dimensions}
\label{ss:diag-est}
In this section, we bound the diagram functions for $d>4$ as follows:

\begin{lem}\label{lem:PtP-bds}
Let $d>4$ and suppose that \refeq{fbd} holds for some $\lamb_0$ and all
$s\leq t$.  Then, there are $\beta_0>0$ and $\Cg<\infty$ such that, for
$\lamb\leq\lamb_0$, $\beta<\beta_0$, $s\in\vep\Zp$ with
$2\vep\leq s\leq t+\vep$, and $q=0,2,4$,
\begin{align}
\sum_x|x|^q\,P_s^{\sss(N)}(x)&\leq\frac{\vep^2(\Cg\beta)^{1\vee N}\sigma^q
 N^{q/2}}{(1+s)^{(d-q)/2}},\qquad\text{for }N\geq0,\lbeq{PN-bd}\\
\sum_x\tilde P_s^{\sss(N;n)}(x)&\leq\frac{\vep^2(\Cg\beta)^N}{(1+s)^{(d-2)
 /2}},\qquad\qquad\text{for }N\geq n\geq1.\lbeq{tPNn-bd}
\end{align}
\end{lem}

Lemma~\ref{lem:piPi-bds}(i) is an immediate consequence of
Lemmas~\ref{lem:BKbds}--\ref{lem:PtP-bds}.  To prove
Lemma~\ref{lem:PtP-bds}, we will use the following three lemmas.

\begin{lem}\label{lem:tau*Dbd}
Assume \refeq{fbd} for $s\leq t=n\vep$ and $\lamb\in I_n$.  Then, there
is a $\Cg=\Cg(d,\lamb)<\infty$ such that the following bounds hold for
$s\leq t$, $q=0,2$ and for that $\lamb$:
\begin{gather}
\sum_x|x|^q\,(\tau_s*D)(x)\leq\Cg\sigma^q(1+s)^{q/2},\lbeq{sumtau*D}\\
\sup_x|x|^q\,(\tau_s*D)(x)\leq\frac{\Cg\sigma^q\beta}{(1+s)^{(d-q)/2}},
 \lbeq{suptau*D}\\
\sup_x|x|^q\,\tau_s(x)\leq(1-\vep)^{s/\vep}\delta_{q,0}+\frac{\Cg\sigma^q
 \beta}{(1+s)^{(d-q)/2}}.\lbeq{suptau*D0}
\end{gather}
\end{lem}

\begin{lem}\label{lem:constr1}
Assume \refeq{fbd} for $s\leq t$.  Let $f_t(x)$ be a diagram function that
satisfies $\sum_xf_t(x)\leq F(t)$ by assigning $l_1$ or $l_\infty$ norm
to each diagram line and using \refeq{fbd} to estimate those norms.
Then, there is a $\Cg=\Cg(d,\lamb)<\infty$ such that
$\sum_xf_t(x,B^\ell(s))\leq\vep\Cg F(t)$ for every $s\leq t$ and
every admissible line $\ell$.
\end{lem}

\begin{lem}\label{lem:conv}
Let $a,b\in\mR$, and let $\kappa$ be a positive number if $a$ or $b$ is 2,
and zero otherwise.  Then, there exists a $C=C(a,b,\kappa)<\infty$ such that
\begin{align}\lbeq{convlem}
\ddsum_{s_1=0}^t\frac{\vep}{(1+s_1)^a}~\ddsum_{s_2=t-s_1}^t\frac{\vep}
 {(1+s_2)^b}\leq\frac{C}{(1+t)^{a\wedge b\wedge(a+b-2)-\kappa}}.
\end{align}
\end{lem}

Lemmas~\ref{lem:tau*Dbd} and \ref{lem:constr1} correspond respectively to
Lemmas~4.3 and 4.6(a) in \cite{hs01}, and Lemma~\ref{lem:conv} corresponds
to \cite[Lemma~3.2]{hs02}.  The result of applying Lemma~\ref{lem:conv} is
the same as \cite[(4.26)]{hs01} when $d>4$ (see also the proof of
Lemma~4.6(b) in \cite{hs01}), but Lemma~\ref{lem:conv} can be applied to
the lower dimensional case as well.  We will use
Lemmas~\ref{lem:tau*Dbd}--\ref{lem:constr1} again in
Sections~\ref{ss:diag-est-lower}--\ref{ss:fc}, and Lemma~\ref{lem:conv} in
Section~\ref{ss:diag-est-lower} and Appendix~\ref{s:adv-ind}.

First, we prove Lemma~\ref{lem:PtP-bds} assuming
Lemmas~\ref{lem:tau*Dbd}--\ref{lem:conv}, and then prove these lemmas.
We will use $c$ to denote a finite positive constant whose exact value
is unimportant and may change from line to line.

\begin{proof}[Proof of Lemma~\ref{lem:PtP-bds} assuming
Lemmas~\ref{lem:tau*Dbd}--\ref{lem:conv}]
For \refeq{PN-bd}, we only consider $q=0$, since the other cases in
\refeq{PN-bd} are proved along the same line of argument as in the last
paragraph of \cite[Section~4.2]{hs01}.

For $s'\leq s<t$, we use Lemma~\ref{lem:tau*Dbd} to obtain
\begin{align}\lbeq{Lsup-bd}
\sup_{u,v}\sum_xL((u,s),(v,s');(x,t))\leq\frac{c\beta\vep}{(1+t-s')^{d/2}}.
\end{align}
Since $P_t^{\sss(0)}(x)=\lamb\vep L(\ovec,\ovec;(x,t))$ for $t\geq2\vep$,
this implies \refeq{PN-bd} with $q=N=0$.  By Lemma~\ref{lem:constr1}, we
also obtain $\sum_\ell\sum_xP_t^{\sss(0)}
(x,B^\ell(s))\leq c\vep[\delta_{0,t}+\vep^2\beta\,(1+t)^{-d/2}]$,
where $\delta_{0,t}$ is the contribution from the first term in
\refeq{Pzero-def}.

For $N\geq1$, we note that, by \refeq{PN-def} we have
\begin{align}\lbeq{PN-prebd}
\sum_xP_t^{\sss(N)}(x)\leq\ddsum_{s,s'}\bigg[\sum_\ell\sum_{u,v}P^{\sss
 (N-1)}((u,s),B^\ell(v,s'))\bigg]\bigg[\sup_{u,v}\sum_xL((u,s),(v,s');
 (x,t))\bigg],
\end{align}
where $\sum_\ell$ is the sum over the $(N-1)^{\rm st}$ admissible lines
in $P_s^{\sss(N-1)}(u)$.  Therefore, $\sum_xP_t^{\sss(1)}(x)$ satisfies
\refeq{PN-bd}, and $\sum_\ell\sum_xP_t^{\sss(1)}(x,B^\ell(s))$ is bounded
by $c\vep^3\Cg\beta\,(1+t)^{-d/2}$.  This initializes the inductive proof
of \refeq{PN-bd} for $N\geq1$ with $q=0$.  Suppose that
$\sum_xP_t^{\sss(N-1)}(x,B^\ell(s))$ is bounded by
$c\vep^3(\Cg\beta)^{N-1}(1+t)^{-d/2}$.  Then, by \refeq{Lsup-bd}
and Lemma~\ref{lem:conv}, $\sum_xP_t^{\sss(N)}(x)$ is bounded
by $\vep^2(\Cg\beta)^N(1+t)^{-d/2}$ if $\Cg$ is sufficiently large.
Note that the factor $\vep^2$ is used in applying Lemma~\ref{lem:conv},
to approximate $\vep\dsum_{s\in\vep\Zp}$ by the Riemann sum.  Using
Lemma~\ref{lem:constr1}, we then obtain
$\sum_xP_t^{\sss(N)}(x,B^\ell(s))\leq c\vep^3(\Cg\beta)^N(1+t)^{-d/2}$.
This completes the proof of \refeq{PN-bd}.

To prove \refeq{tPNn-bd}, we first use Lemma~\ref{lem:constr1} to obtain
$\sum_xP_t^{\sss(N)}(x,\bspat^\ell(s))\leq c\vep^3(\Cg\beta)^N(1+t)^{-d/2}$
for every $s\leq t$, where $\ell$ is an $n^{\rm th}$ admissible line.
Then, we sum the bound over $s\in[0,t]\cap\vep\Zp$ to obtain the desired
bound in \refeq{tPNn-bd}.  Note that the factor $\vep$ is used in an
approximation by the Riemann sum.  This completes the proof.
\end{proof}

\begin{proof}[Proof of Lemma~\ref{lem:tau*Dbd}]
The inequality \refeq{sumtau*D} immediately follows from \refeq{fbd} and
the properties of $D$.  To prove \refeq{suptau*D}--\refeq{suptau*D0}, we use
\begin{align}
(\tau_s*D)(x)\leq(1-\vep)^{s/\vep}D(x)+\lamb\vep\sum_{j=1}^{s/\vep}
 (1-\vep)^{j-1}(D*\tau_{s-j\vep}*D)(x),\lbeq{tau*Dbd}
\end{align}
and
\begin{align}
\tau_s(x)\leq(1-\vep)^{s/\vep}\delta_{o,x}+\lamb s(1-\vep)^{s/\vep-1}
 D(x)+\lamb^2\vep\sum_{j=2}^{s/\vep}(j-1)\vep\,(1-\vep)^{j-2}
 (D*\tau_{s-j\vep}*D)(x).\lbeq{tau*D1bd}
\end{align}
where $(1-\vep)^n$ is the probability that $(o,0)$ is connected to
$(o,n\vep)$ along the temporal axis.  Since the first term in
\refeq{tau*Dbd} and the second term in \refeq{tau*D1bd}, multiplied by
$|x|^q$, are both bounded by $c\sigma^q\beta\,(1+s)^{-(d-q)/2}$ for any
$x\in\Zd$, we only need to consider the last terms in
\refeq{tau*Dbd}--\refeq{tau*D1bd}.

We fix $r\in(0,1)$ and use \refeq{fbd} to bound the second term in
\refeq{tau*Dbd} by
\begin{align}\lbeq{separating}
\lamb\vep\sum_{j=1}^{rs/\vep}\frac{K\beta\,(1-\vep)^{j-1}}{(1+s-j
 \vep)^{d/2}}+\lamb\vep\sum_{j=rs/\vep}^{s/\vep}\frac{K\beta\,(1-
 \vep)^{j-1}}{(1+s-j\vep)^{d/2}}\leq\frac{c\beta}{(1+s)^{d/2}}.
\end{align}
By the same argument, the third term in \refeq{tau*D1bd} can be bounded
by $c\beta\,(1+s)^{-d/2}$.  This completes the proof of \refeq{suptau*D}--\refeq{suptau*D0} for
$q=0$.

For $|x|^2$ times the second term in \refeq{tau*Dbd} or the third term
in \refeq{tau*D1bd}, we have
\begin{align}
|x|^2(D*\tau_s*D)(x)&\leq2\sum_y(|y|^2+|x-y|^2)\,(D*\tau_{s/2})(y)\,
 (\tau_{s/2}*D)(x-y)\nn\\
&\leq4\|D*\tau_{s/2}\|_\infty\sum_x|x|^2(\tau_{s/2}*D)(x).
\end{align}
Applying \refeq{tau*Dbd} to $\|D*\tau_{s/2}\|_\infty$, using \refeq{fbd},
and then separating the sum over $j$ as in \refeq{separating}, we obtain
\refeq{suptau*D}--\refeq{suptau*D0} for $q=2$.  This completes the
proof of Lemma~\ref{lem:tau*Dbd}.
\end{proof}

\begin{proof}[Proof of Lemma~\ref{lem:constr1}]
By the convention used in \refeq{PN-def}, for $s=t$ we have
\begin{align}
\sum_{x,y}f_t(x,B^\ell(y,t))=\sum_xf_t(x,B^\ell(x,t))\leq\sum_x\sum_\ell
 f_t(x,B^\ell(x,t))=\sum_x2\lamb\vep f_t(x)\leq2\lamb\vep F(t),
\end{align}
where $\sum_\ell$ is the sum over the admissible lines arriving at $(x,t)$.
For $s<t$, Construction~$B^\ell(s)$ replaces the diagram line $\ell$, say
$\tau_t(x)$, by $\lamb\vep(\tau_s*D*\tau_{t-s-\vep})(x)+\lamb\vep(1-\vep)
(\tau_s*\tau_{t-s-\vep})(x)$.  By Lemma~\ref{lem:tau*Dbd} and \refeq{fbd},
we obtain
\begin{align}
&\lamb\vep\|\tau_s*D*\tau_{t-s-\vep}\|_1=\lamb\vep\|\tau_s*D\|_1\,\|\tau_{t
 -s-\vep}\|_1\leq\lamb\vep C_{\sss\ref{lem:tau*Dbd}}K,\\
&\lamb\vep\|\tau_s*D*\tau_{t-s-\vep}\|_\infty\leq\lamb\vep\|\tau_{s\vee(t-s
 -\vep)}*D\|_\infty\,\|\tau_{s\wedge(t-s-\vep)}\|_1\leq\lamb\vep\frac{2^{d/
 2}C_{\sss\ref{lem:tau*Dbd}}K}{(1+t)^{d/2}},
\end{align}
where $C_{\sss\ref{lem:tau*Dbd}}$ is the constant in Lemma~\ref{lem:tau*Dbd},
and we use $s\vee(t-s-\vep)\geq t/2$ to obtain the last inequality.  The
$l_1$ and $l_\infty$ norms of $\lamb\vep(1-\vep)(\tau_s*\tau_{t-s-\vep})(x)$
can be estimated in the same way.  Therefore, the effect of
Construction~$B^\ell(s)$ is to obtain, at worst, an additional constant
$\Cg\vep=2^{1+d/2}C_{\sss\ref{lem:tau*Dbd}}\lamb\vep$ in a bound.  This
completes the proof.
\end{proof}

\begin{proof}[Proof of Lemma~\ref{lem:conv}]
We prove \refeq{convlem} for $a\wedge b\geq0$ and for $a\wedge b<0$
separately.

Let $a\wedge b\geq0$.  Separating the sum over $s_1$ into
$\dsum_{0\leq s_1\leq t/2}$ and $\dsum_{t/2<s_1\leq t}$, and using
$s_2\geq t-s_1\geq t/2$ in the former sum, we can bound the left-hand side
of \refeq{convlem} by
\begin{align}\lbeq{convab}
\frac{c}{(1+t)^b}\ddsum_{s_1=0}^{t/2}\frac{\vep}{(1+s_1)^{a-1}}+\frac{c}
 {(1+t)^a}\ddsum_{s_1=t/2}^t~\ddsum_{s_2=t-s_1}^t\frac{\vep^2}{(1+s_2)^b},
\end{align}
where the first term is bounded by $c'(1+t)^{-b-(a-2)\wedge0+\kappa_a}$,
where $\kappa_a$ is an arbitrarily small but positive number if $a=2$,
otherwise $\kappa_a=0$.  Also, the double sum in \refeq{convab} is
\begin{align}
\ddsum_{s_2=0}^t\frac{\vep}{(1+s_2)^b}~\ddsum_{s_1=t/2\vee(t-s_2)}^t\vep~
 \leq~\sum_{s_2=0}^t\frac{\vep}{(1+s_2)^{b-1}}~\leq~\frac{c}{(1+t)^{(b-2)
 \wedge0-\kappa_b}},
\end{align}
where $\kappa_b$ is an arbitrarily small but positive number if $b=2$,
otherwise $\kappa_b=0$.  This completes the proof of \refeq{convlem} for
$a\wedge b\geq0$.

Next, we consider the case $a\wedge b<0$.  Due to the symmetry of the
left-hand side of \refeq{convlem} in terms of $s_1,s_2$, we suppose
$b<a\wedge0$.  Then, we use the trivial inequality
$(1+s_2)^{-b}\leq(1+t)^{-b}$.  The remaining term equals
$\vep\dsum_{0\leq s_1\leq t}s_1\,(1+s_1)^{-a}$ and is bounded by
$c\,(1+t)^{-(a-2)\wedge0+\kappa_a}$.  This completes the proof of
\refeq{convlem} for $a\wedge b<0$, and hence the proof of
Lemma~\ref{lem:conv}.
\end{proof}

\subsection{Estimate of the diagram functions at and below four dimensions}
\label{ss:diag-est-lower}
In this section, we bound the diagram functions for $d\leq4$
by using their inductive construction in \refeq{Pzero-def}
and \refeq{PN-def}--\refeq{tPNn-def} as well as
Lemmas~\ref{lem:tau*Dbd}--\ref{lem:conv}, as in the proof of
Lemma~\ref{lem:PtP-bds} in Section~\ref{ss:diag-est},  but we replace
$\sigma$ and $\beta$ in Lemma~\ref{lem:tau*Dbd} by $\sigma_{\sT}$
and $\beta_{\sT}$, respectively.

Lemma~\ref{lem:piPi-bds}(ii) is an immediate consequence of
Lemma~\ref{lem:BKbds} and the following lemma:

\begin{lem}\label{lem:PtP-bds-lower}
Let $d\leq4$ with $\alpha=bd-\frac{4-d}2>0$, $\mu\in(0,\alpha)$ and
$t\leq T\log T$, and suppose that \refeq{fbd} holds for some $\lamb_0$ and
all $s\leq t$.  Then, there are $\beta_0>0$ and $\Cg<\infty$ such that, for
$\lamb\leq\lamb_0$, $\beta_1<\beta_0$, $s\in\vep\Zp$ with
$2\vep\leq s\leq t+\vep$, and $q=0,2,4$,
\begin{align}
\sum_x|x|^q\,P_s^{\sss(N)}(x)&\leq\frac{\vep^2\Cg\beta_{\sT}(\Cg\hat
 \beta_{\sT})^{0\vee(N-1)}\sigma_{\sT}^qN^{q/2}}{(1+s)^{(d-q)/2}},
 \qquad\text{for }N\geq0,\lbeq{PN-bd-lower}\\
\sum_x\tilde P_s^{\sss(N;n)}(x)&\leq\frac{\vep^2\Cg\beta_{\sT}(\Cg\hat
 \beta_{\sT})^{N-1}}{(1+s)^{(d-2)/2}},\qquad\qquad\qquad\text{ for }N
 \geq n\geq1,\lbeq{tPNn-bd-lower}
\end{align}
where $\beta_{\sT}=\beta_1T^{-bd}$ and $\hat\beta_{\sT}=\beta_1T^{-\mu}$.
\end{lem}

\begin{proof}
The proof is almost the same as that of Lemma~\ref{lem:PtP-bds}.  The
only difference arises when we apply Lemma~\ref{lem:conv}.  Let $N\geq1$
and suppose that the quantity in the first brackets in \refeq{PN-prebd} is
bounded by $c\vep^3\Cg\beta_{\sT}(C\hat\beta_{\sT})^{N-1}(1+t)^{-d/2}$,
where $\hat\beta_{\sT}=\beta_1T^{-\mu}$ with $\mu\in(0,\alpha)$.  Then, by
Lemma~\ref{lem:conv} and \refeq{Lsup-bd} with $\beta$ replaced by
$\beta_{\sT}$, the right-hand side of \refeq{PN-prebd} for $d\leq4$ is
bounded by
\begin{align}\lbeq{conv-result-lower}
\frac{c\vep^2\Cg\beta_{\sT}^2(\Cg\hat\beta_{\sT})^{N-1}}{(1+t)^{\frac{d}2
 \wedge(d-2)-\kappa}}=c\beta_1T^{-bd}(1+t)^{(4-d)/2+\kappa}\frac{\vep^2
 \Cg\beta_{\sT}(\Cg\hat\beta_{\sT})^{N-1}}{(1+t)^{d/2}},
\end{align}
where $\kappa$ is an arbitrarily small but positive number if $d=4$,
otherwise $\kappa=0$.  Since $t\leq T\log T$ and
$-bd+\frac{4-d}2=-\alpha<-\mu$, the factor in front of the fraction in the
right-hand side is bounded by $\Cg\hat\beta_{\sT}$ if $\Cg$ is sufficiently
large, and thus we obtain \refeq{PN-bd-lower} with $q=0$.

For \refeq{PN-bd-lower} with $q=2$ and \refeq{tPNn-bd-lower}, we use
Lemma~\ref{lem:conv} as in \refeq{conv-result-lower}, with
$(a,b)=(\frac{d}2,\frac{d-2}2)$, to obtain the factor
\begin{align}
\frac{\beta_{\sT}}{(1+t)^{(d-2)/2\wedge(d-3)-\kappa}}
 =\frac{\beta_1T^{-bd}(1+t)^{(4-d)/2+\kappa}}{(1+t)^{(d-2)/2}}
 \leq\frac{c\hat\beta_{\sT}}{(1+t)^{(d-2)/2}}.
\end{align}
To prove \refeq{PN-bd-lower} for $q=4$, we apply Lemma~\ref{lem:conv}
as above, with $(a,b)=(\frac{d}2,\frac{d-4}2)$ and
$(\frac{d-2}2,\frac{d-2}2)$.  This completes the proof.
\end{proof}

\subsection{Finite containment}\label{ss:fc}
In this section, we prove that $\pi_t^{\sss(N)}(x)$ can be approximated by
\begin{align}
\pi_t^{\sss(N)}(x\,|\,R)&=\sum_{\vec b_{\sN}}\mP\big(\tilde E_{\vec b_N}
 ^{\sss(N)}(\xvec)\cap\{\bC_{[0,t]}\subset\Box_{\sR}\}\big),\lbeq{pi-fc}
\end{align}
where $\xvec=(x,t)$, $\Box_{\sR}=[-R,R]^d\cap\Zd$ and
\begin{align}\lbeq{bC-gen}
\bC_{[0,t]}=\bigcup_{s=0}^t\bC_s(\ovec),&&
\bC_s(\yvec)=\{z\in\Zd:\yvec\conn(z,s)\}.
\end{align}
We will also use the abbreviation $\bC_s=\bC_s(\ovec)$.  More precisely,
we prove below that
\begin{align}\lbeq{fincont}
\pi_t^{\sss(N)}(x)=\pi_t^{\sss(N)}(x\,|\,R)+o(1)\,\vep^2,
\end{align}
where $o(1)$ is independent of $\vep$ and decays to zero as $R\to\infty$,
by using the estimates for the diagram functions in
Sections~\ref{ss:diag-bds}--\ref{ss:diag-est-lower}.  This is a refined
version of the finite containment argument used in proving the uniformity
of \refeq{uniformconv}, and will be useful in dealing with the continuum
limit in Section~\ref{s:continuum}.

\begin{proof}
First, we note that
\begin{align}
0\leq\pi_t^{\sss(N)}(x)-\pi_t^{\sss(N)}(x\,|\,R)&\leq\sum_{\vec b_{\sN}}
 \ddsum_{s=\vep}^t\mP\big(\tilde E_{\vec b_N}^{\sss(N)}(\xvec)\cap\{
 \bC_{[0,s-\vep]}\subset\Box_{\sR},~\bC_s\not\subset\Box_{\sR}\}\big)\nn\\
&\leq\sum_{\vec b_{\sN}}\ddsum_{s=\vep}^t\sum_{\substack{u\in\Box_R\\
 v\notin\Box_R}}\mP\big(\tilde E_{\vec b_N}^{\sss(N)}(\xvec)\cap\{\ovec
 \conn((u,s-\vep),(v,s))\}\big).
\end{align}
Note that $((u,s-\vep),(v,s))$ is a spatial bond, since $u\in\Box_{\sR}$
and $v\notin\Box_{\sR}$.  The event in the rightmost expression is included
in the union of two events: (i) $\ovec\conn((u,s-\vep),(v,s))\conn\xvec$,
or (ii) there is a vertex $(w,r)\in\Box_{\sR}\times[0,s)$ such that
$\ovec\conn(w,r)\conn\xvec$ and that $(w,r)\conn((u,s-\vep),(v,s))$
disjointly from $\tilde E_{\vec b_{\sN}}^{\sss(N)}(\xvec)$.
The contribution from the case (i) and from the case (ii) with
$(w,r)=(u,s-\vep)$ can be bounded by
$\sum_\ell\dsum_s\sum_{v\notin\Box_{\sR}}P^{\sss(N)}(\xvec,B^\ell(v,s))$,
where $\sum_\ell$ is the sum over all admissible lines (i.e., the sum over
$n=1,\dots,N$ of the sum over the $n^{\rm th}$ admissible lines), and where
we modified Construction~$\bspat^\ell(v,s)$ by fixing the second endpoint
$(v,s)$, instead of fixing the first endpoint $(u,s-\vep)$ as defined in
\refeq{bspat}.  This bound, divided by $\vep^2$, decays as $R\to\infty$
uniformly in $\vep$, since the sum over $x\in\Zd$ of the unrestricted
sum $\sum_\ell\dsum_sP^{\sss(N)}(\xvec,B^\ell(s))$ is bounded, by using
Lemma~\ref{lem:constr1}, by $c\vep^2t\,(1+t)^{-d/2}$.

The contribution from the case (ii) with $(w,r)\ne(u,s-\vep)$
can be bounded by
\begin{align}\lbeq{fincont-ii}
\sum_\ell\ddsum_{s=\vep}^t\ddsum_{r=0}^{s-\vep}\sum_{w\in\Box_R}P^{\sss(N)}
 (\xvec,B^\ell(w,r))\sum_{v\notin\Box_R}(\tau_{s-\vep-r}*\lamb\vep D)(v-w),
\end{align}
where we relabelled the second endpoint of the spatial bond in \refeq{btemp}
as $(w,r)$.  For $w\in\Box_{\sss R/2}$, we use
\begin{align}
&\sum_{w\in\Box_{R/2}}P^{\sss(N)}(\xvec,B^\ell(w,r))\sum_{v\notin\Box_R}
 (\tau_{s-\vep-r}*D)(v-w)\\
&\quad\leq P^{\sss(N)}(\xvec,B^\ell(r))\sup_{w\in\Box_{R/2}}\sum_{v\notin
 \Box_R}(\tau_{s-\vep-r}*D)(v-w)\leq P^{\sss(N)}(\xvec,B^\ell(r))\sum_{z
 \notin\Box_{R/2}}(\tau_{s-\vep-r}*D)(z),\nn
\end{align}
and for $w\in\Box_{\sR}\setminus\Box_{\sss R/2}$, we use
\begin{align}
\sum_{w\notin\Box_{R/2}}P^{\sss(N)}(\xvec,B^\ell(w,r))\sum_{v\notin\Box_R}
 (\tau_{s-\vep-r}*D)(v-w)\leq\|\tau_{s-\vep-r}*D\|_1\sum_{w\notin\Box_{R/2}}
 P^{\sss(N)}(\xvec,B^\ell(w,r)).
\end{align}
By Lemma~\ref{lem:tau*Dbd}, we have $\|\tau_{s-\vep-r}*D\|_1\leq\Cg$, so
that $\sum_{z\notin\Box_{R/2}}(\tau_{s-\vep-r}*D)(z)$ decays to zero as
$R\to\infty$, independently of $\vep$.  In addition, by
Lemma~\ref{lem:constr1}, we have
$\sum_\ell P^{\sss(N)}(\xvec,B^\ell(r))\leq c\vep^3(1+t)^{-d/2}$, so that
$\sum_\ell\sum_{w\notin\Box_{R/2}}P^{\sss(N)}(\xvec,B^\ell(w,r))\leq o(1)
\,\vep^3(1+t)^{-d/2}$.  Since there is another factor of $\vep$ in the
summand of \refeq{fincont-ii}, while there are two summations over
$\vep\Zp$ in \refeq{fincont-ii}, we conclude that \refeq{fincont-ii} is
$o(1)\,\vep^2t^2(1+t)^{-d/2}$.  This completes the proof of
\refeq{fincont}.
%
\end{proof}


\newcommand{\kbar}{k}

\section{Inductive argument}\label{s:ind}
In this section, we prove Proposition~\ref{thm-disc} by applying the inductive
method of \cite{hs02} for self-avoiding walk in $\Zd$ and for oriented
percolation in $\Zd\times\Zp$, to the recursion equation \refeq{tau-fourier}
for oriented percolation in $\Zd\times\vep\Zp$.  To consider the case of
$\vep\ll1$, we will modify the induction hypotheses to incorporate the
dependence on $\vep$.  We expect that a similar method could be used for
{\it continuous-time} weakly self-avoiding walk above its upper critical
dimension.

First, we consider $d>4$ in Sections~\ref{ss:indhyp}--\ref{ss:discr-thm}.
In Section~\ref{ss:indhyp}, we give the modified version of the induction
hypotheses.  In Section~\ref{ss:cons-ind}, we show several consequences
of the induction hypotheses, mainly the bounds in \refeq{fbd}.  In
Section~\ref{ss:discr-thm}, we prove Proposition~\ref{thm-disc}.  We
complete this section by proving the results for $d\leq4$ in
Section~\ref{ss:lowdim}.  Since a similar strategy applies to the
lower-dimensional case, we only discuss the necessary changes.

\subsection{Induction hypotheses}\label{ss:indhyp}
\subsubsection{General assumptions}\label{sss:genass}
In Section~\ref{ss:exp-two}, we derived the recursion equation
\refeq{tau-2exp} for the two-point function.  Taking the Fourier transform
with respect to the spatial component, we obtain \refeq{tauk}, i.e.,
\begin{align}\lbeq{tau-fourier}
\wtau_{t+\vep}(k)=\ddsum_{s=0}^t\wpi_s(k)\,\wpr(k)\,\wtau_{t-s}(k)
 +\wpi_{t+\vep}(k).
\end{align}
The probability distribution $D:\Zd\mapsto[0,1]$ satisfies the assumptions
in Section~\ref{ss:results}.  In addition, we assume that there exists
$\eta>0$ such that
\begin{align}
a(k)\equiv1-\wD(k)\begin{cases}
  \asymp L^2|k|^2,&\mbox{if }\|k\|_\infty\leq L^{-1},\\
  >\eta,&\mbox{if }\|k\|_\infty>L^{-1},
 \end{cases}&&
a(k)<2-\eta\quad\forall k\in[-\pi,\pi]^d,
\end{align}
where $a\asymp b$ means that the ratio $a/b$ is bounded away from zero and
infinity.  These assumptions correspond to Assumption~D in \cite{hs02}.

If we replace $t$ by $n\vep$, and write
\begin{align}\lbeq{indchoices}
f_n(k)=\wtau_{n\vep}(k),&&
e_n(k)=\wpi_{n\vep}(k),&&
g_{n+1}(k)=\wpi_{n\vep}(k)~\wpr(k),
\end{align}
where the dependence on $\lamb$ is left implicit, then \refeq{tau-fourier}
equals
\begin{align}\lbeq{f-renew}
f_{n+1}(k)=\sum_{m=0}^ng_{m+1}(k)\,f_{n-m}(k)+e_{n+1}(k),
\end{align}
with $f_0(k)=1$.  This is equivalent to the recursion relation (1.1) in
\cite{hs02}.  The only difference is
\begin{align}\lbeq{f1g1-def}
f_1(k)=g_1(k)=\wpr(k)=1-\vep+\lamb\vep\wD(k),
\end{align}
whereas in \cite{hs02}, $f_1(k)=g_1(k)=\lamb\wD(k)$.  This change leads to
a modification of the induction hypotheses in \cite{hs02}, the main reason
being that we need to prove uniformity in $\vep$.  Further technical
changes are explained in Section~\ref{sss:indhyp}.

\subsubsection{Statement of the induction hypotheses}\label{sss:indhyp}
Fix $\gamma$, $\delta$ and $\rho$ according to
\begin{align}\lbeq{drgd}
-(2+\rho)<0<\tfrac{d}2-(2+\rho)<\gamma<\gamma+\delta<1\wedge\Delta\wedge
 \tfrac{d-4}2.
\end{align}
We define $\lamb_n$ recursively by $\lamb_0=\lamb_1=1$ and, for $n\geq2$,
\begin{align}\lbeq{lambn-def}
\lamb_n=1-\frac1\vep\sum_{l=2}^ng_l(0;\lamb_{n-1}),
\end{align}
where we explicitly write the dependence on $\lamb_{n-1}$ of $g_l(0)$.  Let
\begin{align}\lbeq{In-def}
I_n=\lamb_n+\frac{K_1\beta}{(1+n\vep)^{(d-2)/2}}[-1,1],
\end{align}
and define $v_n=v_n(\lamb)$ by $v_0=v_1=\lamb$ and, for $n\geq2$,
\begin{align}\lbeq{vn-def}
v_n=\frac{\lamb-\frac1{\sigma^2\vep}\sum_{l=2}^n\nabla^2g_l(0)}{1+\sum_{l=2}^n
 (l-1)\,g_l(0)}.
\end{align}
Let $K_1,\ldots,K_5$ be some positive and finite constants that are independent of
$\beta$ and $\vep$, and are related by
\begin{align}\lbeq{Krel}
K_3\gg K_1\gg K_4\gg1,\qquad\qquad K_2,K_5\gg K_4.
\end{align}

The induction hypotheses are that the following (H1)--(H4) hold for all
$\lamb\in I_n$ and $m=1,\dots,n$.
\begin{enumerate}
\item[]{\bf(H1)--(H2)}
\begin{align}\lbeq{H12}
|\lamb_m-\lamb_{m-1}|\leq\frac{\vep K_1\beta}{(1+m\vep)^{d/2}},&&
|v_m-v_{m-1}|\leq\frac{\vep K_2\beta}{(1+m\vep)^{(d-2)/2}},
\end{align}
\item[]{\bf(H3)}\quad
For $k\in\cA_m\equiv\{k:a(k)\leq\gamma\frac{\log(2+m\vep)}{1+m\vep}\}$,
$f_m(k)$ can be written in the form
\begin{align}\lbeq{fm-expr}
f_m(k)=\prod_{l=1}^m[1-\vep v_l\,a(k)+\vep\,r_l(k)],
\end{align}
where $r_l(k)$ obeys the bounds
\begin{align}\lbeq{H3}
|r_l(0)|\leq\frac{K_3\beta}{(1+l\vep)^{(d-2)/2}},&&
|r_l(k)-r_l(0)|\leq\frac{K_3\beta\,a(k)}{(1+l\vep)^\delta}.
\end{align}
\item[]{\bf(H4)}\quad
For $k\notin\cA_m$, $f_m(k)$ obeys the bounds
\begin{align}\lbeq{H4}
|f_m(k)|\leq\frac{K_4\,a(k)^{-2-\rho}}{(1+m\vep)^{d/2}},&&
|f_m(k)-f_{m-1}(k)|\leq\frac{\vep K_5\,a(k)^{-1-\rho}}{(1+m\vep)^{d/2}}.
\end{align}
\end{enumerate}

Instead of \refeq{fm-expr}, we can alternatively write $f_m(k)$ as
\begin{align}\lbeq{fm-reexpr}
f_m(k)=f_m(0)\prod_{l=1}^m[1-\vep v_l\,a(k)+\vep\,s_l(k)],
\end{align}
where
\begin{align}\lbeq{f0sl-def}
f_m(0)=\prod_{l=1}^m[1+\vep\,r_l(0)],&&
s_l(k)=\frac{\vep v_l\,r_l(0)\,a(k)+[r_l(k)-r_l(0)]}{1+\vep\,r_l(0)}.
\end{align}
The induction hypothesis (H3) implies
\begin{align}\lbeq{sl-bd}
|s_l(k)|\leq\frac{\vep v_l|r_l(0)|a(k)+|r_l(k)-r_l(0)|}{1-\vep|r_l(0)|}\leq
 \frac{(1+\vep v_l)K_3\beta\,a(k)}{(1-\vep K_3\beta)(1+l\vep)^\delta}.
\end{align}
In some cases, we will use \refeq{fm-reexpr}--\refeq{sl-bd}, instead of \refeq{fm-expr}--\refeq{H3}.  Moreover, by \refeq{fm-reexpr} and spatial
symmetry, we obtain
\begin{align}\lbeq{fmdiff-bd}
\nabla^2 f_m(0)=f_m(0)\;\vep\sum_{l=1}^m[-v_l\sigma^2
 +\nabla^2s_l(0)\big].
\end{align}

The advancement of the induction hypotheses is a small modification of
that in \cite{hs02}, which we add to keep this paper self-contained.
The advancement is deferred to Appendix~\ref{s:adv-ind}.

\subsection{Consequences of the induction hypotheses}\label{ss:cons-ind}
We assume $\beta\ll1$ and use $c$ to denote a positive and finite constant
that may depend on $d,\gamma,\delta,\rho$, but not on $K_i,k,n,\beta,\vep$.
The value of $c$ may change from line to line.

The following four lemmas, corresponding respectively to
\cite[Lemmas~2.1, 2.2, 2.4, 2.3]{hs02}, are consequences
of the induction hypotheses (H1)--(H4) for $d>4$.

\begin{lem}\label{lem:I}
Assume (H1) for $m=1,\dots,n$.
Then, $I_0\supset I_1\supset\cdots\supset I_n$.
\end{lem}

\begin{lem}\label{lem:fabs}
Let $\lamb\in I_n$ and assume (H2)--(H3) for $m=1,\dots,n$.
For $k\in\cA_m$,
\begin{align}
|f_m(k)|\leq e^{cK_3\beta}e^{-m\vep[1-c(K_1+K_2+K_3)\beta]\,a(k)}.
\end{align}
\end{lem}

\begin{lem}\label{lem:fdiff}
Let $\lamb\in I_n$ and assume (H2)--(H3) for $m=1,\dots,n$.  Then,
\begin{align}
|\nabla^2f_m(0)|\leq[1+c\,(K_1+K_2+K_3)\beta]\sigma^2m\vep.
\end{align}
\end{lem}

\begin{lem}\label{lem:fD2}
Let $\lamb\in I_n$ and assume (H2)--(H4) for $m=1,\dots,n$.  Then,
\begin{align}
\|\wD^2f_m\|_1\leq\frac{c\,(1+K_4)\beta}{(1+m\vep)^{d/2}}.
\end{align}
\end{lem}

The bounds \refeq{fbd} for $s\leq n\vep$ follow from
Lemmas~\ref{lem:fabs}--\ref{lem:fD2}, if $K\gg K_4$.  The proofs of
Lemmas~\ref{lem:I}--\ref{lem:fD2} are almost identical to those of
\cite[Lemmas~2.1--2.4]{hs02}, and are deferred to Appendix~\ref{s:adv-ind}.

By Lemma~\ref{lem:I}, if $\lamb\in I_m$ for some $m\geq0$,
then $\lamb\in I_0$ and hence, by \refeq{In-def},
\begin{align}\lbeq{lamb-bd}
|\lamb-1|\leq K_1\beta.
\end{align}
It also follows that $I_\infty=\bigcap_{m=0}^\infty I_m$ is a singleton
$\lamb_\infty$.  As discussed in \cite[Theorem 1.2]{hs02}, we obtain
$\lamb_\infty=\lambce$.  Moreover, it follows from the second inequality
of \refeq{H12} that, for $\lamb\in I_m$,
\begin{align}\lbeq{vm-bd}
|v_m-1|\leq\sum_{l=1}^m|v_l-v_{l-1}|+|v_0-1|\leq\sum_{l=1}^m
 \frac{\vep K_2\beta}{(1+l\vep)^{(d-2)/2}}+|\lamb-1|\leq(cK_2+K_1)\beta.
\end{align}
We note that the factor $\vep$ in the numerator is necessary to approximate
the sum by the Riemann sum when $\vep$ is small.  The factors of $\vep$ in
\refeq{H12}--\refeq{H4} are incorporated for the same reason.

\subsection{Proof of Proposition~\ref{thm-disc}}\label{ss:discr-thm}
Fix $\lamb=\lambce$, so that the induction hypotheses (H1)--(H4) and
Lemmas~\ref{lem:I}--\ref{lem:fD2} hold for all $m\in\mN$.  From now
on, we suppress the dependence on $\vep$ and write $\lambc=\lambce$,
$A=A^{\sss(\vep)}$ and $v=v^{\sss(\vep)}$.

Note that, by \refeq{H12}--\refeq{H3}, we have that, for $n<m$,
\begin{align}
|f_n(0)-f_m(0)|&=\prod_{l=1}^n[1+\vep r_l(0)]\;\bigg|1-\prod_{l=n+1}^m
 [1+\vep r_l(0)]\bigg|\leq\frac{cK_3\beta}{(1+n\vep)^{(d-4)/2}},
 \lbeq{fn-fm}\\
|v_n-v_m|&\leq\sum_{l=n+1}^m|v_l-v_{l-1}|\leq\frac{cK_2\beta}{(1+n\vep)
 ^{(d-4)/2}},\lbeq{vn-vm}
\end{align}
so that $\{f_n(0)\}_{n=1}^\infty$ and $\{v_n\}_{n=1}^\infty$ are Cauchy
sequences.  Therefore, the limits $A=\lim_{n\to\infty}f_n(0)$ and
$v=\lim_{n\to\infty}v_n$ exist, and satisfy
\begin{align}\lbeq{Av}
|f_n(0)-A|\leq\frac{cK_3\beta}{(1+n\vep)^{(d-4)/2}},&&
|v_n-v|\leq\frac{cK_2\beta}{(1+n\vep)^{(d-4)/2}}.
\end{align}
In particular, by Lemma~\ref{lem:fabs} and \refeq{vm-bd}, both $A$ and $v$
are equal to $1+O(\beta)$.

Let $t=n\vep$ and $\tilde k=\frac{k}{\sqrt{v\sigma^2t}}\in\cA_n$.
By \refeq{Deltadef},
$a(\tilde k)=\frac{|k|^2}{2dvt}+O(|k|^{2+2\Delta}t^{-1-\Delta})$.
Using \refeq{fm-reexpr}, \refeq{sl-bd}, \refeq{Av} and
$\delta<1\wedge\frac{d-4}2$, we obtain
\begin{align}
f_n(\tilde k)&=\bigg[A+\frac{O(\beta)}{(1+t)^{(d-4)/2}}\bigg]\prod_{l=1}^n
 \bigg[1-\vep\bigg[v+\frac{O(\beta)}{(1+l\vep)^\delta}\bigg]a(\tilde k)
 \bigg]\nn\\
&=\bigg[A+\frac{O(\beta)}{(1+t)^{(d-4)/2}}\bigg]\bigg[1-\frac{vt\,
 a(\tilde k)}n\bigg]^n\prod_{l=1}^n\bigg[1-\frac{\vep O(\beta)\,a(\tilde k)}
 {(1+l\vep)^\delta}\bigg]\nn\\
&=\bigg[A+\frac{O(\beta)}{(1+t)^{(d-4)/2}}\bigg]e^{-\frac{|k|^2}{2d}
 +O(|k|^{2+2\Delta}t^{-\Delta})+O(\vep|k|^4t^{-1})}\bigg[1+\vep
 \sum_{l=1}^{t/\vep}\frac{O(\beta|k|^2t^{-1})}{(1+l\vep)^\delta}\bigg]\nn\\
&=Ae^{-\frac{|k|^2}{2d}}\bigg[1+\frac{O(\beta)}{(1+t)^{(d-4)/2}}+O(|k|^{2+2
 \Delta}t^{-\Delta})+O(\vep|k|^4t^{-1})+\frac{O(\beta|k|^2)}{(1+t)^\delta}
 \bigg],\lbeq{11cpf}
\end{align}
where the last error term follows from
\begin{align}
\vep\sum_{l=1}^{t/\vep}\frac{O(t^{-1})}{(1+l\vep)^\delta}=O(t^{-1})
 \big[(1+t)^{1-\delta}-1\big]=O((1+t)^{-\delta}),\lbeq{error1}
\end{align}
for $\delta<1$.  Using $\frac{|k|^2}t\leq c\frac{\log(2+t)}{1+t}$
for small $\tilde k\in\cA_n$ and $\delta<1\wedge\Delta$, we have
\begin{align}
O(|k|^{2+2\Delta}t^{-\Delta})&\leq O(|k|^2)\Big[\frac{\log(2+t)}
 {1+t}\Big]^\Delta\leq\frac{O(|k|^2)}{(1+t)^\delta},\lbeq{error2}\\
O(\vep|k|^4t^{-1})&\leq O(\vep|k|^2)\frac{\log(2+t)}{1+t}\leq
 \frac{O(\vep|k|^2)}{(1+t)^\delta}.\lbeq{error3}
\end{align}
By \refeq{11cpf}--\refeq{error3}, we obtain \refeq{tauasyvep}.

Let $e_1,\ldots,e_d$ denote the standard basis vectors in $\Rd$.
Then, by \refeq{sl-bd},
\begin{align}\lbeq{sldiff-bd}
|\nabla^2s_l(0)|=\bigg|\sum_{i=1}^d\lim_{h\to0}\frac{s_l(he_i)
 -s_l(0)}{h^2}\bigg|\leq\frac{cK_3\beta}{(1+l\vep)^\delta}
 \bigg|\sum_{i=1}^d\lim_{h\to0}\frac{a(he_i)}{h^2}\bigg|=
 \frac{cK_3\sigma^2\beta}{(1+l\vep)^\delta}.
\end{align}
Since $\delta<1\wedge\frac{d-4}2$, it follows from \refeq{fmdiff-bd},
\refeq{Av}, \refeq{error1} and \refeq{sldiff-bd} that
\begin{align}
-\frac{\nabla^2f_n(0)}{f_n(0)}=v\sigma^2t\,[1+O(\beta)\,(1+t)^{-\delta}],
\end{align}
which is \refeq{taugyrvep}.

The upper bound in \refeq{tausupvep} is an immediate consequence of
\refeq{suptau*D0}.  For the lower bound, we consider the case of $t\geq1$
and the case of $t<1$, separately.  When $t\geq1$, we follow the proof of \cite[Corollary~1.4]{hs02} for
oriented percolation.  In this case, we use \refeq{tauasyvep} and obtain
the lower bound of the form $cL^{-d}t^{-d/2}$.  When $t<1$, we use the
trivial inequality
\begin{align}
\|\tau_t\|_\infty\geq\|p_\vep\|_\infty^{t/\vep}
 \geq[(1-\vep)\vee(\lamb\vep\|D\|_\infty)]^{t/\vep},
\end{align}
which can be bounded from below by an $\vep$-independent multiple
of $L^{-d}(1+t)^{-d/2}$.  This completes the proof of
Proposition~\ref{thm-disc}.
\qed

\bigskip

Finally, we derive the expressions \refeq{lambdacdef}--\refeq{Avdef} for
$\lambc$, $A$ and $v$.  Recall \refeq{indchoices}.  The expressions for
$\lambc$ and $v$ immediately follow from \refeq{lambn-def}, \refeq{vn-def}
and the fact that $\lambc=\lamb_\infty$.  To derive the expression for
$A=\lim_{n\to\infty}f_n(0)$, we follow the same strategy as in
\cite[p.424]{hs02}.  Let $F_{\sN}=\vep\sum_{n=0}^Nf_n(0)$, which can be
approximated by $AN\vep$ as $N\to\infty$.  Summing the recursion equation
\refeq{f-renew} with $k=0$, multiplied by $\vep$, over $n=0,\dots,N-1$,
and using $f_0(0)=1$, $e_1(0)=0$, $g_1(0)=1+(\lambc-1)\vep$ and the
expression for $\lambc$, we have
\begin{align}
F_{\sN}&=g_1(0)\,F_{\sss N-1}+\sum_{n=2}^Ng_n(0)\,F_{\sss N-n}+\vep
 \sum_{n=2}^Ne_n(0)+\vep\nn\\
&=F_{\sss N-1}-\sum_{n=2}^\infty g_n(0)\,F_{\sss N-1}+\sum_{n=2}^N
 g_n(0)\,F_{\sss N-n}+\vep\sum_{n=2}^Ne_n(0)+\vep.
\end{align}
Taking the limit $N\to\infty$ of $f_{\sN}(0)=\frac{F_N-F_{N-1}}\vep$,
with the help of the bound \refeq{pimom}, we obtain
\begin{align}
A=-A\sum_{n=2}^\infty(n-1)g_n(0)+\sum_{n=2}^\infty e_n(0)+1.
\end{align}
which, with the help of \refeq{indchoices}, gives the expression
for $A$ in \refeq{Avdef}.  This completes the derivation of
\refeq{lambdacdef}--\refeq{Avdef}.
\qed

\subsection{Discussion on changes below and at four dimensions}
\label{ss:lowdim}
In dimensions $d\leq4$, the induction analysis in
Sections~\ref{ss:indhyp}--\ref{ss:discr-thm} no longer works as long as
the infection range is fixed, and we need to incorporate the factor
$L_{\sT}=L_1T^b$ into the induction hypotheses.

Recall $\alpha=bd-\frac{4-d}2>0$, and let $\omega\in(\delta,1\wedge\alpha)$
and $\hat\beta_{\sT}=\beta_1T^{-\mu}$ with $\mu\in(0,\alpha-\omega)$.  We
again define $\lamb_n=\lamb_n(T)$ and $v_n=v_n(\lamb)$ by \refeq{lambn-def}
and \refeq{vn-def}, respectively, where we emphasize the dependence on $T$
of $\lamb_n$.  However, we replace \refeq{drgd},
\refeq{In-def}, \refeq{H12} and \refeq{H3}, respectively, by
\begin{gather}
-(2+\rho)<0<\tfrac{d}2-(2+\rho)<\gamma<\gamma+\delta<\omega\wedge\Delta,
 \lbeq{drgd-lower}\\
I_n=\lamb_n+\frac{K_1\hat\beta_{\sT}}{(1+n\vep)^{1+\omega}}[-1,1],
 \lbeq{In-def-lower}\\
|\lamb_m-\lamb_{m-1}|\leq\frac{\vep K_1\hat\beta_{\sT}}{(1+m\vep)^{2+
 \omega}},\qquad\qquad
|v_m-v_{m-1}|\leq\frac{\vep K_2\hat\beta_{\sT}}{(1+m\vep)^{1+\omega}},
 \lbeq{H12-lower}\\
|r_m(0)|\leq\frac{K_3\hat\beta_{\sT}}{(1+m\vep)^{1+\omega}},\qquad\qquad
|r_m(k)-r_m(0)|\leq\frac{K_3\hat\beta_{\sT}}{(1+m\vep)^\delta}\,a(k).
 \lbeq{H3-lower}
\end{gather}
The induction hypotheses are that (H1)--(H4) hold for all $\lamb\in I_n$
and $m=1,\dots,n$, where we assume that $n\vep\leq T\log T$.  It suffices
to prove the main statement for sufficiently small $\beta_1>0$, i.e., for
sufficiently large initial infection ranges $L_1$.

We now discuss the induction hypotheses.  One of the key ingredients in
the induction is the fact that the intervals $I_n$ are decreasing in
$n\leq\frac{T}\vep\log T$.  This implies that we can use the bounds
following from (H1)--(H4) in the advancement of the induction hypotheses.
One would expect that one could choose
$I_n=\lamb_n+K_1\beta_{\sT}(1+n\vep)^{-(d-2)/2}[-1,1]$, i.e., by simply
replacing $\beta$ in \refeq{In-def} by $\beta_{\sT}$.  However, to obtain
a decreasing sequence of $I_n$, it is required for the power exponent
$(d-2)/2$ in the width of $I_n$ to be greater than 1, and it is not the
case when $d\leq4$ (cf., the proof of Lemma~\ref{lem:I} in
Appendix~\ref{s:adv-ind}).  To satisfy this requirement,
we transfer some power exponent of $\beta_{\sT}$ as
\begin{align}
\frac{\beta_{\sT}}{(1+n\vep)^{(d-2)/2}}=\frac{\beta_1T^{-db}(1+n\vep)^{(4-d)
 /2+\omega}}{(1+n\vep)^{1+\omega}}\leq\frac{c\beta_1T^{-\mu}}{(1+n\vep)^{1
 +\omega}},
\end{align}
for $T\geq1$, where we use $n\vep\leq T\log T$ and
$-bd+\frac{4-d}2+\omega=-(\alpha-\omega)<-\mu$.  This is the motivation
of the changes in \refeq{In-def-lower}--\refeq{H3-lower}.

By the above changes, \refeq{lamb-bd}--\refeq{vm-bd} are modified by
replacing $\beta$ with $\hat\beta_{\sT}$.  We have that $\lamb$ and
$v_m(\lamb)$ are both $1+O(\hat\beta_{\sT})$ for $\lamb\in I_m$ and
$m=1,\dots,n$ with $n\vep\leq T\log T$.  Similarly, we replace $\beta$
in Lemmas~\ref{lem:fabs}--\ref{lem:fdiff} and Lemma~\ref{lem:fD2} by
$\hat\beta_{\sT}$ and $\beta_{\sT}$, respectively, although the proofs
of these lemmas remain unchanged.  However, the proof of the main
result does change, due to the fact that the constants $A$ and $v$ for
$d>4$ are replaced by 1, and the fact that there is no unique limit of
$\bigcap_{m=1}^nI_m$.

\begin{proof}[Proof of Proposition~\ref{thm-disc2}]
Let $n\leq\frac{T}\vep\log T$ and $\lamb\in I_n$.  In particular, the
following results hold at $\lamb=\lamb_{\sT}$, which is defined as
\begin{align}\lbeq{lambT-def}
\lamb_{\sT}=\lamb_{\sss\frac{T}\vep\log T}(T)=1-\frac1\vep\sum_{l
 =2}^{\frac{T}\vep\log T}g_l(0;\lamb_{\sss\frac{T}\vep\log T-1}(T)).
\end{align}
By \refeq{H3-lower}, we can bound $|f_n(0)-1|$ by
\begin{align}\lbeq{fn-Alowdim}
\bigg|\prod_{m=1}^n[1+\vep r_m(0)]-1\bigg|\leq\vep\sum_{m=1}^n|r_m(0)|
 \prod_{l=m+1}^n[1+\vep|r_l(0)|]\leq\vep\sum_{m=1}^n\frac{K_3\hat
 \beta_{\sT}\,e^{\vep\sum_{l=m+1}^n|r_l(0)|}}{(1+m\vep)^{1+\omega}}
 \leq cK_3\hat\beta_{\sT},
\end{align}
which proves that the asymptotic expected number of infected individuals
is 1.  Also, using \refeq{vm-bd} and \refeq{H12-lower}, we have
$v_n=1+O(\hat\beta_{\sT})$, which means that the asymptotic diffusion
constant is 1.

Let $n=Tt/\vep$ with $t\leq\log T$ and
$\tilde k=\frac{k}{\sqrt{\sigma_{\sT}^2Tt}}\in\cA_n$.  By \refeq{Deltadef},
$a(\tilde k)=\frac{|k|^2}{2dTt}+O(|k|^{2+2\Delta}(Tt)^{-1-\Delta})$.
Therefore,
\begin{align}
f_n(\tilde k)&=[1+O(\hat\beta_{\sT})]\bigg[1-[1+O(\hat\beta_{\sT})]
 \frac{Tt\,a(\tilde k)}n\bigg]^n\nn\\
&=e^{-\frac{|k|^2}{2d}}\big[1+O(\hat\beta_{\sT})+O(|k|^{2+2\Delta}
 (Tt)^{-\Delta})+O(\vep|k|^4(Tt)^{-1})\big].
\end{align}
By \refeq{error2}--\refeq{error3} for small $\tilde k$, the last two
error terms can be replaced by $O(|k|^2(1+Tt)^{-\delta})$.  This proves \refeq{tauasy-lowdimdis}.

Using \refeq{fmdiff-bd} and \refeq{sldiff-bd} as well as
$v_n=1+O(\hat\beta_{\sT})$, we obtain \refeq{taugyr-lowdimdis}.
The proof of \refeq{tausup-lowdimdis} does not depend on $d$, and is
the same as in Section~\ref{ss:discr-thm}.  This completes the proof
of Proposition~\ref{thm-disc2}.
\end{proof}


\newcommand{\sites}{\Ccal}
\newcommand{\Asites}{\Acal}
\newcommand{\Sone}{{\bf S}_1}
\newcommand{\Stwo}{{\bf S}_2}

\section{Continuum limit}\label{s:continuum}
In this section, we compute the limit of the lace expansion coefficients
as $\vep\daw0$, and prove Proposition~\ref{prop-disc}.

%

We prove below convergence of $\frac1{\vep^2}\pi_{t;\vep}^\lamb(x)$ for
$t/\vep\in[2,\infty)\cap\Zp$ with a {\it fixed} $\lamb\leq\lambc$, and then
extend this to $\frac1{\vep^2}\partial_\lamb\pi_{t;\vep}^\lamb(x)$.  The
proof of the continuity in $\lamb$ of $\partial_\lamb\pi_t^\lamb(x)=\lim
_{\vep\daw0}\frac1{\vep^2}\partial_\lamb\pi_{t;\vep}^\lamb(x)$ is more or
less immediate from its finite containment property that is similar to the
one for the discretized contact process in Section~\ref{ss:fc}, and this
will be discussed briefly at the end of this section.  These statements
imply convergence of $\frac1{\vep^2}\pi_{t;\vep}^{\lamb^{(\vep)}}(x)$
whenever $\lamb^{\sss(\vep)}\to\lamb$ such that
$\lamb^{\sss(\vep)}\leq\lambc^{\sss(\vep)}$ for $\vep$ sufficiently
small.  Indeed, for any $\lamb_0<\lamb\leq\lambc$, we can write
\begin{align}
\bigg|\frac1{\vep^2}\pi_{t;\vep}^{\lamb^{(\vep)}}(x)-\pi_t^{\lamb_0}(x)
 \bigg|\leq\frac1{\vep^2}\big|\pi_{t;\vep}^{\lamb^{(\vep)}}(x)-\pi_{t;
 \vep}^{\lamb_0}(x)\big|+\bigg|\frac1{\vep^2}\pi_{t;\vep}^{\lamb_0}(x)
 -\pi_t^{\lamb_0}(x)\bigg|.\lbeq{triangleineq}
\end{align}
The second term in \refeq{triangleineq} converges to zero by assumption,
while we estimate the first term by
\begin{align}
\frac1{\vep^2}\big|\pi_{t;\vep}^{\lamb^{(\vep)}}(x)-\pi_{t;\vep}^{\lamb_0}
 (x)\big|\leq\int_{\lamb_0}^{\lamb^{(\vep)}}d\lamb'~\frac1{\vep^2}\big|
 \partial_{\lamb'}\pi_{t;\vep}^{\lamb'}(x)\big|,
\end{align}
where we use $\lamb_0\leq\lambc^{\sss(\vep)}$ for sufficiently small
$\vep$, which is due to the fact that $\lambc^{\sss(\vep)}$ converges to
$\lambc>\lamb_0$.  Since the integrand is uniformly bounded (even when we
sum over $x$) by $C\beta (1+t)^{-(d-2)/2}$, the limsup of the integral
when $\vep\downarrow 0$ is bounded by a multiple of $\lamb-\lamb_0$.
By taking the limit $\lamb_0\uaw\lamb$ and using the fact that
$\partial_\lamb\pi_t^\lamb(x)$ is continuous in $\lamb$, the first
expression in \refeq{piconvvep} follows.  Therefore, we are left to prove
convergence of $\frac1{\vep^2}\pi_{t;\vep}^\lamb(x)$ and
$\frac1{\vep^2}\partial_\lamb\pi_{t;\vep}^\lamb(x)$ for every
$\lamb\leq\lambc$ and the continuity in $\lamb$ of
$\partial_\lamb\pi_t^\lamb(x)$.


We prove below that, for every $N\geq0$, $\lamb\leq\lambc$ and $(x,t)$,
there is a $\pi^{\sss(N)}_t(x)$ such that
\begin{align}\lbeq{pointwise}
\lim_{\vep\daw0}\frac1{\vep^2}\pi^{\sss(N)}_{t;\vep}(x)=\pi_t^{\sss(N)}(x),
\end{align}
where we suppress the dependence on $\lamb$.  That is, we will deal with
{\it pointwise convergence}, rather than the uniform bounds in
Section~\ref{s:pibds} which are valid for all $(x,t)$ and $\vep\leq1$,
and hence all terms which are $o(1)$ as $\vep\daw0$ will be estimated
away.  By this pointwise convergence, together with the uniform bounds
in Section~\ref{s:pibds} and the dominated convergence theorem, we have
\begin{align}
\lim_{\vep\daw0}\frac1{\vep^2}\pi_{t;\vep}^\lamb(x)=\lim_{\vep\daw0}
 \sum_{N=0}^\infty(-1)^N\frac1{\vep^2}\pi^{\sss(N)}_{t;\vep}(x)=
 \sum_{N=0}^\infty(-1)^N\pi^{\sss(N)}_t(x)=\pi^\lamb_t(x).
\end{align}
This completes the proof of the pointwise convergence of
$\frac1{\vep^2}\pi_{t;\vep}^\lamb(x)$.  The proof of the convergence of
$\frac1{\vep^2}\partial_\lamb\pi_{t;\vep}^\lamb(x)$ is similar, and we
will only discuss the necessary changes.

The proof of \refeq{pointwise} is divided into several steps.

\paragraph{Statement of the induction hypothesis.}
Given a site set $\cC\subset\Zd$ (which may be an empty set), we define
\begin{align}\lbeq{piCdef}
\pi^{\sss(N)}_{t;\vep}(x;\sites)=\sum_{\vec b_N}\mP_\vep^\lamb\big(\tilde
 E^{\sss(N)}_{\vec b_N}(x,t)\cap\big\{\bC_t(\tb_{\sN})\setminus\{x\}=\cC
 \big\}\big),
\end{align}
where we recall $\tilde E^{\sss(0)}_{\vec b_0}(x,t)=\{(o,0)\db(x,t)\}$,
$\tb_0=(o,0)$ and the notation \refeq{bC-gen} for $\bC_t(\tb_{\sN})$.
We will use induction in $N$ to prove that, for every $t>0$, there is a
$\pi^{\sss(N)}_t(x;\cC)$ such that
\begin{align}
\lim_{\vep\daw0}\frac1{\vep^2}\pi^{\sss(N)}_{t;\vep}(x;\cC)
 =\pi^{\sss(N)}_t(x;\cC).\lbeq{cont-indhyp}
\end{align}
The claim for $\pi_{t;\vep}^{\sss(N)}(x)$ in \refeq{pointwise} then follows
by summing over $\cC\subset\Zd$, together with the fact that the main
contribution comes from $\cC\subset\Box_{\sR}$ by the finite containment
property in Section~\ref{ss:fc}.

\paragraph{Initialization of the induction.}
First, we investigate $N=0$.  For $\Sone,\Stwo,\bA\subset\Zd\times\vep\Zp$,
we denote
\begin{align}
\{\Sone\conn\Stwo\}=\bigcup_{\substack{\svec_1\in\Sone\\ \svec_2\in\Stwo}}
 \{\svec_1\conn\svec_2\},&&
\{\Sone\db\Stwo\}=\!\!\bigcup_{\substack{\svec_1,\svec_1'\in\Sone\\ \svec_2,
 \svec_2'\in\Stwo}}\{\svec_1\conn\svec_2\}\circ\{\svec_1'\conn\svec_2'\},
\end{align}
and define
\begin{align}\lbeq{CtA}
\bC_t(\bA)=\{x\in\Zd:\bA\conn(x,t)\}=\bigcup_{\avec\in\bA}
 \bC_t(\avec),&& \bC(\bA)=\bigcup_{t\geq0}\bC_t(\bA).
\end{align}
Using the Markov property at time $\vep$, we arrive at
\begin{align}
\pi^{\sss(0)}_{t;\vep}(x;\cC)&=\sum_{\Asites\subset\Zd:|\cA|\geq2}\bigg[
 \prod_{a\in\Asites}p_\vep(a)\bigg]\bigg[\prod_{a\notin\Asites}[1-p_{\vep}
 (a)]\bigg]\nn\\
&\qquad\qquad\times\mP_\vep^\lamb\left(\!\begin{array}{c}
 \exists\,a,a'(\ne a)\in\Asites:\{(a,\vep)\conn(x,t)\}\circ\{(a',\vep)
  \conn(x,t)\}\\
 \bC_t(\Asites\times\{\vep\})\setminus\{x\}=\cC
 \end{array}\!\right),\lbeq{cont-prepi01st}
\end{align}
Since every $p_\vep(a)$ for $a\neq o$ gives rise to a factor of $\vep$,
we immediately see that the main contribution comes from $\Asites=\{o,y\}$
for some $y\ne o$.  Therefore, we obtain
\begin{align}
\frac1\vep\pi^{\sss(0)}_{t;\vep}(x;\cC)&=\sum_{y\in\Zd}\lamb D(y)~
 \mP_\vep^\lamb\left(\!\begin{array}{c}
 \{(o,\vep)\conn(x,t)\}\circ\{(y,\vep)\conn(x,t)\}\\
 \bC_t(\{(o,\vep),(y,\vep)\})\setminus\{x\}=\cC
 \end{array}\!\right)+o(1)\nn\\
&=\sum_{y\in\Zd}\lamb D(y)~\mP_\vep^\lamb\left(\!\begin{array}{c}
 \{(o,0)\conn(x,t)\}\circ\{(y,0)\conn(x,t)\}\\
 \bC_t(\{(o,0),(y,0)\})\setminus\{x\}=\cC
 \end{array}\!\right)+o(1)\nn\\
&=\sum_{y\in\Zd}\lamb D(y)~\mP_\vep^\lamb\left(\!\begin{array}{c}
 \{(o,0),(y,0)\}\db(x,t)\\ \bC_t(\{(o,0),(y,0)\})\setminus\{x\}=\cC
 \end{array}\!\right)+o(1),\lbeq{cont-pi01st}
\end{align}
where the second equality is due to the fact that $((o,0),(o,\vep))$ or
$((y,0),(y,\vep))$ is vacant (with probability $(2-\vep)\vep$) in the
symmetric difference between the events on both sides of the equality,
and the third equality is due to the fact that the double connection
from $(o,\vep)$ or $(y,\vep)$ gives rise to an extra factor of $\vep$.

We repeat the same observation around $(x,t)$ and obtain that, for every
$y\in\Zd$ and $\cC\subset\Zd$,
\begin{align}
&\frac1\vep\mP_\vep^\lamb\left(\!\begin{array}{c}
 \{(o,0),(y,0)\}\db(x,t)\\ \bC_t(\{(o,0),(y,0)\})\setminus\{x\}=\cC
 \end{array}\!\right)\nn\\
&\qquad=\sum_{z\in\Zd}\lamb D(x-z)~\mP_\vep^\lamb\left(\!\begin{array}{c}
 \{(o,0),(y,0)\}\db\{(x,t),(z,t)\}\\
 \bC_t(\{(o,0),(y,0)\})\setminus\{x\}=\cC
 \end{array}\!\right)+o(1).\lbeq{cont-pi02nd}
\end{align}
Substituting \refeq{cont-pi02nd} into \refeq{cont-pi01st} and using the
weak convergence of $\mP_\vep^\lamb$ towards $\mP^\lamb$ as formulated in
\cite[Proposition 2.7]{bg91}, we obtain
\begin{align}\lbeq{pi0lim}
\lim_{\vep\daw0}\frac1{\vep^2}\pi^{\sss(0)}_{t;\vep}(x;\cC)=\sum_{y,z\in\Zd}
 \lamb^2 D(y)\,D(x-z)~\mP^\lamb\left(\!\begin{array}{c}
 \{(o,0),(y,0)\}\db\{(x,t),(z,t)\}\\
 \bC_t(\{(o,0),(y,0)\})\setminus\{x\}=\cC
 \end{array}\!\right).
\end{align}
Here and in the rest of this section, we use ``$\conn$'' and ``$\db$'' to
denote connections in $\Zd\times\mR_+$, via the graphical representation
in Section~\ref{ss:discretize}.  This completes the proof of
\refeq{cont-indhyp} for $N=0$, with $\pi^{\sss(0)}_t(x;\cC)$ for $t>0$
defined to be the right-hand side of \refeq{pi0lim}.

\paragraph{Preliminaries for the advancement.}
To advance the induction hypothesis, we first note that, by using the
finite containment property of Section~\ref{ss:fc} and the Markov property
at the time component of $\bb_{\sN}$, we have
\begin{align}\lbeq{cont-dec}
&\frac1{\vep^2}\pi^{\sss(N)}_{t;\vep}(x;\cC)=\frac1{\vep^2}\sum_{\vec b_N}
 \mP_\vep^\lamb\Big(\tilde E_{\vec b_{N-1}}^{\sss(N-1)}(\bb_{\sN})\cap E(
 b_{\sN},(x,t);\tilde\bC^{b_N}(\tb_{\sss N-1}))\cap\big\{\bC_t(\tb_{\sN})
 \setminus\{x\}=\cC\big\}\Big)\nn\\
&\qquad=\frac1{\vep^2}\ddsum_{s=\vep}^t\sum_{y\in\Zd}\;\sum_{\vec b_N:\bb_N
 =(y,s-\vep)}\;\sum_{\cC'\subset\Box_R:\cC'\ni y}\mP_\vep^\lamb\big(\tilde
 E_{\vec b_{N-1}}^{\sss(N-1)}(y,s-\vep)\cap\{\bC_{s-\vep}(\tb_{\sss N-1})=
 \cC'\}\big)\nn\\
&\qquad\qquad\qquad\times\mP_\vep^\lamb\Big(E\big(b_{\sN},(x,t);\tilde
 \bC^{b_N}(\cC'\times\{s-\vep\})\big)\cap\big\{\bC_t(\tb_{\sN})\setminus
 \{x\}=\cC\big\}\Big)+o(1),
\end{align}
where, similarly to \refeq{CtA}, we write $\tilde\bC^{b_N}(\cC'\times
\{s-\vep\})=\bigcup_{v\in\cC'}\tilde\bC^{b_N}(v,s-\vep)$, and $o(1)$
is independent of $\vep$ and decays to zero as $R\to\infty$.  We now
investigate the second probability in \refeq{cont-dec} when $\cC'=\{y\}$
and when $\cC'\supsetneq\{y\}$, separately.

When $\cC'=\{y\}$, we recall the definitions \refeq{E'def}--\refeq{Edef}
and use the Markov property at time $s$, similarly to the discussion
around \refeq{cont-prepi01st}--\refeq{cont-pi01st}, to obtain that
\begin{align}\lbeq{|cC'|=1}
&\sum_{b_N:\bb_N=(y,s-\vep)}\frac1\vep\mP_\vep^\lamb\Big(E(b_{\sN},(x,t);
 \tilde\bC^{b_N}(y,s-\vep))\cap\big\{\bC_t(\tb_{\sN})\setminus\{x\}=\cC
 \big\}\Big)\nn\\
&\qquad=\sum_{y'\in\Zd\setminus\{y\}}\lamb D(y'-y)\,\bigg[\mP_\vep^\lamb
 \Big(E'\big((y,s),(x,t);\bC(y',s)\big)\cap\big\{\bC_t(y,s)\setminus\{x\}
 =\cC\big\}\Big)\nn\\
&\qquad\qquad\qquad\qquad+\mP_\vep^\lamb\Big(E'\big((y',s),(x,t);\bC(y,s)
 \big)\cap\big\{\bC_t(y',s)\setminus\{x\}=\cC\big\}\Big)\bigg]+o(1),
\end{align}
where $o(1)$ decays to zero as $\vep\daw0$, and the first probability in
the brackets is the contribution from the case in which $b_{\sN}$ is the
temporal bond $((y,s-\vep),(y,s))$, while the second probability is the
contribution from the case in which $b_{\sN}$ is the spatial bond
$((y,s-\vep),(y',s))$.  In \refeq{|cC'|=1}, we also use the fact that, with
probability $1-o(1)$, $\tilde\bC^{b_N}(y,s-\vep)\cap(\Zd\times[s,\infty))$
equals $\bC(y',s)$ when $b_{\sN}=((y,s-\vep),(y,s))$, and equals
$\bC(y,s)$ when $b_{\sN}=((y,s-\vep),(y',s))$.

When $\cC'\supsetneq\{y\}$, we again use the Markov property at time $s$,
and then we use the fact that, with probability $1-o(1)$, every temporal
bond growing from each site in $\cC'\times\{s-\vep\}$ is occupied, and
all the spatial bonds growing from the sites in $\cC'\times\{s-\vep\}$
are vacant.  Therefore, with probability $1-o(1)$, the subset of
$\tilde\bC^{b_N}(\cC'\times\{s-\vep\})$ after time $s$ equals
$\bC\big((\cC'\setminus\{y\})\times\{s\}\big)$, and we have
\begin{align}\lbeq{|cC'|>1}
&\sum_{b_N:\bb_N=(y,s-\vep)}\mP_\vep^\lamb\Big(E\big(b_{\sN},(x,t);\tilde
 \bC^{b_N}(\cC'\times\{s-\vep\})\big)\cap\big\{\bC_t(\tb_{\sN})\setminus
 \{x\}=\cC\big\}\Big)\nn\\
&\qquad=\mP_\vep^\lamb\Big(E'\big((y,s),(x,t);\bC\big((\cC'\setminus\{y\})
 \times\{s\}\big)\big)\cap\big\{\bC_t(y,s)\setminus\{x\}=\cC\big\}\Big)+o(1).
\end{align}

To deal with the event $E'((y,s),(x,t);\cA\times\{s\})$ for
$\cA\subset\Zd\setminus\{y\}$ in \refeq{|cC'|=1}--\refeq{|cC'|>1}, we
introduce some notation.  We define the set of sites that are connected
from $(y,s)$ via a path which does not go through $\vvec$ by
\begin{align}
\tilde\bC^{\vvec}(y,s)=\bigcap_{b=(\,\cdot\,,\vvec)}\tilde\bC^b(y,s).
\end{align}
We also define
\begin{align}
\cE_{s,t}(y,x;\cA)&=\bigcup_{\vvec}\Big\{\big\{\vvec\notin\bC(\cA\times
 \{s\})\big\}\cap\big\{(y,s)\conn\vvec\db(x,t)\in\bC(\cA\times\{s\})
 \setminus\tilde\bC^{\vvec}(y,s)\big\}\Big\},\lbeq{cE-def}\\
\cR_{s,t}(y,x;\cA)&=\bigcup_{\vvec}\Big\{\big\{\cA\times\{s\}\conn\vvec
 \big\}\circ\big\{(y,s)\conn\vvec\db(x,t)\big\}\Big\}.\lbeq{cR-def}
\end{align}
By this notation, it is not hard to see that
$E'((y,s),(x,t);\cA\times\{s\})$ is rewritten as
\begin{align}\lbeq{Erewr}
E'((y,s),(x,t);\cA\times\{s\})=\cE_{s,t}(y,x;\cA)\Dcup{}{}\cR_{s,t}
 (y,x;\cA).
\end{align}
The contribution from $\cR_{s,t}(y,x;\cA)$ has an extra factor of $\vep$,
due to the fact that there are at least {\it two} spatial bonds at $\vvec$
(one before and one after $\vvec$), which leads to an error term as
$\vep\daw0$.  Therefore, we only need to focus on the contribution from
$\cE_{s,t}(y,x;\cA)$, i.e.,
\begin{align}
\mP_\vep^\lamb\big(\cE_{s,t}(y,x;\cA)\cap\big\{\bC_t(y,s)\setminus
 \{x\}=\cC\big\}\big).
\end{align}
Generalizing the definition \refeq{cE-def} from a single end $(x,t)$
to a pair $\{(x,t),(z,t)\}$ with $x\ne z$ as
\begin{align}\lbeq{cEz-def}
&\cE_{s,t}(y,\{x,z\};\cA)\\
&\quad=\bigcup_{\vvec}\Big\{\big\{\vvec\notin\bC(\cA\times\{s\})\big\}\cap
 \big\{(y,s)\conn\vvec\db\{(x,t),(z,t)\}\subset\bC(\cA\times\{s\})\setminus
 \tilde\bC^{\vvec}(y,s)\big\}\Big\},\nn
\end{align}
and following the argument in \refeq{cont-pi02nd} (see also the discussion
around \refeq{cont-prepi01st}--\refeq{cont-pi01st}), we obtain
\begin{align}\lbeq{lim1vep2}
&\frac1\vep\mP_\vep^\lamb\big(\cE_{s,t}(y,x;\cA)\cap\big\{\bC_t(y,s)
 \setminus\{x\}=\cC\big\}\big)\nn\\
&\qquad=\sum_{z\in\Zd}\lamb D(x-z)~\mP_\vep^\lamb\big(\cE_{s,t}(y,\{x,z\};
 \cA)\cap\big\{\bC_t(y,s)\setminus\{x\}=\cC\big\}\big)+o(1).
\end{align}

\paragraph{Advancement of the induction hypothesis.}
Now we advance the induction hypothesis in $N\geq1$ by using
\refeq{cont-dec}--\refeq{|cC'|>1} and \refeq{lim1vep2}.

First, we consider the contribution from $\cC'=\{y\}$ in \refeq{cont-dec},
which equals
\begin{align}\lbeq{advind=1}
&\ddsum_{s=\vep}^t\sum_{y\in\Box_R}\,\sum_{b_N=(\,\cdot\,,(y,s-\vep))}
 \pi_{s-\vep;\vep}^{\sss(N-1)}(y;\varnothing)~\frac1{\vep^2}\mP_\vep^\lamb
 \big(E(b_{\sN},(x,t);\tilde\bC^{b_N}(y,s-\vep))\cap\big\{\bC_t(\tb_{\sN})
 \setminus\{x\}=\cC\big\}\big)\nn\\
&\quad=\delta_{N,1}\sum_{b=(\ovec,\,\cdot\,)}\frac1{\vep^2}\mP_\vep^\lamb
 \big(E(b,(x,t);\tilde\bC^b(\ovec))\cap\big\{\bC_t(\tb)\setminus\{x\}=\cC
 \big\}\big)\\
&\qquad+\vep^2\ddsum_{s=2\vep}^t\sum_{y\in\Box_R}\frac1{\vep^2}\pi_{s-\vep;
 \vep}^{\sss(N-1)}(y;\varnothing)\shift\sum_{b=((y,s-\vep),\,\cdot\,)}\frac1
 {\vep^2}\mP_\vep^\lamb\big(E(b,(x,t);\tilde\bC^b(y,s-\vep))\cap\big\{
 \bC_t(\tb)\setminus\{x\}=\cC\big\}\big),\nn
\end{align}
where we use
$\pi_{0;\vep}^{\sss(N-1)}(y;\varnothing)=\delta_{N,1}\,\delta_{o,y}$ to
obtain the first term in the right-hand side.  We note that, by using the
induction hypothesis, as well as \refeq{|cC'|=1} and \refeq{lim1vep2},
the second term is $O(\vep)=o(1)$.  Therefore, the first term is the main
contribution.  By using \refeq{|cC'|=1} and \refeq{lim1vep2} again, as well
as the weak convergence of $\mP_\vep^\lamb$, the limit in $R\uaw\infty$ of
the continuum limit of \refeq{advind=1} equals
\begin{gather}
\delta_{N,1}\sum_{y,z\in\Zd}\lamb^2D(y)\,D(x-z)\,\Big[\mP^\lamb\big(
 \cE_{0,t}(o,\{x,z\};\{y\})\cap\big\{\bC_t(\ovec)\setminus\{x\}=\cC
 \big\}\big)\nn\\
+\mP^\lamb\big(\cE_{0,t}(y,\{x,z\};\{o\})\cap\big\{\bC_t(y,0)\setminus
 \{x\}=\cC\big\}\big)\Big].\lbeq{contlim=1}
\end{gather}

Next, we consider the contribution from $\cC'\supsetneq\{y\}$ in
\refeq{cont-dec}, which equals
\begin{align}\lbeq{advind>1}
&\vep\ddsum_{s=2\vep}^t\sum_{y\in\Box_R}\,\sum_{\cA\subset\Box_R\setminus
 \{y\}:\cA\ne\varnothing}\frac1{\vep^2}\pi_{s-\vep;\vep}^{\sss(N-1)}(y;\cA)
 \nn\\
&\qquad\times\sum_{z\in\Zd}\lamb D(x-z)~\mP_\vep^\lamb\big(\cE_{s,t}(y,\{x,
 z\};\cA)\cap\big\{\bC_t(y,s)\setminus\{x\}=\cC\big\}\big)+o(1),
\end{align}
where we use \refeq{|cC'|>1} and \refeq{lim1vep2}, as well as
$\pi_{0;\vep}^{\sss(N-1)}(y;\cA)=0$ for $\cA\ne\varnothing$ (so that
the sum over $s$ starts from $s=2\vep$).  By the dominated convergence
theorem, as well as the induction hypothesis and the weak convergence
of $\mP_\vep^\lamb$, the limit in $R\uaw\infty$ of the continuum limit
of \refeq{advind>1} equals
\begin{align}\lbeq{contlim>1}
\int_0^tds\sum_{y,z\in\Zd}\lamb D(x-z)\!\sum_{\substack{\cA\subset\Zd
 \setminus\{y\}\\ \cA\ne\varnothing}}\pi_s^{\sss(N-1)}(y;\cA)~\mP^\lamb\big(
 \cE_{s,t}(y,\{x,z\};\cA)\cap\big\{\bC_t(y,s)\setminus\{x\}=\cC\big\}\big).
\end{align}

Therefore, the limit $\pi_t^{\sss(N)}(x;\cC)$ for $t>0$ exists and equals
the sum of \refeq{contlim=1} and \refeq{contlim>1}.  This advances the
induction hypothesis.

\paragraph{Bounds on $\pi_t^\lamb$ in \refeq{boundspiCP} and convergence
of $A^{\smallsup{\vep}}$ and $v^{\smallsup{\vep}}$.}
The bound on $\sum_{x\in \Z^d} |x|^q \pi_t^{\lamb}(x)$ follow immediately
from the pointwise convergence of $\frac1{\vep^2}\pi_{t;\vep}^\lamb(x)$,
together with the uniform bounds in Proposition~\ref{lem-Pibd.smp} and
dominated convergence for the sum over $x$.

To prove convergence of $A^{\sss(\vep)}$ and $v^{\sss(\vep)}$, we first
note that by \cite[Section~3.1]{s01}, $\lambc^{\sss(\vep)}\to\lambc$.
Convergence of $A^{\sss(\vep)}$ and $v^{\sss(\vep)}$ follows by dominated
convergence, together with the identification of $A^{\sss(\vep)}$ and
$v^{\sss(\vep)}$ in \refeq{Avdef}.  Thus, we obtain that
\begin{align}
A=\bigg[1+\int_0^\infty\!dt~ t\,\hat\pi_t^{\lambc}(0)\bigg]^{-1},&&
v=A\bigg[\lambc-\frac1{\sigma^2}\int_0^\infty\!dt~\hat\nabla^2
 \pi_t^{\lambc}(0)\bigg].
\end{align}

\paragraph{Convergence of $\frac1{\vep^2}\partial_\lamb\pi_{t;\vep}^\lamb
(x)$ and the bound on $\partial_\lamb\pi_t^\lamb$ in \refeq{boundspiCP}.}
The only difference between $\frac1{\vep^2}\partial_\lamb\pi_{t;\vep}^\lamb
(x)$ and $\frac1{\vep^2}\pi_{t;\vep}^\lamb(x)$ is the occurrence of the sum
over spatial bonds $b$ and the indicator of the event
$b\in\{b_n\}\cup\piv[\tb_n,\bb_{n+1}]$.  Clearly, the main term in the
above comes from $b\in\piv[\tb_n,\bb_{n+1}]$.  The extra inclusion of this
event gives rise to an extra integral over the time variable $r'$ and an
indicator that the arrow $((w,r'),(w',r'))$ is pivotal for the connection
from $\tb_n$ to $\bb_{n+1}$ (see \cite[p.61]{Ligg99} for the definition of
a pivotal arrow).  Apart from this minor modification, the proof remains
unchanged.  The bound on $\sum_{x\in\Zd}|\partial_\lamb\pi_t^\lamb(x)|$ in
\refeq{boundspiCP} follows immediately from the pointwise convergence of
$\frac1{\vep^2}\partial_\lamb\pi_{t;\vep}^\lamb(x)$, together with the
uniform bounds in Proposition~\ref{lem-Pibd.smp} and dominated convergence
for the sum over $x\in\Zd$.

\paragraph{Continuity in $\lamb$ of $\partial_\lamb\pi_t^\lamb(x)$.}
Following the same strategy as above, we may obtain an explicit
expression for $\partial_\lamb\pi_t^\lamb(x)$, similar to the expression
obtained for $\pi_t^\lamb(x)$ from \refeq{pi0lim}, \refeq{contlim=1} and
\refeq{contlim>1}.  Let $\partial_\lamb\pi_t^\lamb(x\,|\,R)$ be equal to
$\partial_\lamb\pi_t^\lamb(x)$ with the extra condition
$\{\bC_{[0,t]}\subset\Box_{\sR}\}$ being imposed, as in \refeq{pi-fc}
for the discretized contact process.  Note that, as explained above,
$\partial_\lamb\pi_t^\lamb(x)=\partial_\lamb\pi_t^\lamb(x\,|\,R)+o(1)$,
where $o(1)$ decays to zero as $R\to\infty$, and that
$\partial_\lamb\pi_t^\lamb(x\,|\,R)$ is continuous in $\lamb$
since it depends only on events in the finite space-time box
$\Box_{\sR}\times[0,t]$.  Therefore, $\partial_\lamb\pi_t^\lamb(x)$ is
also continuous in $\lamb$.  This completes the proof.
\qed


\appendix
\section{Advancement of the induction hypotheses}\label{s:adv-ind}
In this appendix, we prove Lemmas~\ref{lem:I}--\ref{lem:fD2} and we
advance the induction hypotheses.  We discuss the case of $d>4$ in
Appendix~\ref{ss:adv-high}, which is quite similar to the argument in
\cite{hs02}.  The main difference is due to the required uniformity in
$\vep$.  We will explain in detail how to use the factors of $\vep$
contained in the induction hypotheses and in the bounds
\refeq{pimom}--\refeq{pider}, in order to obtain this uniformity.  The
argument for $d\leq4$ is almost identical, except for modifications due
to the factors $\beta_{\sT}$ and $\hat\beta_{\sT}$ in
\refeq{PN-bd-lower}--\refeq{tPNn-bd-lower} and
\refeq{In-def-lower}--\refeq{H3-lower}.
We discuss the necessary changes for $d\leq4$ in Appendix~\ref{ss:adv-low}.

\subsection{Advancement above four dimensions}\label{ss:adv-high}
\subsubsection{Proofs of Lemmas~\ref{lem:I}--\ref{lem:fD2}}
Recall the induction hypotheses (H1)--(H4) and the definitions of
$\lamb_n$, $I_n$ and $v_n$ in Section~\ref{sss:indhyp}.  We now
prove Lemmas~\ref{lem:I}--\ref{lem:fD2} using the induction hypotheses.
\begin{proof}[Proof of Lemma~\ref{lem:I}]
We prove $\lamb\in I_{m-1}$ assuming $\lamb\in I_m$.  By \refeq{In-def} and
\refeq{H12},
\begin{align}
|\lamb-\lamb_{m-1}|\leq|\lamb-\lamb_m|+|\lamb_m-\lamb_{m-1}|\leq K_1\beta
 \frac{1+(m+1)\vep}{(1+m\vep)^{d/2}}\leq\frac{K_1\beta}{[1+(m-1)\vep]^{(d-2)/2}},
\end{align}
where the last inequality is due to the fact that
$f(\vep)=(c+\vep)(c-\vep)^a$ is decreasing in $\vep\geq0$ if $c>0$ and
$a\geq1$, so that $f(\vep)\leq f(0)=c^{1+a}$ (in the above inequality, $c=1+m\vep$ and $a=\frac{d-2}2$).
This completes the proof of $I_m\subset I_{m-1}$.
\end{proof}

\begin{proof}[Proof of Lemma~\ref{lem:fabs}]
By \refeq{fm-expr}--\refeq{H3} and the trivial inequality $1+x\leq e^x$,
\begin{align}\lbeq{fm0-bd}
|f_m(0)|=\bigg|\prod_{l=1}^m[1+\vep r_l(0)]\bigg|\leq e^{\vep\sum_{l=1}^m
 |r_l(0)|}\leq e^{cK_3\beta}.
\end{align}
By \refeq{fm-reexpr}--\refeq{sl-bd} and \refeq{vm-bd},
$|f_m(k)/f_m(0)|$ is bounded by
\begin{align}\lbeq{fmnz-bd}
\bigg|\prod_{l=1}^m[1-\vep v_l\,a(k)+\vep\,s_l(k)]\bigg|\leq e^{-\vep
 \sum_{l=1}^m[v_la(k)-|s_l(k)|]}\leq e^{-m\vep[1-c\,(K_1+K_2+K_3)\beta]
 \,a(k)}.
\end{align}
This completes the proof.
\end{proof}

\begin{proof}[Proof of Lemma~\ref{lem:fdiff}]
This is an immediate consequence of \refeq{fmdiff-bd}, \refeq{vm-bd},
\refeq{sldiff-bd} and \refeq{fm0-bd}.
\end{proof}

\begin{proof}[Proof of Lemma~\ref{lem:fD2}]
Recalling $\cA_m\equiv\{k:a(k)\leq\gamma\frac{\log(2+m\vep)}{1+m\vep}\}$,
we define
\begin{align*}
R_1=\{k\in\cA_m:\|k\|_\infty\leq L^{-1}\},&&
 R_2=\{k\in\cA_m:\|k\|_\infty>L^{-1}\},\\
R_3=\{k\notin\cA_m:\|k\|_\infty\leq L^{-1}\},&&
 R_4=\{k\notin\cA_m:\|k\|_\infty>L^{-1}\},
\end{align*}
where $R_2$ is empty if $m\gg1$.  Then,
\begin{align}
\|\wD^2f_m\|_1=\sum_{i=1}^4\int_{R_i}\frac{d^dk}{(2\pi)^d}~\wD(k)^2|f_m(k)|.
\end{align}

On $R_1$, we consider the cases of $m\vep<1$ and $m\vep\geq1$ separately.
If $m\vep<1$, we use Lemma~\ref{lem:fabs} and obtain
\begin{align}
\int_{R_1}\frac{d^dk}{(2\pi)^d}~\wD(k)^2|f_m(k)|\leq c\int_{R_1}\frac{d^dk}
 {(2\pi)^d}~\wD(k)^2\leq\frac{c\beta}{(1+m\vep)^{d/2}}.
\end{align}
If $m\vep\geq1$, we use the inequality $\wD^2(k)\leq1$, Lemma~\ref{lem:fabs},
and then the assumption $a(k)\asymp L^2|k|^2$ for $\|k\|_\infty\leq L^{-1}$,
and obtain
\begin{align}
\int_{R_1}\frac{d^dk}{(2\pi)^d}~\wD(k)^2|f_m(k)|\leq c\int_{R_1}\frac{d^dk}
 {(2\pi)^d}~e^{-cm\vep L^2|k|^2}\leq\frac{c\beta}{(1+m\vep)^{d/2}}.
\end{align}
Summarizing both cases, we obtain the desired bound on the contribution
from $R_1$.

On $R_2$, we use Lemma~\ref{lem:fabs} and the assumption $a(k)>\eta$ for
$\|k\|_\infty>L^{-1}$ to conclude that there exists an $r>1$ independently
of $\beta$ such that
\begin{align}
\int_{R_2}\frac{d^dk}{(2\pi)^d}~\wD(k)^2|f_m(k)|\leq c\int_{R_2}\frac{d^dk}
 {(2\pi)^d}~\wD(k)^2r^{-m\vep}\leq c\beta\,r^{-m\vep}.
\end{align}
Since $r^{-m\vep}\leq c\,(1+m\vep)^{-d/2}$, we obtain the desired bound on
the contribution from $R_2$.

On $R_3$ and $R_4$, we use (H4).  Then, the contribution from these two
regions is bounded by
\begin{align}
\frac{K_4}{(1+m\vep)^{d/2}}\sum_{i=3}^4\int_{R_i}\frac{d^dk}{(2\pi)^d}~
 \frac{\wD(k)^2}{a(k)^{2+\rho}}.
\end{align}
It thus suffices to bound the integral by $c\beta$.  On $R_3$, we use the
inequality $\wD(k)^2\leq1$ and the assumption $a(k)\asymp L^2|k|^2$ for
$\|k\|_\infty\leq L^{-1}$.  Since $d>2(2+\rho)$ (cf., \refeq{drgd}), we
obtain
\begin{align}
\int_{R_3}\frac{d^dk}{(2\pi)^d}~\frac{\wD(k)^2}{a(k)^{2+\rho}}\leq\frac
 {c}{L^{4+2\rho}}\int_{\|k\|_\infty\leq L^{-1}}\frac{d^dk}{|k|^{4+2\rho}}
 \leq c\beta.
\end{align}
On $R_4$, we use the assumption $a(k)>\eta$ for $\|k\|_\infty>L^{-1}$ and
the fact that $\int\frac{d^dk}{(2\pi)^d}~\wD(k)^2\leq\beta$, to obtain the
desired bound $c\beta$ on the integral over $R_4$.  This completes the proof.
\end{proof}

\subsubsection{Initialization and advancement of the induction hypotheses}
First we verify that the induction hypotheses hold for $n=1$.
\begin{description}
\item{\bf(H1)--(H2)}\quad
By definition, $|\lamb_1-\lamb_0|=|v_1-v_0|=0$.
\item{\bf(H3)}\quad
By \refeq{f1g1-def} and \refeq{fm-expr}, $r_1(k)\equiv\lamb-1$ and thus
$|r_1(k)-r_1(0)|\equiv0$.  Together with $\lamb\in I_1$, we obtain
$|r_1(0)|\leq K_1\beta/(1+\vep)^{(d-2)/2}$.  Therefore, (H3) holds, if
$K_3\geq K_1$.
\item{\bf(H4)}\quad
By \refeq{f1g1-def}, $|f_1(k)|\leq1+3\vep$ and $|f_1(k)-f_0(k)|\leq3\vep$
for $\beta\ll1$. Together with the trivial bound $a(k)\leq2$, (H4) is
proved to hold, if $K_4\geq(1+3\vep)2^{2+\rho}(1+\vep)^{d/2}$ and $K_5\geq3\cdot2^{1+\rho}(1+\vep)^{d/2}$.
\end{description}

Next we advance the induction hypotheses for $\lamb\in I_{n+1}$ under
the assumption that (H1)--(H4) hold for all $m\leq n$.  As mentioned
below Lemma~\ref{lem:fD2}, this assumption implies \refeq{fbd} for all
$s\leq n\vep$ if $K\gg K_4$, and thus implies \refeq{pimom}--\refeq{pider}
for all $s\leq n\vep+\vep$.  By \refeq{indchoices}, these bounds are
translated into the following bounds for all $m\leq n+1$: there is a
$\Cg<\infty$ such that
\begin{align}
\lbeq{ebd}|e_m(k)|&\leq\frac{\vep^2\Cg\beta}{(1+m\vep)^{d/2}},&
|e_m(k)-e_m(0)|&\leq\frac{\vep^2\Cg\beta\,a(k)}{(1+m\vep)^{(d-2)/2}},\\
\lbeq{g1bd}|g_m(k)|&\leq\frac{\vep^2\Cg\beta}{(1+m\vep)^{d/2}},&
|\nabla^2g_m(0)|&\leq\frac{\vep^2\Cg\sigma^2\beta}{(1+m\vep)^{(d-2)/2}},
\end{align}
\vspace{-2pc}
\begin{align}
\lbeq{g3bd}\Big|g_m(k)-g_m(0)-\frac{a(k)}{\sigma^2}\,\nabla^2g_m(0)\Big|
&\leq\frac{\vep^2\Cg\beta\,a(k)^{1+\Delta'}}{(1+m\vep)^{(d-2)/2-\Delta'}},\\
\lbeq{g2bd}|\partial_\lamb g_m(0)|&\leq\frac{\vep^2\Cg\beta}{(1+m\vep)^{(d
 -2)/2}}.
\end{align}
We note that $\Cg$ depends on $K$ and that, by
Lemmas~\ref{lem:fabs}--\ref{lem:fD2}, $K$ depends only on $K_4$ when
$\beta\ll1$.  Therefore, we can choose $\Cg$ large depending only on
$K_4$ when $\beta\ll1$.

\paragraph{Advancement of (H1).}
By \refeq{lambn-def} and the mean-value theorem,
\begin{align}
\lamb_{n+1}-\lamb_n&=-\frac1\vep\,g_{n+1}(0;\lamb_n)-\frac1\vep\sum_{m=2}^n
 [g_m(0;\lamb_n)-g_m(0;\lamb_{n-1})]\nnmb\\
&=-\frac1\vep\,g_{n+1}(0;\lamb_n)-\frac{\lamb_n-\lamb_{n-1}}\vep\sum_{m=2}^n
 \partial_\lamb g_m(0;\lamb_*),\lbeq{lambndiff}
\end{align}
for some $\lamb_*$ between $\lamb_n$ and $\lamb_{n-1}$.  Since
$\lamb_{n-1}\in I_n$ (cf., \refeq{In-def} and \refeq{H12}),
$\lamb_*$ is also in $I_n$.  By \refeq{g1bd}, \refeq{g2bd} and (H1),
\begin{align}
|\lamb_{n+1}-\lamb_n|\leq\frac{\vep\Cg\beta}{[1+(n+1)\vep]^{d/2}}+|\lamb_n
 -\lamb_{n-1}|\vep\sum_{m=2}^n\frac{\Cg\beta}{(1+m\vep)^{(d-2)/2}}\leq\frac
 {\vep\Cg(1+cK_1\beta)\beta}{[1+(n+1)\vep]^{d/2}}.
\end{align}
Therefore, (H1) holds for $n+1$, if $\beta\ll1$ and $K_1>\Cg$.
\qed

\paragraph{Advancement of (H2).}
Let $1+M_n$ be the denominator of \refeq{vn-def}, and let $N_n$ be the
numerator of \refeq{vn-def}.  Then,
\begin{align}\lbeq{prevdiff}
v_{n+1}-v_n=\frac{\frac{-1}{\sigma^2\vep}\nabla^2g_{n+1}(0)}{1+M_{n+1}}
 -\frac{N_n\,n\,g_{n+1}(0)}{(1+M_{n+1})(1+M_n)}.
\end{align}
By \refeq{g1bd}, we obtain that, for $m\leq n+1$,
\begin{align}\lbeq{Nnbd}
|M_m|\leq\vep\sum_{l=2}^{m}\frac{(l-1)\vep\Cg\beta}{(1+l\vep)^{d/2}}
 \leq c\Cg\beta,&&
|N_m-\lamb|\leq\vep\sum_{l=2}^m\frac{\Cg\beta}{(1+l\vep)^{(d-2)/2}}
 \leq c\Cg\beta,
\end{align}
and
\begin{align}
\Big|\frac{-1}{\sigma^2\vep}\,\nabla^2g_{n+1}(0)\Big|\leq\frac{\vep\Cg
 \beta}{[1+(n+1)\vep]^{(d-2)/2}},&&
|n\,g_{n+1}(0)|\leq\frac{n\vep^2\Cg\beta}{[1+(n+1)\vep]^{d/2}}.
\end{align}
Therefore,
\begin{align}
&|v_{n+1}-v_n|\leq\frac{\vep\Cg\beta}{(1-c\Cg\beta)[1+(n+1)\vep]^{(d-2)/2}}
 +\frac{(\lamb+c\Cg\beta)n\vep^2\Cg\beta}{(1-c\Cg\beta)^2[1+(n+1)\vep]
 ^{d/2}}\\
&\quad=\frac{1-c\Cg\beta+(\lamb+c\Cg\beta)\frac{n\vep}{1+(n+1)\vep}}{(1-c
 \Cg\beta)^2}\,\frac{\vep\Cg\beta}{[1+(n+1)\vep]^{(d-2)/2}}\leq\frac{1+
 \lamb}{(1-c\Cg\beta)^2}\,\frac{\vep\Cg\beta}{[1+(n+1)\vep]^{(d-2)/2}}.\nn
\end{align}
Since $\lamb\in I_{n+1}$, (H2) holds for $n+1$, if $\beta\ll1$ and
$K_2>2\Cg$.
\qed

\paragraph{Advancement of (H3).}
First, we derive expressions for $r_{n+1}(0)$ and $r_{n+1}(k)-r_{n+1}(0)$.
By dividing both sides of \refeq{f-renew} by $f_n(k)$ and using
$g_1(k)=1-\vep+\lamb\vep\wD(k)$,
\begin{align}
&\frac{f_{n+1}(k)}{f_n(k)}=g_1(k)+\sum_{m=1}^ng_{m+1}(k)\frac{f_{n-m}(k)}
 {f_n(k)}+\frac{e_{n+1}(k)}{f_n(k)}\nn\\
&\qquad=1-\vep v_{n+1}a(k)+\vep\bigg[v_{n+1}a(k)-1+\lamb\wD(k)+\frac1\vep
 \sum_{m=1}^ng_{m+1}(k)\frac{f_{n-m}(k)}{f_n(k)}+\frac{e_{n+1}(k)}{\vep
 f_n(k)}\bigg].\lbeq{rnk-def}
\end{align}
Therefore, $r_{n+1}(k)$ equals the expression in the above brackets.
In particular,
\begin{align}\lbeq{rn0}
r_{n+1}(0)&=-1+\lamb+\frac1\vep\sum_{m=1}^ng_{m+1}(0)\frac{f_{n-m}(0)}{f_n
 (0)}+\frac{e_{n+1}(0)}{\vep f_n(0)}\nn\\
&=\bigg[\lamb-1+\frac1\vep\sum_{m=2}^{n+1}g_m(0)\bigg]+\frac1\vep\sum_{m=2}
 ^{n+1}g_m(0)\bigg[\frac{f_{n+1-m}(0)}{f_n(0)}-1\bigg]+\frac{e_{n+1}(0)}
 {\vep f_n(0)}\\[7pt]
&=r_{n+1}^{\sss(1)}(0)+r_{n+1}^{\sss(2)}(0)+r_{n+1}^{\sss(3)}(0),\nn
\end{align}
where we denote the first, second and third terms in \refeq{rn0} by
$r_{n+1}^{\sss(1)}(0)$, $r_{n+1}^{\sss(2)}(0)$ and $r_{n+1}^{\sss(3)}(0)$,
respectively.  Similarly, we can obtain an expression for
$r_{n+1}(k)-r_{n+1}(0)$.  To do so, we note that, by \refeq{vn-def},
\begin{align}
v_{n+1}=\lamb-\frac1{\sigma^2\vep}\sum_{m=2}^{n+1}\nabla^2g_m(0)-v_{n+1}
 \sum_{m=2}^{n+1}(m-1)\,g_m(0).
\end{align}
Using this identity, we obtain
\begin{align}\lbeq{rn-rn0}
r_{n+1}(k)-r_{n+1}(0)&=(v_{n+1}-\lamb)\,a(k)+\frac1\vep\sum_{m=2}^{n+1}
 [g_m(k)-g_m(0)]\frac{f_{n+1-m}(k)}{f_n(k)}\nn\\
&\qquad+\frac1\vep\sum_{m=2}^{n+1}g_m(0)\bigg[\frac{f_{n+1-m}(k)}{f_n(k)}
 -\frac{f_{n+1-m}(0)}{f_n(0)}\bigg]+\frac1\vep\bigg[\frac{e_{n+1}(k)}
 {f_n(k)}-\frac{e_{n+1}(0)}{f_n(0)}\bigg]\nn\\
&=\frac1\vep\sum_{m=2}^{n+1}\bigg[[g_m(k)-g_m(0)]\,\frac{f_{n+1-m}(k)}
 {f_n(k)}-\frac{a(k)}{\sigma^2}\,\nabla^2g_m(0)\bigg]\nn\\
&\qquad+\frac1\vep\sum_{m=2}^{n+1}g_m(0)\bigg[\frac{f_{n+1-m}(k)}{f_n(k)}
 -\frac{f_{n+1-m}(0)}{f_n(0)}-\vep v_{n+1}\,(m-1)\,a(k)\bigg]\nn\\
&\qquad+\frac1\vep\,\bigg[\frac{e_{n+1}(k)}{f_n(k)}-\frac{e_{n+1}(0)}
 {f_n(0)}\bigg]\\[7pt]
&=\Delta r_{n+1}^{\sss(1)}(k)+\Delta r_{n+1}^{\sss(2)}(k)+\Delta r_{n+1}
 ^{\sss(3)}(k),\nn
\end{align}
where we denote the first, second and third terms in \refeq{rn-rn0} by
$\Delta r_{n+1}^{\sss(1)}(k)$, $\Delta r_{n+1}^{\sss(2)}(k)$ and
$\Delta r_{n+1}^{\sss(3)}(k)$, respectively.

Therefore, to advance (H3), we are left to investigate
$r_{n+1}^{\sss(i)}(0)$ and $\Delta r_{n+1}^{\sss(i)}(k)$ for $i=1,2,3$.

\begin{proof}[\it Advancement of the first inequality in \refeq{H3}]
We recall that $r_{n+1}(0)$ has been decomposed, as in \refeq{rn0}, into
$r_{n+1}^{\sss(i)}(0)$ for $i=1,2,3$.  First, we investigate
$r_{n+1}^{\sss(1)}(0)$.  By \refeq{lambn-def} and the mean-value theorem,
we have
\begin{align}
|r_{n+1}^{\sss(1)}(0)|&\leq|\lamb-\lamb_n|+|\lamb_n-\lamb_{n+1}|+\bigg|
 \lamb_{n+1}-1+\frac1\vep\sum_{m=2}^{n+1}g_m(0;\lamb)\bigg|\nn\\
&=|\lamb_n-\lamb_{n+1}|+|\lamb-\lamb_n|+\bigg|\frac1\vep\sum_{m=2}^{n+1}
 [g_m(0;\lamb)-g_m(0;\lamb_n)]\bigg|\nn\\
&\leq|\lamb_n-\lamb_{n+1}|+|\lamb-\lamb_n|\bigg[1+\frac1\vep\sum_{m=2}
 ^{n+1}|\partial_\lamb g_m(0;\lamb_*)|\bigg],
\end{align}
for some $\lamb_*$ between $\lamb$ and $\lamb_n$.  Since
$\lamb\in I_{n+1}\subset I_n$, $\lamb_*$ is also in $I_n$.
By \refeq{H12}, \refeq{In-def} and \refeq{g2bd},
\begin{align}\lbeq{rn0-1stbd}
|r_{n+1}^{\sss(1)}(0)|\leq\frac{\vep K_1\beta}{[1+(n+1)\vep]^{d/2}}+\frac
 {K_1\beta}{(1+n\vep)^{(d-2)/2}}\Bigg[1+\sum_{m=2}^{n+1}\frac{\vep\Cg\beta}
 {(1+m\vep)^{(d-2)/2}}\Bigg]\leq\frac{cK_1\beta}{[1+(n+1)\vep]^{(d-2)/2}}.
\end{align}
Therefore, we need $K_3\gg K_1$.

Next we investigate $r_{n+1}^{\sss(2)}(0)$.  We will use the following
results of Taylor's theorem applied to $h(t)=\prod_i(1+c_it)^{-1}$ with
$|c_i|<1$ for all $i$:
\begin{align}
&|h(1)-h(0)|\leq\sup_{t\in(0,1)}|h'(t)|\leq\sum_i\frac{|c_i|}{1-|c_i|}\,
 e^{\sum_j\frac{|c_j|}{1-|c_j|}},\lbeq{Taylor1}\\
&|h(1)-h(0)-h'(0)|\leq\frac12\sup_{t\in(0,1)}|h''(t)|\leq\bigg(\sum_i\frac
 {|c_i|}{1-|c_i|}\bigg)^2\,e^{\sum_j\frac{|c_j|}{1-|c_j|}}\lbeq{Taylor2},
\end{align}
By \refeq{f0sl-def}, \refeq{g1bd} and \refeq{Taylor1},
\begin{align}\lbeq{rn0-2nd}
|r_{n+1}^{\sss(2)}(0)|\leq\frac1\vep\sum_{m=2}^{n+1}|g_m(0)|\bigg|\prod_{l
 =n+2-m}^n[1+\vep\,r_l(0)]^{-1}-1\bigg|\leq\sum_{m=2}^{n+1}\frac{\vep\Cg
 \beta}{(1+m\vep)^{d/2}}~\phi_m\,e^{\phi_m},
\end{align}
where, by \refeq{H3},
\begin{align}\lbeq{phi-def}
\phi_m=\sum_{l=n+2-m}^n\frac{\vep|r_l(0)|}{1-\vep|r_l(0)|}\leq
 ~\vep\!\sum_{l=n+2-m}^n\frac{cK_3\beta}{(1+l\vep)^{(d-2)/2}},
\end{align}
and thus $e^{\phi_m}\leq e^{cK_3\beta}$ for all $m\leq n+1$.  Substituting
\refeq{phi-def} into \refeq{rn0-2nd} and using Lemma~\ref{lem:conv} with
$(a,b)=(\frac{d}2,\frac{d-2}2)$, we obtain
\begin{align}\lbeq{rn0-2ndbd}
|r_{n+1}^{\sss(2)}(0)|\leq\frac{c\Cg K_3\beta^2}{[1+(n+1)\vep]^{(d-2)/2}}.
\end{align}

Finally, we investigate $r_{n+1}^{\sss(3)}(0)$.  As in \refeq{rn0-2nd},
$|f_n(0)^{-1}-1|$ is bounded by
\begin{align}\lbeq{fn0bd}
\bigg|\prod_{l=1}^n[1+\vep r_l(0)]^{-1}-1\bigg|\leq\phi_{n+1}\,
 e^{\phi_{n+1}}\leq cK_3\beta.
\end{align}
Using \refeq{ebd}, we obtain
\begin{align}\lbeq{rn0-3rdbd}
|r_{n+1}^{\sss(3)}(0)|\leq\frac{\vep\Cg(1+cK_3\beta)\beta}
 {[1+(n+1)\vep]^{d/2}}.
\end{align}

The advancement of the first inequality in \refeq{H3} is now completed by
\refeq{rn0}, \refeq{rn0-1stbd}, \refeq{rn0-2ndbd} and \refeq{rn0-3rdbd},
if $\beta\ll1$ and $K_3\gg K_1\vee\Cg$.
\end{proof}

\begin{proof}[\it Advancement of the second inequality in \refeq{H3}]
Recall that $k\in\cA_{n+1}$, and that $r_{n+1}(k)-r_{n+1}(0)$ has been
decomposed, as in \refeq{rn-rn0}, into $\Delta r_{n+1}^{\sss(i)}(k)$ for
$i=1,2,3$.

First, we investigate $\Delta r_{n+1}^{\sss(1)}(k)$, which is bounded as
\begin{align}
|\Delta r_{n+1}^{\sss(1)}(k)|&\leq\frac1\vep\sum_{m=2}^{n+1}\bigg|
 g_m(k)-g_m(0)-\frac{a(k)}{\sigma^2}\,\nabla^2g_m(0)\bigg|+\frac1\vep\sum
 _{m=2}^{n+1}|g_m(k)-g_m(0)|\bigg|\frac{f_{n+1-m}(0)}{f_n(0)}-1\bigg|\nn\\
&\qquad+\frac1\vep\sum_{m=2}^{n+1}|g_m(k)-g_m(0)|\bigg|\frac{f_{n+1-m}(k)}
 {f_n(k)}-\frac{f_{n+1-m}(0)}{f_n(0)}\bigg|.\lbeq{rn-rn0-pre1stbd}
\end{align}
By \refeq{g3bd} with $\delta<\Delta'<\frac{d-4}2$, the first sum is bounded
by
\begin{align}\lbeq{1stbd-pre1stbd}
\vep\sum_{m=2}^{n+1}\frac{\Cg\beta\,a(k)^{1+\Delta'}}{(1+m\vep)^{(d-2)/2-
 \Delta'}}\leq c\Cg\beta\,a(k)\bigg[\frac{\log[2+(n+1)\vep]}{[1+(n+1)\vep]}
 \bigg]^{\Delta'}\leq\frac{c\Cg\beta\,a(k)}{[1+(n+1)\vep]^\delta},
\end{align}
while the second sum in \refeq{rn-rn0-pre1stbd} is first bounded similarly
to \refeq{rn0-2nd}, and then bounded, by using \refeq{g3bd} with $\Delta'=0$
and \refeq{phi-def}, as well as Lemma~\ref{lem:conv} with $a=b=\frac{d-2}2$,
by
\begin{align}\lbeq{2ndbd-pre1stbd}
\sum_{m=2}^{n+1}\frac{2\vep\Cg\beta\,a(k)}{(1+m\vep)^{(d-2)/2}}\sum_{l=n+2
 -m}^n\frac{\vep cK_3\beta}{(1+l\vep)^{(d-2)/2}}\leq\frac{c\Cg K_3\beta^2
 a(k)}{[1+(n+1)\vep]^{(d-2)/2\wedge(d-4)}}\leq\frac{c\Cg K_3\beta^2a(k)}
 {[1+(n+1)\vep]^{2\delta}},
\end{align}
where we use
$\frac{d-2}2\wedge(d-4)=\frac{d-4}2+1\wedge\frac{d-4}2\geq2\delta$.  By
using \refeq{g3bd} with $\Delta'=0$ again and \refeq{Taylor1}, the third
sum in \refeq{rn-rn0-pre1stbd} is bounded similarly to \refeq{rn0-2nd} by
\begin{gather}
\sum_{m=2}^{n+1}\frac{2\vep\Cg\beta\,a(k)}{(1+m\vep)^{(d-2)/2}}
 \,\bigg|\frac{f_{n+1-m}(0)}{f_n(0)}\bigg|\,\bigg|\prod_{l=n+2-
 m}^n[1-\vep v_l\,a(k)+\vep s_l(k)]^{-1}-1\bigg|\nnmb\\
\leq\sum_{m=2}^{n+1}\frac{2\vep\Cg\beta\,a(k)}{(1+m\vep)^{(d-2)
 /2}}\,\big(1+\phi_m\,e^{\phi_m}\big)\,\psi_m(k)\,
 e^{\psi_m(k)},\lbeq{3rd-pre1stbd}
\end{gather}
where $\phi_m\,e^{\phi_m}\leq cK_3\beta$ as discussed below
\refeq{phi-def}, and
\begin{align}\lbeq{psi-def}
\psi_m(k)=\sum_{l=n+2-m}^n\frac{\vep[v_l a(k)+|s_l (k)|]}{1-\vep[v_l
 a(k)+|s_l(k)|]}.
\end{align}
By \refeq{sl-bd} and \refeq{vm-bd},
\begin{align}
v_la(k)+|s_l(k)|\leq\bigg[v_l+\frac{(1+\vep v_l)K_3\beta}{(1-\vep K_3
 \beta)(1+l\vep)^\delta}\bigg]a(k)\leq[1+c\,(K_1+K_2+K_3)\beta]\,a(k)
 \equiv q\,a(k).\lbeq{q-def}
\end{align}
Since $k\in\cA_{n+1}$, $\psi_m(k)$ is bounded by
\begin{align}\lbeq{psibd}
\psi_m(k)\leq\frac{(m-1)\vep q\,a(k)}{1-\vep q\,a(k)}\leq[1+c\vep\,a(k)]
 (m-1)\vep q\,a(k),
\end{align}
which is further bounded by $\gamma q\big[1+c\vep\frac{\log[2+(n+1)\vep]}
{1+(n+1)\vep}\big]\log[2+(n+1)\vep]$, and hence
\begin{align}\lbeq{epsibd}
e^{\psi_m(k)}\leq c\,e^{\gamma q\log[2+(n+1)\vep]} \leq c\,[1+(n+1)
 \vep]^{\gamma q}.
\end{align}
Substituting \refeq{psibd}--\refeq{epsibd} into \refeq{3rd-pre1stbd},
and using $a(k)\leq\gamma\frac{\log[2+(n+1)\vep]}{1+(n+1)\vep}$ and
$\gamma q+\delta<1\wedge\frac{d-4}2$ for $\beta\ll1$ (cf., \refeq{drgd}
and \refeq{q-def}), we can bound \refeq{3rd-pre1stbd} by
\begin{align}\lbeq{3rdbd-pre1stbd}
c\Cg\beta\,a(k)\frac{\log[2+(n+1)\vep]}{[1+(n+1)\vep]^{1-\gamma q}}~\vep
 \sum_{m=2}^{n+1}\frac{(m-1)\vep}{(1+m\vep)^{(d-2)/2}}\leq\frac{c\Cg\beta
 \,a(k)}{[1+(n+1)\vep]^\delta}.
\end{align}
By \refeq{rn-rn0-pre1stbd}--\refeq{2ndbd-pre1stbd} and
\refeq{3rdbd-pre1stbd}, if $\beta\ll1$ and $K_3\gg\Cg$, we obtain
\begin{align}
|\Delta r_{n+1}^{\sss(1)}(k)|\leq\frac{\frac13K_3\beta a(k)}{[1+(n+1)
 \vep]^\delta}.\lbeq{rn-rn0-1stbd}
\end{align}

Next, we investigate $|\Delta r_{n+1}^{\sss(2)}(k)|$, which is bounded, by
using \refeq{g1bd} and \refeq{rn0-2nd}, as
\begin{align}
|\Delta r_{n+1}^{\sss(2)}(k)|&\leq\frac1\vep\sum_{m=2}^{n+1}|g_m(0)|\bigg|
 \frac{f_{n+1-m}(0)}{f_n(0)}\bigg|\bigg|\prod_{l=n+2-m}^n[1-\vep v_la(k)
 +\vep s_l(k)]^{-1}-1-(m-1)\vep v_{n+1}a(k)\bigg|\nn\\
&\qquad+\frac1\vep\sum_{m=2}^{n+1}|g_m(0)|\bigg|\frac{f_{n+1-m}(0)}{f_n(0)}
 -1\bigg|\,(m-1)\vep v_{n+1}a(k)\nn\\
&\leq\vep\sum_{m=2}^{n+1}\frac{\Cg\beta(1+\phi_me^{\phi_m})}{(1+m\vep)^{d
 /2}}\bigg|\prod_{l=n+2-m}^n[1-\vep v_la(k)+\vep s_l(k)]^{-1}-1-(m-1)\vep
 v_{n+1}a(k)\bigg|\nn\\
&\qquad+\vep\sum_{m=2}^{n+1}\frac{\Cg\beta v_{n+1}a(k)\,(m-1)\vep}{(1+m\vep)
 ^{d/2}}~\phi_me^{\phi_m}.\lbeq{rn-rn0-2brackets}
\end{align}
Using \refeq{vm-bd}, \refeq{phi-def}, Lemma~\ref{lem:conv} with
$a=b=\frac{d-2}2$ and $\delta<1\wedge\frac{d-4}2$, we can bound the
second sum by $c\Cg\beta a(k)[1+(n+1)\vep]^{-2\delta}$.  The first sum
in \refeq{rn-rn0-2brackets} is bounded, by using \refeq{Taylor2}, by
\begin{align}
&\vep\sum_{m=2}^{n+1}\frac{c\Cg\beta}{(1+m\vep)^{d/2}}\Bigg[\bigg|
 \prod_{l=n+2-m}^n[1-\vep v_la(k)+\vep s_l(k)]^{-1}-1-\sum_{l=n+2-m}^n
 \vep[v_la(k)-s_l(k)]\bigg|\nn\\
&\hspace{10pc}+\bigg|\sum_{l=n+2-m}^n\vep[(v_l-v_{n+1})a(k)-s_l(k)]
 \bigg|\Bigg]\nn\\
&\quad\leq\vep\sum_{m=2}^{n+1}\frac{c\Cg\beta}{(1+m\vep)^{d/2}}\Bigg[
 \psi_m(k)^2e^{\psi_m(k)}+\sum_{l=n+2-m}^n\vep\bigg[\sum_{j=l+1}^{n+1}
 |v_j-v_{j-1}|a(k)+|s_l(k)|\bigg]\Bigg].
\end{align}
Similarly to \refeq{3rdbd-pre1stbd}, the contribution from $\psi_m(k)^2
e^{\psi_m(k)}$ is bounded by $c\Cg\beta a(k)[1+(n+1)]^{-\delta}$.  By
\refeq{H12}, \refeq{sl-bd} and Lemma~\ref{lem:conv} with $a=\frac{d}2$
and $b=\delta(<\frac{d-4}2)$, the other contribution is
bounded by
\begin{align}
\vep\sum_{m=2}^{n+1}\frac{c\Cg\beta}{(1+m\vep)^{d/2}}\sum_{l=n+2-m}^n\vep
 \bigg[\frac{cK_2\beta a(k)}{(1+l\vep)^{(d-4)/2}}+\frac{cK_3\beta a(k)}
 {(1+l\vep)^\delta}\bigg]\leq\frac{c\Cg(K_2+K_3)\beta^2a(k)}{[1+(n+1)\vep]
 ^\delta}.
\end{align}
Therefore, if $\beta\ll1$ and $K_3\gg\Cg$, we obtain
\begin{align}
|\Delta r_{n+1}^{\sss(2)}(k)|\leq\frac{\frac13K_3\beta a(k)}{[1+(n+1)
 \vep]^\delta}.\lbeq{rn-rn0-2ndbd}
\end{align}

Finally, we investigate $|\Delta r_{n+1}^{\sss(3)}(k)|$, which is bounded as
\begin{align}
|\Delta r_{n+1}^{\sss(3)}(k)|\leq\bigg|\frac{e_{n+1}(k)}{\vep f_n(0)}\bigg|
 \bigg|\frac{f_n(0)}{f_n(k)}-1\bigg|+\bigg|\frac{e_{n+1}(k)-e_{n+1}(0)}{\vep
 f_n(0)}\bigg|.
\end{align}
By \refeq{fn0bd}, $|f_n(0)|\geq1-cK_3\beta$.  As in \refeq{3rd-pre1stbd},
$\big|\frac{f_n(0)}{f_n(k)}-1\big|$ is bounded, by using
\refeq{psibd}--\refeq{epsibd}, by
$\psi_{n+1}(k)e^{\psi_{n+1}(k)}\leq c[1+(n+1)\vep]^{1+\gamma q}a(k)$.
Therefore, using \refeq{ebd} and taking $\beta$ sufficiently small such
that $\gamma q+\delta<\frac{d-4}2$ (cf., \refeq{drgd} and \refeq{q-def}),
we obtain
\begin{align}\lbeq{rn-rn0-3rdbd}
|\Delta r_{n+1}^{\sss(3)}(k)|\leq\frac{c\vep\Cg\beta a(k)}{[1+(n+1)\vep]^{
 (d-2)/2-\gamma q}}\leq\frac{c\vep\Cg\beta a(k)}{[1+(n+1)\vep]^{\delta+1}}
 \leq\frac{\frac13K_3\beta a(k)}{[1+(n+1)\vep]^\delta},
\end{align}
if $K_3\gg\Cg$.

The advancement of the second inequality in \refeq{H3} is now completed by
\refeq{rn-rn0-1stbd}, \refeq{rn-rn0-2ndbd} and \refeq{rn-rn0-3rdbd}, if
$\beta\ll1$ and $K_3\gg\Cg$.
\end{proof}

\paragraph{Advancement of (H4).}
To advance (H4), we rewrite \refeq{f-renew} as
\begin{align}
f_{n+1}(k)=\bigg[g_1(k)+\sum_{m=2}^{n+1}g_m(k)\bigg]f_n(k)+W_{n+1}(k)
 +e_{n+1}(k),
\end{align}
where
\begin{align}\lbeq{W-def}
W_{n+1}(k)=\sum_{m=2}^{n+1}g_m(k)[-f_n(k)+f_{n+1-m}(k)]=\sum_{m=2}^{n+1}
 g_m(k)\sum_{l=n+2-m}^n[f_{l-1}(k)-f_l(k)].
\end{align}
Furthermore, using
$g_1(k)=1-\vep+\lamb\vep\hat D(k)=1-\lamb\vep a(k)+(\lamb-1)\vep$, we have
\begin{align}\lbeq{H4-fncoeff}
g_1(k)+\sum_{m=2}^{n+1}g_m(k)&=1-\bigg[\lamb-\frac1{\sigma^2\vep}
 \sum_{m=2}^{n+1}\nabla^2g_m(0)\bigg]\vep a(k)+\vep\bigg[\lamb-1
 +\frac1\vep\sum_{m=2}^{n+1}g_m(0)\bigg]\nn\\
&\qquad+\sum_{m=2}^{n+1}\bigg[g_m(k)-g_m(0)-\frac{a(k)}{\sigma^2}
 \nabla^2g_m(0)\bigg]\\
&=1-N_{n+1}\vep a(k)+\vep r_{n+1}^{\sss(1)}(0)+X_{n+1}(k),\nn
\end{align}
where we recall $N_n$ and $r_{n+1}^{\sss(1)}(0)$ in \refeq{prevdiff} and
\refeq{rn0}, respectively, and denote the last sum in \refeq{H4-fncoeff}
by $X_{n+1}(k)$.  Therefore,
\begin{align}\lbeq{fn-rewritten}
f_{n+1}(k)=f_n(k)\,\big[1-N_{n+1}\,\vep\,a(k)+\vep r_{n+1}^{\sss(1)}(0)
 +X_{n+1}(k)\big]+W_{n+1}(k)+e_{n+1}(k).
\end{align}
We have already obtained $|N_{n+1}-\lamb|\leq c\Cg\beta$ in \refeq{Nnbd}
and $|r_{n+1}^{\sss(1)}(0)|\leq cK_1\beta[1+(n+1)\vep]^{-(d-2)/2}$ in
\refeq{rn0-1stbd}, while $X_{n+1}(k)$ equals $\vep$ times the first sum of
\refeq{rn-rn0-pre1stbd} and is bounded, by using the leftmost expression
of \refeq{1stbd-pre1stbd} with $\Delta'<\frac{d-4}2$, by
$c\vep\Cg\beta\,a(k)^{1+\Delta'}$.  We prove below that,
for $k\notin\cA_{n+1}$,
\begin{align}\lbeq{Wbd}
|W_{n+1}(k)|\leq\vep\frac{c\Cg(1+K_3\beta+K_5)\beta\,a(k)^{-1-\rho}}
 {[1+(n+1)\vep]^{d/2}}.
\end{align}

Assuming \refeq{Wbd}, we first advance the second inequality in
\refeq{H4}, and then advance the first inequality in \refeq{H4}.
To advance these inequalities, we will use the first inequality
in \refeq{H4} for $m=n$ in the extended region
$\cA_{n+1}^{\rm c}=\cA_n^{\rm c}\Dcup{}{}(\cA_n\setminus\cA_{n+1})$.  We
now verify the use of this inequality for $k\in\cA_n\setminus\cA_{n+1}$.
When $n\vep\leq T$ for some large $T$, we can choose $K_4\gg1$ (depending
on $T$) such that, for all $k\in[-\pi,\pi]^d$,
\begin{align}\lbeq{fn-naivebd}
|f_n(k)|\leq\|\tau_{n\vep}\|_1\leq\|p_\vep^{*n}\|_1=(1-\vep+\lamb\vep)^n
 \leq e^{(\lamb-1)n\vep}\leq\frac{2^{-2-\rho}K_4}{(1+n\vep)^{d/2}}\leq
 \frac{K_4a(k)^{-2-\rho}}{(1+n\vep)^{d/2}}.
\end{align}
When $n\vep>T$, we use Lemma~\ref{lem:fabs} and
$k\in\cA_n\setminus\cA_{n+1}$ (so that $\gamma\frac{\log[2+(n+1)\vep]}
{1+(n+1)\vep}<a(k)\leq\gamma\frac{\log(2+n\vep)}{1+n\vep}$) to obtain
\begin{align}\lbeq{H3toH4}
|f_n(k)|\leq ce^{-n\vep qa(k)}&\leq c(2+n\vep)^{-\frac{n\vep}{1+(n+1)
 \vep}\frac{\log[2+(n+1)\vep]}{\log(2+n\vep)}q\gamma}\nn\\
&\leq c(1+n\vep)^{-q'\gamma}=\frac{c}{(1+n\vep)^{d/2}}\frac{(1+n\vep)
 ^{2+\rho}}{(1+n\vep)^{q'\gamma-[\frac{d}2-(2+\rho)]}}\leq\frac{K_4a
 (k)^{-2-\rho}}{(1+n\vep)^{d/2}},
\end{align}
if $K_4\gg1$, where we use $q'\gamma>\frac{d}2-(2+\rho)$ for $\beta\ll1$
and $T\gg1$ (cf., \refeq{drgd} and \refeq{q-def}).

Therefore, by using the first inequality in \refeq{H4} with $m=n$ for
$k\notin\cA_{n+1}$, together with \refeq{lamb-bd} and
\refeq{fn-rewritten}--\refeq{Wbd}, we obtain
\begin{align}\lbeq{H4-2nd}
|f_{n+1}(k)-f_n(k)|&\leq\vep\frac{K_4a(k)^{-2-\rho}}{(1+n\vep)^{d/2}}\bigg[
 (\lamb+c\Cg\beta)a(k)+\frac{cK_1\beta}{[1+(n+1)\vep]^{(d-2)/2}}+c\Cg\beta
 a(k)^{1+\Delta'}\bigg]\nn\\
&\quad+\vep\frac{c\Cg(1+K_3\beta+K_5)\beta a(k)^{-1-\rho}}{[1+(n+1)\vep]^{d
 /2}}+\frac{\vep^2\Cg\beta}{[1+(n+1)\vep]^{d/2}}\nn\\
&\leq\vep\frac{cK_4[1+O(\beta)]\,a(k)^{-1-\rho}}{[1+(n+1)\vep]^{d/2}}+\vep
 \frac{O(\beta)\,a(k)^{-1-\rho}}{[1+(n+1)\vep]^{d/2}}+\vep^2\frac{O(\beta)\,
 a(k)^{-1-\rho}}{[1+(n+1)\vep]^{d/2}},
\end{align}
where we use
$[1+(n+1)\vep]^{-(d-2)/2}\leq a(k)^{(d-2)/2}\leq 2^{(d-4)/2}a(k)$ to
obtain the first term, and use $2^{-1-\rho}\leq a(k)^{-1-\rho}$ for the
third term.  This completes the advancement of the second inequality in
\refeq{H4}, if $\beta\ll1$ and $K_5\gg K_4$, under the hypotheses that
\refeq{Wbd} holds for $k\notin\cA_{n+1}$.

Since \refeq{fn-naivebd} holds for $n\leq T/\vep$ independently of $k$,
it remains to advance the first inequality of (H4) for $n>T/\vep$.
Similarly to \refeq{H4-2nd}, we have
\begin{align}
|f_{n+1}(k)|&\leq\frac{K_4a(k)^{-2-\rho}}{(1+n\vep)^{d/2}}\bigg[|1-N_{n+1}
 \vep a(k)|+\frac{c\vep K_1\beta}{[1+(n+1)\vep]^{(d-2)/2}}+c\vep \Cg\beta
 a(k)^{1+\Delta'}\bigg]\nn\\
&\quad+\frac{c\vep\Cg(1+K_3\beta+K_5)\beta a(k)^{-1-\rho}}{[1+(n+1)\vep]
 ^{d/2}}+\frac{\vep^2\Cg\beta}{[1+(n+1)\vep]^{d/2}}.\lbeq{H4-1st}
\end{align}
Again, by $a(k)\leq2$, the sum of the last two terms is bounded by
$\vep O(\beta)\,a(k)^{-2-\rho}[1+(n+1)\vep]^{-d/2}$.  To prove the first
inequality in \refeq{H4} with $m=n+1$, it thus suffices to show that
\begin{align}\lbeq{H4strineq}
\bigg[\frac{1+(n+1)\vep}{1+n\vep}\bigg]^{d/2}\bigg[|1-N_{n+1}\vep a(k)|
 +\frac{c\vep K_1\beta}{[1+(n+1)\vep]^{(d-2)/2}}+c\vep\Cg\beta a(k)^{1+
 \Delta'}\bigg]<1.
\end{align}
To achieve this inequality uniformly in $\vep\leq1$, we consider the case
in which $a(k)\leq1/2$ and the other case in which $a(k)>1/2$ separately.

When $a(k)\leq1/2$, since $N_{n+1}=1+O(\beta)$ (cf., \refeq{lamb-bd} and
\refeq{Nnbd}), we have $|1-N_{n+1}\vep a(k)|=1-N_{n+1}\vep a(k)$ for
$\beta\ll1$.  Using $a(k)^{\Delta'}\leq2^{\Delta'}$ and
then $a(k)>\gamma\frac{\log[2+(n+1)\vep]}{1+(n+1)\vep}$, we can bound
\refeq{H4strineq} by
\begin{align}
&\Big(1+\frac{c\vep}{1+n\vep}\Big)\bigg[1-(1-c\beta)\vep a(k)+\frac{c\vep
 \beta}{[1+(n+1)\vep]^{(d-2)/2}}\bigg]\\
&\quad\leq1-\vep\bigg[(1-c\beta)\frac{\gamma\log[2+(n+1)\vep]}{1+(n+1)\vep}
 -\frac{c}{1+n\vep}-\Big(1+\frac{c\vep}{1+n\vep}\Big)\frac{c\beta}{[1+(n+1)
 \vep]^{(d-2)/2}}\bigg]<1,\nn
\end{align}
if $\beta\ll1$ and $T\gg1$.

Since the above argument also applies to the case in which
$1-\vep a(k)>1-(2-\eta)\vep>0$ (and $\beta\ll1$, depending on $\eta$),
it thus remains to consider the other case in which $1-(2-\eta)\vep\leq0$
and $a(k)>1/2$.  In this case, since $\vep\leq1$, we have
\begin{align}
|1-\vep\,a(k)|\leq[(2-\eta)\vep-1]\vee\Big(1-\frac\vep2\Big)
 \leq1-\Big(\eta\wedge\frac\vep2\Big).
\end{align}
Since $N_{n+1}=1+O(\beta)$, \refeq{H4strineq} is bounded by
\begin{gather}
\Big(1+\frac{c\vep}{1+n\vep}\Big)\bigg[1-\Big(\eta\wedge\frac\vep2\Big)
 +c\vep\beta a(k)+\frac{c\vep\beta}{[1+(n+1)\vep]^{(d-2)/2}}\bigg]\nn\\
\quad\leq1-\bigg[\Big(\eta\wedge\frac\vep2\Big)-\frac{c\vep}{1+n\vep}
 -\Big(1+\frac{c\vep}{1+n\vep}\Big)c\vep\beta\bigg]<1,
\end{gather}
if $\beta\ll1$ and $T\gg1$, depending on $\eta$.  This completes the proof
of \refeq{H4strineq}, and hence the advancement of the first inequality in
\refeq{H4}, if $\beta\ll1$, $T\gg1$ and $K_4\gg1$, under the hypotheses
that \refeq{Wbd} holds for $k\notin\cA_{n+1}$.

\begin{proof}[Proof of \refeq{Wbd}]
Given $k\notin\cA_{n+1}$, let $\mu=\mu(k)=\max\{l\in\mN:k\in\cA_l\}$.  For
$l\leq\mu$, $f_l$ is in the domain of (H3), while for $\mu<l\leq n$, $f_l$
is in the domain of (H4).  We separate the sum over $l$ in \refeq{W-def}
into two parts, corresponding respectively to $l\leq\mu$ and $\mu<l\leq n$,
yielding $W_{n+1}(k)=W_{n+1}^\leq(k)+W_{n+1}^>(k)$, where
\begin{align}
&|W_{n+1}^\leq(k)|\leq\sum_{m=n+2-\mu}^{n+1}\frac{\vep^2\Cg\beta}{(1+m\vep)
 ^{d/2}}\sum_{l=n+2-m}^\mu|f_{l-1}(k)-f_l(k)|,\lbeq{W3}\\
&|W_{n+1}^>(k)|\leq\sum_{m=2}^{n+1}\frac{\vep^2\Cg\beta}{(1+m\vep)^{d/2}}
 \sum_{l=\mu\vee(n+1-m)+1}^n|f_{l-1}(k)-f_l(k)|.\lbeq{W4}
\end{align}
By (H4) and Lemma~\ref{lem:conv} with $a=b=\frac{d}2$, we easily obtain
\begin{align}
|W_{n+1}^>(k)|\leq\sum_{m=2}^{n+1}\frac{\vep^2\Cg\beta}{(1+m\vep)^{d/2}}
 \sum_{l=n+2-m}^n\frac{\vep K_5\,a(k)^{-1-\rho}}{(1+l\vep)^{d/2}}\leq\vep
 \frac{c\Cg K_5\beta\,a(k)^{-1-\rho}}{[1+(n+1)\vep]^{d/2}}.\lbeq{W4bd}
\end{align}

It remains to consider $|W_{n+1}^\leq(k)|$.  By \refeq{H3}, \refeq{vm-bd}
and Lemma~\ref{lem:fabs}, we have
\begin{align}\lbeq{fldiff}
|f_{l-1}(k)-f_l(k)|&=|f_{l-1}(k)|\,\Big|1-\big[1-\vep v_la(k)+\vep[r_l(k)
 -r_l(0)]+\vep r_l(0)\big]\Big|\nn\\
&\leq ce^{-(l-1)\vep qa(k)}\,\vep\bigg[a(k)+\frac{K_3\beta}{(1+l\vep)^{(d
 -2)/2}}\bigg],
\end{align}
where $q=1-O(\beta)$.  We fix a small $r>0$ and separate the sum over $m$
in \refeq{W3} into $\sum_{m>r(n+1)}$ and $\sum_{m\leq r(n+1)}$ (the latter
sum may be empty depending on $\mu$).  The contribution due to the former
sum is bounded by
\begin{align}\lbeq{W31bd}
&\frac{\vep^2\Cg\beta}{[1+(n+1)\vep]^{d/2}}\sum_{m=r(n+1)+1}^{n+1}~\sum_{l
 =n+2-m}^\mu ce^{-(l-1)\vep qa(k)}\,\vep\bigg[a(k)+\frac{K_3\beta}{(1+l\vep)
 ^{(d-2)/2}}\bigg]\nn\\
&\qquad\leq\frac{\vep^2\Cg\beta}{[1+(n+1)\vep]^{d/2}}\sum_{m=r(n+1)+1}^{n+1}
 ce^{-(n+1-m)\vep qa(k)}(1+K_3\beta)\nn\\
&\qquad\leq\vep\frac{c\Cg(1+K_3\beta)\beta a(k)^{-1}}{[1+(n+1)\vep]^{d/2}}
 \leq\vep\frac{c\Cg(1+K_3\beta)\beta a(k)^{-1-\rho}}{[1+(n+1)\vep]^{d/2}}.
\end{align}

To investigate the contribution from $\sum_{m\leq r(n+1)}$, we use the
inequality
\begin{align}
e^{-(l-1)\vep qa(k)}\bigg[a(k)+\frac{K_3\beta}{(1+l\vep)^{(d-2)/2}}\bigg]
 \leq\frac{c(1+K_3\beta)a(k)^{-1-\rho}}{(1+l\vep)^{d/2}},\lbeq{fldiffbd}
\end{align}
which we will prove below.  Assuming this inequality and using
Lemma~\ref{lem:conv} with $a=b=\frac{d}2$, we obtain that the contribution
from $\sum_{m\leq r(n+1)}$ is bounded by
\begin{align}\lbeq{W32bd}
\sum_{m=n+2-\mu}^{r(n+1)}\frac{\vep^2\Cg\beta}{(1+m\vep)^{d/2}}\sum_{l=n+2-
 m}^\mu\vep\frac{c(1+K_3\beta)a(k)^{-1-\rho}}{(1+l\vep)^{d/2}}\leq\vep\frac
 {c\Cg(1+K_3\beta)\beta a(k)^{-1-\rho}}{[1+(n+1)\vep]^{d/2}}.
\end{align}
This, together with \refeq{W4bd} and \refeq{W31bd}, completes the proof
of \refeq{Wbd} for $k\notin\cA_{n+1}$.

It remains to prove \refeq{fldiffbd}.  First, we note that, by
$m\leq r(n+1)$ and $n+2-m\leq l\leq\mu\leq n$, as well as
$a(k)>\gamma\frac{\log[2+(\mu+1)\vep]}{1+(\mu+1)\vep}$, we have
\begin{align}
(1+l\vep)^{-(d-2)/2}\leq[1+(1-r)(n+1)\vep]^{-(d-2)/2}\leq\frac
 {(1-r)^{-(d-2)/2}}{[1+(\mu+1)\vep]^{1+(d-4)/2}}\leq c\,a(k).
\end{align}
Therefore, the left-hand side of \refeq{fldiffbd} is bounded by
$c(1+K_3\beta)a(k)\,e^{-(l-1)\vep qa(k)}$.  Similarly to \refeq{H3toH4},
we have
\begin{align}\lbeq{fldiffbde}
e^{-(l-1)\vep qa(k)}\leq(2+l\vep)^{-\frac{(l-1)\vep q\gamma}{1+(\mu+1)\vep}
 \frac{\log[2+(\mu+1)\vep]}{\log(2+l\vep)}q\gamma}\leq(1+l\vep)^{-q'\gamma},
\end{align}
where $q'=\frac{(1-r)(n+1)\vep}{1+(n+1)\vep}q$.  To bound \refeq{fldiffbde},
we fix $T\gg1$ and consider the case in which $n\vep\leq T$ and the other
case separately.  When $n\vep\leq T$, since $l\leq n$, $a(k)\leq2$ and
$2+\rho>0$, \refeq{fldiffbde} is bounded as
\begin{align}
\frac{(1+l\vep)^{d/2-q'\gamma}}{(1+l\vep)^{d/2}}\leq\frac{(1+T)^{d/2}}
 {(1+l\vep)^{d/2}}\leq\frac{2^{2+\rho}(1+T)^{d/2}}{(1+l\vep)^{d/2}}\,
 a(k)^{-2-\rho}.
\end{align}
When $n\vep>T$, since $a(k)\leq\gamma\frac{\log(2+l\vep)}{1+l\vep}$,
\refeq{fldiffbde} is bounded as
\begin{gather}
\frac1{(1+l\vep)^{d/2}}\frac{(1+l\vep)^{2+\rho}}{(1+l\vep)^{q'\gamma-
 [\frac{d}2-(2+\rho)]}}\leq\frac{c}{(1+l\vep)^{d/2}}\,a(k)^{-2-\rho},
\end{gather}
where we use $q'\gamma>\frac{d}2-(2+\rho)$ for $\beta\ll1$, $T\gg1$ and
$r\ll1$.  This completes the proof.
\end{proof}

Finally, we summarize the relations among the constants $K_1, \ldots, K_5$
that have been necessary in advancing the induction hypotheses.
We have taken $\beta\ll1$ and  have chosen the constants
$K_1,\dots,K_5$ such that
\begin{align}
K_1>\Cg,\qquad K_2>2\Cg,\qquad K_3\gg K_1,\qquad K_5\gg K_4,
\end{align}
where, as stated below \refeq{g2bd}, $\Cg$ depends only on $K_4$ (when
$\beta\ll1$).  This gives \refeq{Krel}.

\subsection{Advancement below and at four dimensions}\label{ss:adv-low}
The proofs of Lemmas~\ref{lem:I}--\ref{lem:fD2} and the advancement of
the induction hypotheses for $d\leq4$ remain almost unchanged,
except for the factors $\beta_{\sT}$ and $\hat\beta_{\sT}$ in
\refeq{PN-bd-lower}--\refeq{tPNn-bd-lower} and
\refeq{In-def-lower}--\refeq{H3-lower}.  We simply replace $\beta$ by
$\hat\beta_{\sT}$ in the proofs of Lemmas~\ref{lem:I}--\ref{lem:fdiff},
and by $\beta_{\sT}$ in the proof of Lemma~\ref{lem:fD2} (we also replace
$d/2$ by $2+\omega$ in the proof of Lemma~\ref{lem:I}).  In the advancement
of the first inequality in \refeq{H12-lower}, we use \refeq{lambndiff}
together with \refeq{g1bd} and \refeq{g2bd} with $\beta$ replaced by
$\beta_{\sT}$.  Since $n\vep\leq T\log T$, we obtain
\begin{align}
|\lamb_{n+1}-\lamb_n|\leq\frac{\vep\Cg\beta_{\sT}}{[1+(n+1)\vep]^{d/2}}
 +\frac{\vep K_1\hat\beta_{\sT}}{(1+n\vep)^{2+\omega}}\sum_{m=2}^n\frac
 {\vep\Cg\beta_{\sT}}{(1+m\vep)^{(d-2)/2}}\leq\vep\frac{c\Cg(1+K_1\hat
 \beta_{\sT})\hat\beta_{\sT}}{[1+(n+1)\vep]^{2+\omega}},
\end{align}
where we use $\mu\in(0,\alpha-\omega)$.  Similarly, we can advance
the second inequality in \refeq{H12-lower} and the first inequality
in \refeq{H3-lower}.

We need a little more care in advancing the second inequality in
\refeq{H3-lower} and the inequalities in \refeq{H4}.  The second
inequality in \refeq{H3-lower} is rewritten as in \refeq{rn-rn0},
and each term is bounded as in \refeq{rn-rn0-1stbd}, \refeq{rn-rn0-2ndbd}
and \refeq{rn-rn0-3rdbd} when $d>4$.  We can follow the same line when
$d\leq4$, except that, e.g., the factor $\beta^2$ in \refeq{2ndbd-pre1stbd}
is replaced by $\beta_{\sT}\hat\beta_{\sT}$, and that we use $\beta_{\sT}$
to control the convolution in \refeq{2ndbd-pre1stbd}, where the power of
$1+l\vep$ is replaced by $1+\omega$.  If we have only one factor $\beta$
as in \refeq{3rdbd-pre1stbd}, then we use $q=1+O(\hat\beta_{\sT})$ with
$\beta_1=L_1^{-d}\ll1$, as well as $\gamma+\delta<\omega$ and
$\mu<\alpha-\omega$, to obtain
\begin{align}
&c\Cg\beta_{\sT}a(k)\frac{\log[2+(n+1)\vep]}{[1+(n+1)\vep]^{1-\gamma q}}
 [1+(n+1)\vep]^{\frac{6-d}2}\nn\\
&\quad\leq\frac{c\Cg\hat\beta_{\sT}a(k)}{[1+(n+1)\vep]^\delta}T^{-bd+\mu}
 [1+(n+1)\vep]^{\frac{4-d}2+\gamma q+\delta}\log[2+(n+1)\vep]\nn\\
&\quad\leq\frac{c\Cg\hat\beta_{\sT}a(k)}{[1+(n+1)\vep]^\delta}T^{-(\alpha
 -\mu-\omega)}\leq\frac{c\Cg\hat\beta_{\sT}a(k)}{[1+(n+1)\vep]^\delta}.
\end{align}

A similar argument applies to the advancement of the inequalities in
\refeq{H4}.  However, since $-\rho>\frac{4-d}2\geq0$ (cf.,
\refeq{drgd-lower}), we cannot use the trivial inequality $a(k)\leq2$ to
obtain, e.g., the low-dimensional version of \refeq{W31bd}.  To overcome
this difficulty, we use the factor $\beta_{\sT}$ in the bound on $g_m(k)$
and $a(k)\leq\gamma\frac{\log[2+(n+1)\vep]}{1+(n+1)\vep}$ in \refeq{W31bd},
as well as $\mu<\alpha-\omega=bd+\frac{d-4}2-\omega$ and
$\frac{d}2-(2+\rho)<\gamma<\omega$, to obtain
\begin{align}
\beta_{\sT}\leq\beta_1T^{-bd}\bigg[\frac{a(k)}{\gamma\frac{\log[2+(n+1)
 \vep]}{1+(n+1)\vep}}\bigg]^{-\rho}\leq c\hat\beta_{\sT}T^{\alpha-\omega-
 bd-\rho}a(k)^{-\rho}\leq c\hat\beta_{\sT}T^{-(\omega-\gamma)}a(k)^{-\rho}
 \leq c\hat\beta_{\sT}a(k)^{-\rho}.
\end{align}
This completes the advancement of the induction hypotheses for $d\leq4$.
\qed

\section*{Acknowledgements}
The work of RvdH and AS was supported in part by
the Netherlands Organisation for Scientific Research (NWO).
The work of RvdH was carried out in part at Delft University
of Technology, the Netherlands. The work of AS was supported
in part by NSERC of Canada and was carried out in part at the
University of British Columbia, Vancouver,
Canada. This project was initiated during an extensive
visit of RvdH to the University of British Columbia, Vancouver,
Canada. We thank Gordon Slade and Ed Perkins for stimulating
discussions during various stages of the project.


\begin{thebibliography}{99}

\bibitem{an84}M. Aizenman and C. M. Newman.
\newblock Tree graph inequalities and critical behavior in percolation models.
\newblock {\it J. Statist. Phys.} {\bf 36} (1984): 107--143.

\bibitem{ba91}D. J. Barsky and M. Aizenman.
\newblock Percolation critical exponents under the triangle condition.
\newblock {\it Ann. Probab.} {\bf 19} (1991): 1520--1536.

\bibitem{bw98}D. J. Barsky and C. C. Wu.
\newblock Critical exponents for the contact process under the triangle condition.
\newblock {\it J. Statist. Phys.} {\bf 91} (1998): 95--124.


\bibitem{bg91}C. Bezuidenhout and G. Grimmett.
\newblock Exponential decay for subcritical contact and percolation processes.
\newblock {\it Ann. Probab.} {\bf 19} (1991): 984--1009.

\bibitem{BR01}
E.~Bolthausen and C.~Ritzmann.
\newblock Strong pointwise estimates for the weakly self-avoiding walk.
\newblock To appear in {\em Ann. Probab.}

\bibitem{bcr99}
C.~Borgs, J. T. Chayes and D. Randall.
\newblock The van den Berg-Kesten-Reimer inequality.
\newblock {\it Perplexing Problems in Probability: Festschrift in honor of
Harry Kesten} (eds., M. Bramson and R. Durrett). Birkh\"auser (1999): 159-173.


\bibitem{BS85}
D.C. Brydges and T.~Spencer.
\newblock Self-avoiding walk in 5 or more dimensions.
\newblock {\em Commun. Math. Phys.}, {\bf 97} (1985):125--148.



\bibitem{dp99}R. Durrett and E. Perkins.
\newblock Rescaled contact processes converge to super-Brownian motion in two
or more dimensions.
\newblock {\it Probab. Th. Rel. Fields} {\bf 114} (1999): 309--399.

\bibitem{g99}G. Grimmett.
\newblock {\it Percolation}.  Springer, Berlin (1999).

\bibitem{gh01}G. Grimmett and P. Hiemer.
\newblock Directed percolation and random walk.
\newblock {\it In and Out of Equilibrium} (ed., V. Sidoravicius).
Birkh\"auser (2002): 273-297.

\bibitem{HS90a}
T.~Hara and G.~Slade.
\newblock Mean-field critical behaviour for percolation in high dimensions.
\newblock {\em Commun. Math. Phys.}, {\bf 128} (1990):333--391.

\bibitem{HS92b}
T.~Hara and G.~Slade.
\newblock Self-avoiding walk in five or more dimensions. I.  The critical behaviour.
\newblock {\em Commun.\ Math.\ Phys.}, {\bf 147} (1992):101--136.

\bibitem{HS00a}T. Hara and G. Slade.
\newblock The scaling limit of the incipient infinite cluster in high-dimensional
percolation. I. Critical exponents.
\newblock {\it J. Statist. Phys.} {\bf 99} (2000): 1075--1168.

\bibitem{HS00b}T. Hara and G. Slade.
\newblock The scaling limit of the incipient infinite cluster in high-dimensional
percolation. II. Integrated super-Brownian excursion.
\newblock {\it J. Math. Phys.} {\bf 41} (2000): 1244--1293.

\bibitem{hhs98}R. van der Hofstad, F. den Hollander and G. Slade.
\newblock A new inductive approach to the lace expansion for self-avoiding walks.
\newblock {\it Probab. Th. Rel. Fields} {\bf 111} (1998): 253--286.

\bibitem{hhs01}R. van der Hofstad, F. den Hollander and G. Slade.
\newblock Construction of the incipient infinite cluster for the spread-out
oriented percolation above 4+1 dimensions.
\newblock {\it Commun. Math. Phys.} {\bf 231} (2002): 435--461.

%

\bibitem{hsa04b}R. van der Hofstad and A. Sakai.
\newblock Critical points for spread-out self-avoiding walk, percolation and
the contact process above the upper critical dimensions.
\newblock Preprint (2004).

\bibitem{hsa04}R. van der Hofstad and A. Sakai.
\newblock Convergence of the critical finite-range contact process to
super-Brownian motion above 4 spatial dimensions.
\newblock In preparation.

\bibitem{hs02}R. van der Hofstad and G. Slade.
\newblock A generalised inductive approach to the lace expansion.
\newblock {\it Probab. Th. Rel. Fields} {\bf 122} (2002): 389--430.

\bibitem{hs01}R. van der Hofstad and G. Slade.
\newblock Convergence of critical oriented percolation to super-Brownian motion
above 4+1 dimensions.
\newblock {\it Ann. Inst. H. Poincar\'e Probab. Statist.} {\bf 39} (2003): 413--485.

\bibitem{hs03}R. van der Hofstad and G. Slade.
\newblock The lace expansion on a tree with application to networks of self-avoiding
walks.
\newblock {\it Adv. Appl. Math.} {\bf 30} (2003): 471--528.

\bibitem{Ligg99}T. Liggett.
\newblock {\em Stochastic Interacting Systems: Contact, Voter and Exclusion Processes}.
\newblock Springer, Berlin (1999).

\bibitem{ms93}N. Madras and G. Slade.
\newblock {\it The Self-Avoiding Walk}.
\newblock Birkh\"auser, Boston (1993).

\bibitem{ny93}B. G. Nguyen and W.-S. Yang.
\newblock Triangle condition for oriented percolation in high dimensions.
\newblock {\it Ann. Probab.} {\bf 21} (1993): 1809--1844.

\bibitem{ny95}B. G. Nguyen and W.-S. Yang.
\newblock Gaussian limit for critical oriented percolation in high dimensions.
\newblock {\it J. Statist. Phys.} {\bf 78} (1995): 841--876.

\bibitem{s00}A. Sakai.
\newblock Analyses of the critical behavior for the contact process based on
a percolation structure.
\newblock {\it Ph.D. thesis} (2000).

\bibitem{s01}A. Sakai.
\newblock Mean-field critical behavior for the contact process.
\newblock {\it J. Statist. Phys.} {\bf 104} (2001): 111--143.

\bibitem{s02}A. Sakai.
\newblock Hyperscaling inequalities for the contact process and oriented
percolation.
\newblock {\it J. Statist. Phys.} {\bf 106} (2002): 201--211.

%

\bibitem{Scho98}R. Schonmann.
\newblock The triangle condition for contact processes on homogeneous trees.
\newblock {\it J. Statist. Phys.} {\bf 90} (1998): 1429--1440.

\bibitem{Slad87}
G.~Slade.
\newblock The diffusion of self-avoiding random walk in high
dimensions.
\newblock {\em Commun.\ Math.\ Phys.}, {\bf 110} (1987): 661-683.

\bibitem{Slad88}
G.~Slade.
\newblock Convergence of self-avoiding random walk
to {B}rownian motion in high dimensions.
\newblock  {\em J. Phys. A:  Math. Gen.}, {\bf 21} (1988):L417-L420.

\bibitem{Slad91} G.~Slade.
\newblock The lace expansion and the upper critical dimension
for percolation.
\newblock {\it Lectures in Applied Mathematics} {\bf 27} (1991):53--63.

\bibitem{Slad04}
G.~Slade.
\newblock The lace expansion and its applications.
\newblock Saint-Flour lecture notes.  Preprint (2004).

\bibitem{wu95}C. C. Wu.
\newblock The contact process on a tree: Behavior near the first phase transition.
\newblock {\it Stochastic Process. Appl.} {\bf 57} (1995): 99--112.

\end{thebibliography}
\end{document}